\documentclass[12pt]{article}
\usepackage{amsthm,amsmath,amsfonts,epsfig,amssymb,latexsym, bbm}
\usepackage{mathrsfs} 
\usepackage[T1]{fontenc}
\usepackage{graphicx}
\usepackage{enumerate}
\usepackage{xy,xypic}
\usepackage{graphics}
\usepackage{footnote}
\usepackage{scalefnt}
\usepackage{tikz}
\usepackage{color}
\usetikzlibrary{shapes,arrows,fit,calc,positioning}
\usetikzlibrary{matrix}
\xyoption{all}
\usepackage{imakeidx}
\makeindex

\title{Chow-Witt Theory in the Arrowtic Paradigm} %Attempt to Transfer Maps} %\\ Example of Mohan} 
\author{
 Satya Mandal \\ {\small University of Kansas,  Lawrence, Kansas 66045, USA,  mandal@ku.edu}\\
\TCP{\bf ArrowticKTheory.com}\\
 }

\begin{document}

\pagenumbering{roman}
\setcounter{page}{0}

\renewcommand{\baselinestretch}{1.255}
\setlength{\parskip}{1ex plus0.5ex}
\date{13 August 2026} %27 July 2026}%{14 July 2026} %5 July 2026\\13 June 2026 (combining), 28 May 26} % 5 May 26} 
\newcommand{\iso}{\stackrel{\sim}{\longrightarrow}}
\newcommand{\sur}{\twoheadrightarrow}

\newcommand{\eop}{\hfill{\TCP{\rule{2mm}{2mm}}}}
\newcommand{\pf}{\noindent{\bf Proof.~}}
\newcommand{\outlinePf}{\noindent{\bf Outline of the Proof.~}}

\newcommand{\PD}{\dim_{{\SA}}}
\newcommand{\PDV}{\dim_{{\SV}}}

\newcommand{\ra}{\rightarrow}
\newcommand{\lra}{\longrightarrow}
\newcommand{\hra}{\hookrightarrow} 
\newcommand{\Lra}{\Longrightarrow}
\newcommand{\Lla}{\Longleftarrow}
\newcommand{\Llra}{\Longleftrightarrow}
\newcommand{\llra}{\longleftrightarrow} 
\newcommand{\ora}{\overrightarrow}

\newcommand{\pic}{The proof is complete.~}
\newcommand{\Sp}{\mathrm{Sp}}
\newcommand{\BiSp}{\mathrm{BiSp}}
\newcommand{\colim}{\mathrm{colim}}
\newcommand{\codim}{\mathrm{co}\dim}
\newcommand{\coh}{\mathrm{Coh}}
\newcommand{\qcoh}{\mathrm{QCoh}}

\newcommand{\bul}{\bullet}

\newcommand{\bE}{\begin{enumerate}}
\newcommand{\eE}{\end{enumerate}}

\newtheorem{theorem}{Theorem}[section]
\newtheorem{proposition}[theorem]{Proposition}
\newtheorem{lemma}[theorem]{Lemma}
\newtheorem{definition}[theorem]{Definition}
\newtheorem{corollary}[theorem]{Corollary}
\newtheorem{construction}[theorem]{Construction}
\newtheorem{notation}[theorem]{Notation}
\newtheorem{notations}[theorem]{Notations}
\newtheorem{remark}[theorem]{Remark}
\newtheorem{question}[theorem]{Question}
\newtheorem{example}[theorem]{Example} 
\newtheorem{examples}[theorem]{Examples} 
\newtheorem{exercise}[theorem]{Exercise} 
\newtheorem{setup}[theorem]{Setup} 
\newtheorem{openProblem}[theorem]{Open Problem} 
\newtheorem{clarification}[theorem]{Clarification}

\newtheorem{problem}[theorem]{Problem} 
\newtheorem{conjecture}[theorem]{Conjecture} 

\newcommand{\bD}{\begin{definition}}
\newcommand{\eD}{\end{definition}}
\newcommand{\bP}{\begin{proposition}}
\newcommand{\eP}{\end{proposition}}
\newcommand{\bL}{\begin{lemma}}
\newcommand{\eL}{\end{lemma}}
\newcommand{\bT}{\begin{theorem}}
\newcommand{\eT}{\end{theorem}}
\newcommand{\bC}{\begin{corollary}}
\newcommand{\eC}{\end{corollary}} 
%%%%%%%%%%%
\newcommand{\bA}{\begin{application}} 
\newcommand{\eA}{\end{application}} 

\newcommand{\bN}{\begin{notation}} 
\newcommand{\eN}{\end{notation}} 

\newcommand{\bR}{\begin{remark}} 
\newcommand{\eR}{\end{remark}} 

\newcommand{\bExr}{\begin{exercise}} 
\newcommand{\eExr}{\end{exercise}} 

\newcommand{\bExm}{\begin{example}} 
\newcommand{\eExm}{\end{example}}
%%%
\newcommand{\TCP}{\textcolor{purple}}
\newcommand{\TCM}{\textcolor{magenta}}
\newcommand{\TCR}{\textcolor{red}}
\newcommand{\TCB}{\textcolor{blue}}
\newcommand{\TCG}{\textcolor{green}}
%%%%%%%%%%%%%%%%SSSSSSSSSSSSSSSSSSSSSSSSSSSSSSSSSSSSS I've added these
\def\spec#1{\mathrm{Spec}\left(#1\right)}
\def\proj#1{\mathrm{Proj}\left(#1\right)}
\def\supp#1{\mathrm{Supp}\left(#1\right)}
\def\Sym#1{{\CS}\mathrm{ym}\left(#1\right)} %Symm algebra
\def\norm#1{\parallel #1\parallel}
\def\LRf#1{\left( #1\right)}
\def\LRs#1{\left\{ #1\right\}}
\def\LRt#1{\left[ #1\right]}
\def\LR|#1{\left| #1\right|}

\def\LBrace#1{\left\{\begin{array} #1\end{array}\right.}
\def\matrix#1{\left(\begin{array} #1\end{array}\right)}

\def\Hom#1{\underline{\mathrm{Hom}}\left(#1\right)}
\def\Obj#1{\underline{\mathrm{Obj}}\left(#1\right)}

\def\rank#1{\mathrm{rank}\left(#1\right)}

%%%%%
\def\Cat{\underline{\bf Cat}}
\def\eCat{\underline{\bf eCat}}
\def\dgCat{\underline{\bf dgCat}}
\def\SSet{\underline{\bf SSet}}
\def\Top{\underline{\bf Top}}

\def\0{\mathfrak {0}}
\def\1{\mathfrak {1}}
\def\2{\mathfrak {2}}
\def\3{\mathfrak {3}}
\def\4{\mathfrak {4}}
\def\5{\mathfrak {5}}
\def\6{\mathfrak {6}}
\def\7{\mathfrak {7}}
\def\8{\mathfrak {8}}
\def\9{\mathfrak {9}}

\def\a{\mathfrak {a}}
\def\b{\mathfrak {b}}
\def\c{\mathfrak {c}}

\def\d{\mathfrak {d}}
\def\e{\mathfrak {e}}
\def\f{\mathfrak {f}}
\def\g{\mathfrak {g}}
\def\h{\mathfrak {h}}
\def\i{\mathfrak {i}}
\def\j{\mathfrak {j}}
\def\k{\mathfrak {k}}
\def\l{\mathfrak {l}}
\def\m{\mathfrak {m}}
\def\n{\mathfrak {n}}
\def\p{\mathfrak {p}}
\def\q{\mathfrak {q}}
\def\r{\mathfrak {r}}
\def\s{\mathfrak {s}}
\def\t{\mathfrak {t}}
\def\u{\mathfrak {u}}
\def\v{\mathfrak {v}}

\def\w{\mathfrak {w}}

\def\x{\mathfrak {x}}
\def\y{\mathfrak {y}}
\def\z{\mathfrak {z}}

\def\A{\mathfrak {A}}
\def\B{\mathfrak {B}}
\def\C{\mathfrak {C}}
\def\D{\mathfrak {D}}
\def\E{\mathfrak {E}}

\def\F{\mathfrak {F}}

\def\G{\mathfrak {G}}
\def\H{\mathfrak {H}}
\def\I{\mathfrak {I}}
\def\J{\mathfrak {J}}
\def\K{\mathfrak {K}}
\def\L{\mathfrak {L}}
\def\M{\mathfrak {M}}
\def\N{\mathfrak {N}}
\def\O{\mathfrak {O}}
\def\P{\mathfrak {P}}

\def\Q{\mathfrak {Q}}

\def\R{\mathfrak {R}}

\def\Sf{\mathfrak {S}} %Do not use
\def\T{\mathfrak {T}}
\def\U{\mathfrak {U}}
\def\V{\mathfrak {V}}
\def\W{\mathfrak {W}}

%%%%%
\def\CA{\mathcal {A}}
\def\CB{\mathcal {B}}
\def\CP{\mathcal {P}}
\def\CC{\mathcal {C}}
\def\CD{\mathcal {D}}
\def\CE{\mathcal {E}}
\def\CF{\mathcal {F}}
\def\CG{\mathcal {G}}
\def\CH{\mathcal{H}}
\def\CI{\mathcal{I}}
\def\CJ{\mathcal{J}}
\def\CK{\mathcal{K}}
\def\CL{\mathcal{L}}
\def\CM{\mathcal{M}}
\def\CN{\mathcal{N}}
\def\CO{\mathcal{O}}
\def\CP{\mathcal{P}}
\def\CQ{\mathcal{Q}}

\def\CR{\mathcal{R}}

\def\CS{\mathcal{S}}
\def\CT{\mathcal{T}}
\def\CU{\mathcal{U}}
\def\CV{\mathcal{V}}
\def\CW{\mathcal{W}} 
\def\CX{\mathcal{X}}
\def\CY{\mathcal{Y}}
\def\CZ{\mathcal{Z}}

\newcommand{\smallcirc}[1]{\scalebox{#1}{$\circ$}}

\def\BA{\mathbb{A}}
\def\BB{\mathbb{B}}
\def\BC{\mathbb{C}}
\def\BD{\mathbb{D}}
\def\BE{\mathbb{E}}
\def\BF{\mathbb{F}}
\def\BG{\mathbb{G}}
\def\BH{\mathbb{H}}
\def\BI{\mathbb{I}}
\def\BJ{\mathbb{J}}
\def\BK{\mathbb{K}}
\def\BL{\mathbb{L}}
\def\BM{\mathbb{M}}
\def\BN{\mathbb{N}} 
\def\BO{\mathbb{O}}
\def\BP{\mathbb{P}} 
\def\BQ{\mathbb{Q}}
\def\BR{\mathbb{R}}
\def\BS{\mathbb{S}}
\def\BT{\mathbb{T}}
\def\BU{\mathbb{U}}
\def\BV{\mathbb{V}}
\def\BW{\mathbb{W}}
\def\BX{\mathbb{X}}
\def\BY{\mathbb{Y}}
\def\BZ{\mathbb {Z}}

\def\SA{\mathscr {A}} 
\def\SB{\mathscr {B}}
\def\SC{\mathscr {C}}
\def\SD{\mathscr {D}}
\def\SE{\mathscr {E}}
\def\SF{\mathscr {F}}
\def\SG{\mathscr {G}}
\def\SH{\mathscr {H}}
\def\SI{\mathscr {I}}
\def\SJ{\mathscr {J}}
\def\SK{\mathscr {K}} 
\def\SL{\mathscr {L}} 
\def\SM{\mathscr {M}}
\def\SN{\mathscr {N}}
\def\SO{\mathscr {O}}
\def\SP{\mathscr {P}}
\def\SQ{\mathscr {Q}}
\def\SR{\mathscr {R}}
\def\SS{\mathscr {S}}
\def\ST{\mathscr {T}}
\def\SU{\mathscr {U}}
\def\SV{\mathscr {V}}
\def\SW{\mathscr {W}}
\def\SX{\mathscr {X}}
\def\SY{\mathscr {Y}}
\def\SZ{\mathscr {Z}} 
 
 %\def\ST{\mathscr {T}}
 %%%
\def\bfA{\bf A}
\def\bfB{\bf B}
\def\bfC{\bf C}
\def\bfD{\bf D}
\def\bfE{\bf E}
\def\bfF{\bf F}
\def\bfG{\bf G}
\def\bfH{\bf H}
\def\bfK{\bf K}
\def\bfI{\bf I}
\def\bfJ{\bf J}
\def\bfK{\bf K}
\def\bfL{\bf L}
\def\bfM{\bf M}
\def\bfN{\bf N} 
\def\bfO{\bf O}
\def\bfP{\bf P} 
\def\bfQ{\bf Q}
\def\bfR{\bf R}
\def\bfS{\bf S}
\def\bfT{\bf T}
\def\bfU{\bf U}
\def\bfV{\bf V}
\def\bfW{\bf W}
\def\bfX{\bf X}
\def\bfY{\bf Y}
\def\bfZ{\bf Z}

 \def\bfw{\bf w}
 %%%
 %%%%%%%%%%%%%%%%%%%%%%%%%%%%%%%%%
 
\maketitle
\thispagestyle{empty}
%\vspace*{\fill}
\begin{flushright} %center}
\itshape
In memory of Prema and M. Pavaman Murthy!
\end{flushright}%center}
%\vspace*{\fill}
%\clearpage

\begin{abstract}
In this article we propose a theory of Chow-Witt groups, in the paradigm of Quillen's arrow based approach to $K$-theory. We refer to this style of arguments as arrowtic paradigm. Other than establishing the machinery needed to work in this paradigm, we establish the homotopy invariance property of the Arrowtic Chow-Witt groups. We also compute the arrowtic Chow-Witt groups of the projective spaces, over fields. 
\end{abstract}
%%%%

\tableofcontents

%%%%%%%%%%%%%%%%%%%%%%%%%%%%%%%%%%%%%%%%%%%%%%%%%%%%%%%%%%%%%%%%%%%%%%%%%%%
%%%%%%%%%%%%%%%%%%%%%%%%%%%%%%%%%%%%%%%%%%%%%%%%%%%%%%%%%%%%%%%%%%%%%%%%%%%

%%%%%%%%%%%%%%%%%%%%%%%%%%%%%%%%%%%%%%%%%%%%%%%%%%%%%%%%%%%%%%%%%%%%%%%%%%%

%

%\vspace{10mm}

%\tableofcontents
%\setcounter{chapter}{-1}
%\chapter{Introduction}
\pagenumbering{arabic}
%%%%%%%%%%%%%%%%%%%%%%%%%%%%%%%%%%%%%%%%%%%%%%%%%%%%%%%%%
%%%%%%%%%%%%%%%%%%%%%%%%%%%%%%%%%%%%%%%%%%%%%%%%%%%%%%%%%
%%%%%%%%%%%%%%%%%%%%%%%%%%%%%%%%%%%%%%%%%%%%%%%%%%%%%%%%%
%%%%%%%%%%%%%%%%%%%%%%%%%%%%%%%%%%%%%%%%%%%%%%%%%%%%%%%%%
%%%%%%%%%%%%%%%%%%%%%%%%%%%%%%%%%%%%%%%%%%%%%%%%%%%%%%%%%
%%%%%%%%%%%%%%%%%%%%%%%%%%%%%%%%%%%%%%%%%%%%%%%%%%%%%%%%%
\section{Introduction} \label{SecIntro}
In this article we propose an alternate theory of Chow-Witt groups, in the arrowtic paradigm. By the expression 
"arrowtic" we generally mean Quillen's arrow based homotopy approach to $K$-Theory. It is an alternate to
the paradigm of Milnor $K$-theory approach, which has come to have been encapsulated within the Motivic $K$-theory. While Quillen $K$-theory does not agree with Milnor $K$-theory, 
the nature of arguments used in these two paradigms remain disjoint. 
 In the Milnor $K$-theory  the definitions and arguments are based on generators and relations, while  in Quillen $K$-theory arguments are based on arrows and diagrams. 
%%%%%%%%%%%%%%%%%%%%%%%%%%%%%%%%%%%%%%%%%%%%%%%%%%%%%%%%%

The common denominators of higher ({\it inclusive of lower}) Chow groups and Chow-Witt groups are the existence of respective Gersten complexes. In a more formal setting,  expositions of such commonality  are given in 
other sources \cite[Chap 4]{M12}, \cite{R96}. For our purpose, 
let $X$ be a noetherian scheme (regular), with $\dim X=d$. For each $k\in {\BZ}$, there may be a Gersten complex $C^{{\bul}k}\LRf{X}$ of length $d+1$. Now, the cohomologies $H^p\LRf{C^{{\bul}k}\LRf{X}}$ are of our interest, in this context. For example,  the Milnor $K$-theory complex ${\bfK}M^{{\bul}k}\LRf{X}$ is defined as follows:
 \begin{equation}\label{MilGer}
 \diagram 
 0\ar[r] &  K_k\LRf{X_0} \ar[r] & \bigoplus_{x\in X^{\LRf{1}}} K_{k-1}\LRf{X_x} \ar[r] & \cdots \ar[r] &
 \bigoplus_{x\in X^{\LRf{k}}} K_0\LRf{X_x} \ar[r] & 0\\
 \enddiagram 
 \end{equation}
 where for $x\in X$ we denote $X_x=\spec{{\CO}_{X, x}}$, and subsequently, the residue field is denoted by $\kappa(x)$. 
 The Milnor complex (\ref{MilGer}) can be referred to as the Gersten complex of Milnor $K$-theory. 
 The cohomologies $H^p\LRf{{\bfK}M^{{\bul}k}\LRf{X}}=:CH^{pk}\LRf{X}$ are defined to be the higher Chow groups,
 of Milnor $K$-theory.  Similar Gersten complexes in Witt theory are also available in \cite{BW2}.
 Further on the hermitian side, the Gersten complex of the Milnor-Witt $K$-theory  
 ${\bfK}MW^{{\bul}k}\LRf{X}$ is defined as follows \cite{M12}:
 \begin{equation}\label{MWGerValo}
 \diagram 
 0\ar[r] &  K^{MW}_k\LRf{\kappa\LRf{0}} \ar[r] & \bigoplus_{x\in X^{\LRf{1}}} K^{MW}_{k-1}\LRf{\kappa\LRf{x}} \ar[r] & \cdots \ar[r] &
 \bigoplus_{x\in X^{\LRf{k}}} K^{MW}_0\LRf{\kappa\LRf{x}}\ar[dllll] \\
\cdots\ar[rr] && \bigoplus_{x\in X^{\LRf{d}}} K^{MW}_{k-d}\LRf{\kappa\LRf{x}}\ar[r]&0&\\
 \enddiagram 
 \end{equation}
 Homologies of these complexes  ${\bfK}MW^{{\bul}k}\LRf{X}$ are denoted by $H^p\LRf{X, {\bfK}_k^{MW}}$, with 
 $p=0, 1, \ldots, d$, 
 in \cite{M12}. These are also viewed as the sheaf cohomologies. 
 The key to the construction of these Gersten complexes (\ref{MWGerValo}) is the definition, due to Hopkins and Morel, of the Milnor-Witt $K$-theory 
 $K_{\bul}^{MW}\LRf{\BF}$ of fields ${\BF}$. The ${\BZ}$-graded ring $K_{\bul}^{MW}\LRf{\BF}$ is defined by generators and relations. 
 The celebrated Chow-Witt groups
 are defined as $\widetilde{CW}^k\LRf{X}:= H^k\LRf{{\bfK}MW^{{\bul}k}\LRf{X}}$. While Chow-Witt groups 
 $\widetilde{CW}^k\LRf{X}$ stole the limelight, all the cohomologies $\widetilde{CW}^{pk}\LRf{X}:= H^p\LRf{{\bfK}MW^{{\bul}k}\LRf{X}}$ are studied together in the literature. We introduce this notation $\widetilde{CW}^{pk}\LRf{X}$ for subsequent  references, for the purpose of analogies. The original definition of Chow-Witt groups is due to Barge and Morel
 \cite{BM0}, which precedes  the book  \cite{M12} and other related work on Milnor-Witt complexes. The theory was  developed more comprehensively in \cite{F8}, and Chow-Witt groups $\widetilde{CW}^{pk}\LRf{X, {\CL}}$ were defined, for any line bundle ${\CL}$ on $X$. For $0\leq k\leq d$, a complex $BM^{{\bul}k}\LRf{X, {\CL}}$ was defined, in \cite{BM0, F8}, by patching the Minor complex (\ref{MilGer}) and a subcomplex ${\bf I}^k\LRf{X, {\CL}}$ of the Witt complex \cite{BW2}. It was later  established that $\widetilde{CW}^{pk}\LRf{X, {\CL}}:=H^p\LRf{BM^{{\bul}k}\LRf{X, {\CL}}}$ \cite[pp 147]{M12}. All these belong to the Motivic paradigm.

 Same  could be done with the Gersten complexes of Quillen $K$-theory \cite[pp 252]{M23}. The   area remained more dormant than desired. Following the work of Quillen \cite{Q73, G76}, among the benchmark work 
 in this paradigm are the work of Waldhausen \cite{W83}, Thomason \cite{TT90}, Gillet-Grayson \cite{GG87}, and others. In the hermitian side, the theory was mostly developed by Schlichting in multiple articles, culminating  in \cite{S17}. The hermitian theory is usually referred to as the Grothendieck-Witt theory, ({\it  henceforth, abbreviated as} $GW$-theory). The $GW$-theory supersedes $K$-theory, 
 in the sense that $K$-theory coincides with the $GW$-theory of the corresponding hyperbolic category. 
 
 In the arrowtic paradigm there are two versions of the theory - connective and the non-connective theory.  The two of them coincide when applied to  the regular schemes and the latter is applicable for non-regular schemes. 
% We ignore $K$-theory complexes because we do not have any new results({Homotopy invariance works and new}).  
We mostly marvel the $GW$ theory and  outline the results on $K$-theory. 
 The hermitian ($GW$) theory is further enriched by incorporating  an $n$-shifted duality,
 for each $n\in {\BZ}$, which has periodicity 4. Accordingly, for each $n\in {\BZ}$, 
  corresponding to the Minor-Witt Gersten complexes  (\ref{MWGerValo}),   two versions  (connective and non-connective) of $GW$-Gersten complexes \cite[pp 474]{M23}, of 
  connective $G{\CW}$ groups and of non-connetive ${\BG}W$-groups, 
 are available in the arrowtic paradigm. %We mostly, ignore the non-connective version, and 
 For regular schemes two versions coincide. So, for a regular quasi projective scheme $X$ over an affine scheme $\spec{A}$, with 
 $\frac{1}{2}\in A$, and a line bundle ${\CL}$ on $X$, there are $GW$-Gersten complexes
 $$
 C^{\LRt{n}{\bul}k}\LRf{X, {\CL}}:=
 \LRs{
 C^{\LRt{n}pk}\LRf{X, {\CL}}= \bigoplus_{x\in X^{\LRf{p}}}{G}{\CW}_{k-p}^{\LRt{n+k-p}}\LRf{C{\BM}^{p}\LRf{X_x, {\CL}_x}}: p\in {\BZ}
 }\quad \forall n, k\in {\BZ}
 $$
 Here $n$ corresponds to the shift and $k$ indicates the $k$-complex (see Definition \ref{notaSptrCx}). 
 It has been customary to present such Gersten complexes as converging spectral sequences, which we avoid.
 However, the notation $C^{\LRt{n}pk}\LRf{X, {\CL}}$ was designed with spectral sequence notations in mind ({\it shift} $n$ {\it being fixed}).
 Taking cohomologies, we define 
 (Arrowtic) Chiow-Witt groups 
 $$
 \ora{CW}^{\LRt{n}pk}\LRf{X, {\CL}}:= H^p\LRf{ C^{\LRt{n}pk}\LRf{X, {\CL}}}
 $$
 At the outset, we purposefully use the notation $\ora{CW}^{\LRt{n}pk}\LRf{X, {\CL}}$, to distinguish them from motivic Chow-Witt groups  $\widetilde{CW}^{pk}\LRf{X, {\CL}}$. 
In the motivic paradigm, the aspect of the  $n$-shifted duality is fully ignored in the literature. 
The shift zero case $\ora{CW}^{\LRt{0}pk}\LRf{X, {\CL}}$  would correspond to $\widetilde{CW}^{pk}\LRf{X, {\CL}}$.
Any speculation about the nature of correspondence, between the two sets of Chow-Witt groups,  may be premature, at this time. However, for a local 
regular ring $\LRf{R, \m}$, with an infinite residue field, a natural map 
$$
\diagram 
\bigoplus_{p\in {\BZ}}K^{MW}_{p}\LRf{{R}} \ar[r] & \bigoplus_{p\in {\BZ}}G{\CW}_{p}^{\LRt{p}}\LRf{{R}}
\enddiagram 
$$
of graded rings, was established in \cite{AF17, SS25}, which are isomorphic at low degrees. 

We prove homotopy invariance of the arrowtic Chow-Witt groups (\ref{VectorBunCW}), as follows. 
%%%%%%%%%%%%%%%%%%%%%%%%%%%%%%%%%%%%%%%%%%%%%%%%%%
\bT\label{homoThmIntro}{\rm 
Suppose $X$ is a quasi projective scheme over a noetherian affine scheme $\spec{A_0}$. Assume   $X$ is regular,
with $\frac{1}{2}\in A_0$.  Let 
$\pi: Y \lra X$ be a vector bundle. %, meaning, $Y=\spec{\Sym{{\SE}}}$, where ${\SE}\in {\SV}\LRf{X}$. 
Let ${\CL}\in {\SV}\LRf{X}$ be a line bundle. All the dualities are induced by 
the duality $P \mapsto {\CH}om\LRf{-, {\CL}}$ on ${\SV}\LRf{X}$.
 Then $\forall p, k, n$, the pullback map
$$
\diagram 
\ora{CW}^{\LRt{n}pk} \LRf{X,  {\CL}} \ar[rr]^{\pi^*\quad}_{\sim\quad} && \ora{CW}^{\LRt{n}pk}\LRf{Y,  \pi^*{\CL}} \\
\enddiagram 
$$
is  an isomorphism.  
}
\eT

The other main result in this article is the description of the Chow-Witt groups of the projective spaces (\ref{mainPSpace}), as follows. 
\bT\label{ThmIntoPd}{\rm 
Let ${\BF}$ be a field, with $\frac{1}{2}\in {\BF}$. For $m, k\in {\BZ}$, define the map $\b^{\LRt{m}}_k$, by the commutative diagram:  
\begin{equation}\label{betakmdefDia}
\diagram 
G{\CW}^{\LRt{m}}_k\LRf{{\BF}} \ar@{^(->}[r]^{\pi} \ar@{-->}[d]_{\b^{\LRt{m}}_k}& G{\CW}^{\LRt{m}}_k\LRf{{\BF}\LRf{T}} \ar[d]^{\partial}\\
G{\CW}^{\LRt{m-1}}_{k-1}\LRf{{\BF}}\ar[r]_{\D\qquad\quad}^{\sim\qquad\quad}& G{\CW}^{\LRt{m-1}}_{k-1}\LRf{C{\BM}^1{\BF}\LRt{T}_{\LRf{T}}} \\
\enddiagram
\quad 
\LBrace{{l}
{\rm Here,}~\pi~{\rm is~the~pullback~map}\\
\partial ~{\rm is~the ~residue ~map~(\ref{defResidue})}\\ %XXXfrom~localization}\\
\D~{\rm is~the~d\'{e}vissage~isomorphism}\\
C{\BM}^1{\BF}\LRt{T}_{\LRf{T}}~{\rm is~the ~category}\\
{\rm of~finite~length ~modules}\\
}
\end{equation}
({\it So, $\b^{\LRt{m}}_k$ is the restriction of the residue map $\partial$, up to the d\'{e}vissage isomorphism} $\D$.)
Let ${\BP}^d=\proj{{\BF}\LRt{X_0, X_1, \ldots, X_d}}$ be the projective space. 
Let $k, n\in {\BZ}$.
Then 
\begin{equation}\label{GanesherAshirvadIntro}
\begin{array}{l}
\ora{CW}^{\LRt{n}pk}\LRf{{\BP}^d,  {\CO}\LRf{1}} =
\LBrace{{ll}
\ker\LRf{\b_{k-2r}^{\LRt{n+k-2r}}} & 0\leq p=2r\leq d-1\\
co\ker\LRf{\b_{k-2r}^{\LRt{n+k-2r}}} & 1\leq p=2r+1\leq d\\
G{\CW}^{\LRt{n+k-d}}_{k-d}\LRf{{\BF}} & p=d~even\\
%\ora{CW}^{\LRt{n}(p-1)(k-1)}\LRf{{\BP}^{d-1}, 0} &p=2, \ldots, d\\
0 & Otherwise\\
} \\
%$$
%and 
%$$
\ora{CW}^{\LRt{n}pk}\LRf{{\BP}^d,  {\CO}\LRf{0}} =
\LBrace{{ll}
G{\CW}^{\LRt{n+k}}_k\LRf{{\BF}} & p=0\\
\ora{CW}^{\LRt{n}(p-1)(k-1)}\LRf{{\BP}^{d-1}, {\CO}\LRf{1}} & p=1, \ldots, d\\
0 & Otherwise\\
} \\
\end{array}
\end{equation} 
}
\eT
%%%%%%%%%%%%%%%%%%%%%%%%%%%%%%%%%%%%%%%%%%%%%%%%%%%%%%%%%
%\TCP{XXX}
%
%}
%\eT
%%%%%%%%%%%%%%%%%%%%%%%%%%%%%%%%%%%%%%%%%%%%%%%%%%%%%%%%%
%%%%%%%%%%%%%%%%%%%%%%%%%%%%%%%%%%%%%%%%%%%%%%%%%%%%%%%%%
%%%%%%%%%%%%%%%%%%%%%%%%%%%%%%%%%%%%%%%%%%%%%%%%%%%%%%%%%
%%%%%%%%%%%%%%%%%%%%%%%%%%%%%%%%%%%%%%%%%%%%%%%%%%%%%%%%%
In fact, we provide (\S \ref{secOtheDes}) other descriptions of (\ref{GanesherAshirvadIntro}), by establishing  equalities (\ref{betaBohu})
$$
\LBrace{{l}
\ker\LRf{\b^{\LRt{n}}_k} = \ker\LRf{\beta^{\LRt{n}}_k}\\
co\ker\LRf{\b^{\LRt{n}}_k} = co\ker\LRf{\beta^{\LRt{n}}_k}\\
}
%$$
\quad {\rm where}\quad
%$$
\diagram
G{\CW}_{k}^{\LRt{n}}\LRf{{\BF}} \ar[rr]^{\beta_{k}^{\LRt{n}}} && G{\CW}_{k-1}^{\LRt{n-1}}\LRf{{\BF}}\\
\enddiagram
$$ 
are the connecting homomorphisms corresponding to  the Algebraic Bott homotopy fibration 
\begin{equation}\label{BotIntrrrAug11}
\diagram
G{\CW}^{\LRt{n-1}}\LRf{{\BF}} \ar[r]^{\quad F} & {\BK}\LRf{{\BF}} \ar[r]^{H\quad} & G{\CW}^{\LRt{n}}\LRf{{\BF}}\\
\enddiagram
\end{equation}
of $G{\CW}$-spectra. This establishes the centrality to the Algebraic Bott fibration (\ref{BotIntrrrAug11}),
 in this theory
(see Sec \ref{Secmoidyam}).
%%%%%%%%%%%%%%%%%%%%%%%%%%%%%%%%%%%%%%%%%%%%%%%%%%%%%%%%%
%%%%%%%%%%%%%%%%%%%%%%%%%%%%%%%%%%%%%%%%%%%%%%%%%%%%%%%%%
%%%%%%%%%%%%%%%%%%%%%%%%%%%%%%%%%%%%%%%%%%%%%%%%%%%%%%%%%
%%%%%%%%%%%%%%%%%%%%%%%%%%%%%%%%%%%%%%%%%%%%%%%%%%%%%%%%%

Other than the main results  (\ref{homoThmIntro}, \ref{ThmIntoPd}(\ref{GanesherAshirvadIntro})), stated above, the  foundation and the related machinery we developed, in this article, are of fundamental interest and importance.   In section \ref{SecFunda}, we establish some of the basic definitions and fundamentals of Gersten complexes and Chow-Witt groups. 
 In section \ref{secResidue}, for local integral domains $\LRf{R, \m}$ with $\dim R=1$, we define
the residue homomorphism 
\begin{equation}\label{resiIntro}
\diagram {\BG}W_{p+1}^{\LRt{n+1}}\LRf{K\LRf{R}}\ar[r]^{\partial} &  {\BG}W_{p}^{\LRt{n}}\LRf{C{\BM}^1\LRf{R}}\\ \enddiagram
\end{equation}
in the arrowtic setup, where $K\LRf{R}$ is the  field of fraction of $R$. This parallels the corresponding residue maps,
for discrete valuation rings $\LRf{R, \m}$,
of Milnor  in $K$-theory \cite[pp 322]{M69}, of Witt theory \cite[pp 334]{M69} and 
of Milnor-Witt $K$-theory \cite[pp 55]{M12}. Interestingly, this definition 
(\ref{resiIntro}), works for non-regular rings, while applications in the singular case is yet to be explored. Parallel to 
the (fundamental) split exact sequences in Milnor $K$-theory \cite[pp 325]{M69} , in Witt theory   \cite[pp 335]{M69} and in
Milnor-Witt $K$-theory \cite[pp 60]{M12}, we establish that the following  
\begin{equation}\label{primeIntr0}
\diagram
0\ar[r] &
{G}{\CW}^{\LRt{n+1}}_{p+1}\LRf{{\BF}} \ar[r]^{\pi\quad} &
{G}{\CW}^{\LRt{n+1}}_{p+1}\LRf{{\BF}(T)} \ar[r]^{\bigoplus \partial_y\qquad }  & \bigoplus_{y\in Y^{\LRf{1}}}
{G}{\CW}^{\LRt{n}}_{p}\LRf{C{\BM}^y\LRf{Y_y}}\ar[r] &0 \\
\enddiagram 
\end{equation}  
is a split exact sequence (\ref{arrowticExactAll}), %(\ref{alsoDfme}), 
where ${\BF}$ is a field with $\frac{1}{2}\in {\BF}$ and $Y=\spec{{\BF}\LRt{T}}=:{\BA}^1$.  We define the transfer maps, similar to that in \cite{M69, M12} .
  While the proofs in \cite{M69, M12} are based on generators and relations, our proofs are arrow based. A good
 part of the proof of the splitting of (\ref{primeIntr0}) is sourced out in section  \ref{secHomCom}, which involves
 some  basic homotopy theory. The section 
 \ref{secResidue} roughly parallels \cite[Chap 3]{M12}. This work is a culmination of previous work of the author on the same theme \cite{M25, M25l, M23, M17, M15}. 
 
 In section \ref{SecHomo}, we prove the homotopy invariance theorem (\ref{homoThmIntro}). The description (\ref{ThmIntoPd}) of the Chow-Witt groups is dealt with in section \ref{SecPdDesc}. In section \ref{secChow}, we outline results in arrowtic Chow groups, without proofs. 

%%%%%%%%%%%%%%%%%%%%%%%%%%%%%%%%%%%%%%%%%%%%%%%%%%%%%%%%%
%%%%%%%%%%%%%%%%%%%%%%%%%%%%%%%%%%%%%%%%%%%%%%%%%%%%%%%%%
%%%%%%%%%%%%%%%%%%%%%%%%%%%%%%%%%%%%%%%%%%%%%%%%%%%%%%%%%
%%%%%%%%%%%%%%%%%%%%%%%%%%%%%%%%%%%%%%%%%%%%%%%%%%%%%%%%%
%%%%%%%%%%%%%%%%%%%%%%%%%%%%%%%%%%%%%%%%%%%%%%%%%%%%%%%%%
%%%%%%%%%%%%%%%%%%%%%%%%%%%%%%%%%%%%%%%%%%%%%%%%%%%%%%%%%
%%%%%%%%%%%%%%%%%%%%%%%%%%%%%%%%%%%%%%%%%%%%%%%%%%%%%%%%%
%%%%%%%%%%%%%%%%%%%%%%%%%%%%%%%%%%%%%%%%%%%%%%%%%%%%%%%%%

\section{Fundamentals of Chow Witt groups} \label{SecFunda}%
First, we establish some notations. 
\bN\label{NotaConv}{\rm
We  regard the notations in \cite[App C]{M25l} and \cite{M23} as standard, for the purpose of this article. 
For the convenience of our discussions, we recapitulate some of the notations.
The Grothendiek-Witt ($GW$) groups may be considered as the basic ingredient of our study, in this article. 
Let ${\SE}:=\LRf{{\SE}, {\bfw}, ^{\vee}, \varpi}$  be an exact category with weak equivalences and duality. We assume 
$\frac{1}{2}\in {\SE}$. 
The dg-category of ${\SE}$ is denoted by ${\bf dg}{\SE}$ ({\it this is, basically, the category of bounded chain complexes } 
$Ch^b\LRf{{\SE}}$, {\it with inherited weak equivalences and shifted dualities}). Also, ${\bf dg}^{\LRt{n}}{\SE}$ will denote ${\bf dg}{\SE}$, together with the $n$-shifted duality. 
\bE 
\item The $GW$-groups of the $GW$ space ${\bf GW}\LRf{{\SE}}$, of ${\SE}$, is denoted by $GW_p\LRf{{\SE}}$. 
\item The $GW$-groups of the connected spectrums $G{\CW}^{\LRt{n}}\LRf{{\bf dg}{\SE}}$, of ${\bf dg}{\SE}$, with $n$-shifted duality, is denoted by $G{\CW}_p^{\LRt{n}}\LRf{{\bf dg}{\SE}}$. % or $G{\CW}_p^{\LRt{n}}\LRf{{\SE}}$. 
\item The $GW$-groups of the non-connected bi-spectrums of ${\bf dg}{\SE}$, with $n$-shifted duality, is denoted by ${\BG}{W}_p^{\LRt{n}}\LRf{{\bf dg}{\SE}}$. 

\item 
 Usually, we drop the pre-script ${\bf dg}$ from ${\bf dg}{\SE}$, and   simplify notations:
$$
\LBrace{{l} 
G{\CW}^{\LRt{n}}\LRf{{\SE}}:=G{\CW}^{\LRt{n}}\LRf{{\bf dg}{\SE}}\\
{\BG}{W}^{\LRt{n}}\LRf{{\SE}}={\BG}{W}^{\LRt{n}}\LRf{{\bf dg}{\SE}}\\
}
\quad
\LBrace{{l} 
G{\CW}_p^{\LRt{n}}\LRf{{\SE}}:=G{\CW}_p^{\LRt{n}}\LRf{{\bf dg}{\SE}}\\
{\BG}{W}_p^{\LRt{n}}\LRf{{\SE}}={\BG}{W}_p^{\LRt{n}}\LRf{{\bf dg}{\SE}}\\
}
$$
Since $n$-shifted duality does not make sense for ${\SE}$, by itself, there is no scope of confusion because of  the presence of the super script in $\LRf{-}^{\LRt{n}}$. 
\item Throughout this article, we work with quasi projective schemes $X$ over a noetherian affine scheme 
$\spec{A_0}$, $\frac{1}{2}\in A_0$, with $\dim X=d$. Recall \cite{M23, M25l}, for integers $0\leq k\leq d$,
we  denote %the category of perfect $X$-modules $M$, with $grade(M)=k$, is denoted by $C{\BM}^k\LRf{X}$. So, 
$$
C{\BM}^k\LRf{X}=\LRs{
M\in Coh\LRf{X}: grade\LRf{M}=\dim_{{\SV}\LRf{X}}\LRf{M}=k}
%\qquad where\qquad 
$$
the full subcategory of $Coh\LRf{X}$, of perfect modules of grade $k$,
where   the locally-free dimension of $M$ is denoted by $\dim_{{\SV}\LRf{X}}\LRf{M}$. In particular, $C{\BM}^0\LRf{X}={\SV}\LRf{X}$ is the category of locally free sheaves on $X$. 
Further, ${\CL}$ will denote an invertible sheaf. Note, $M \mapsto M^{\vee}:={\SE}xt^k\LRf{M, {\CL}}$ defines a duality 
on $C{\BM}^k\LRf{X}$ for $k=0, 1, \ldots, d$.  We consider $C{\BM}^k\LRf{X}:=\LRf{C{\BM}^k\LRf{X}, ^{\vee}, \varpi}$ as an exact category with duality. 

\eE
By definition, provided subsequently, the Chow-Witt theory greatly concerns the Gersten complexes of the $GW$-groups of schemes $X$. For non-regular schemes $X$, the Gersten complexes of the $G{\CW}$-groups are not definable, while that of ${\BG}W$-groups are defined \cite[pp 474]{M23}. Because of this,
 we use  ${\BG}W$-notations,  when something works  for non-regular schemes ({\it which are mainly definitions}). 
Our final results in this article concerns only regular schemes, and all such $GW$-groups  coincide (\ref{consSpan}). 
In the regular case, we mostly use the  $G{\CW}$ notation,  sometimes use ${\BG}W$, interchangeably. 

}
\eN

For the purpose of discussions on Chow-Witt groups, following \cite[pp 474]{M23}, we define
the  Gersten ${\BG}W$-complexes.
%we write down the Gersten complexes, as follows. % 
%
\bD\label{notaSptrCx}{\rm 
Let  $X$ be a quasi projective scheme over an affine scheme $\spec{A_0}$, $\frac{1}{2}\in A_0$  and ${\CL}$ be a line bundle on $X$. Let $d=\dim X$. 
%\TC
Then $\forall n, k \in {\BZ}$, with cohomology $\deg=p=0, 1, 2, \ldots, d$, % $\deg=p=p=k, k-1, \ldots, k-d$, 
and for closed subsets $Y \subseteq X$,
define
%$$
%$
%
\begin{equation}\label{labelNOta}
\LBrace{{l}
C^{\LRt{n}pk}\LRf{X, {\CL}}= \bigoplus_{x\in X^{\LRf{p}}}{\BG}W_{k-p}^{\LRt{n+k-p}}\LRf{C{\BM}^{p}\LRf{X_x, {\CL}_x}} \\
C^{\LRt{n}pk}\LRf{X, Y, {\CL}}= \bigoplus_{x\in X^{\LRf{p}}\cap Y}{\BG}W_{k-p}^{\LRt{n+k-p}}\LRf{C{\BM}^{p}\LRf{X_x, {\CL}_x}} \\
}
\end{equation}
For conveniences, write 
$$
C^{\LRt{n}pk}\LRf{X, {\CL}}=C^{\LRt{n}pk}\LRf{X, Y, {\CL}}=0 \quad \forall p\in {\BZ}, \quad {\rm unless}\quad 0\leq p\leq d.
$$
For a fixed $n$, $\LRs{C^{\LRt{n}pk}\LRf{X, {\CL}}: p, k\in {\BZ}}$ is, in deed, a spectral sequence \cite[pp 478]{M23}. In particular,  $C^{\LRt{n}{\bul}k}\LRf{X, {\CL}}$ is a complex. When $X$ is regular, ${\BG}W$-notation will, mostly, be substituted by $G{\CW}$.
}
\eD
We define the (Arrowtic) Chow-Witt groups, as follows.
\bD\label{NotaArrCHWDef}{\rm 
Consider the setup (\ref{labelNOta}). For $k, n\in {\BZ}$, %$p= k, k-1, \ldots, k-d$,} 
define
the  $n$-shifted {\bf Arrowtic Chow-Witt groups}, of codimension $p$ cycles, %and degree $p$, 
as follows
$$
\ora{CW}^{\LRt{n}pk}\LRf{X, {\CL}}=\LBrace{{ll}
 H^p\LRf{C^{\LRt{n}{\bul}k}\LRf{X, {\CL}}}  & p=1, 2, \ldots, d\\
0 & {\rm otherwise} \\
}
$$
Further, for closed $Y \subseteq X$ subsets, it will be convenient to have the following
$$
\ora{CW}^{\LRt{n}pk}\LRf{X, Y, {\CL}}=
\LBrace{{ll}
 H^p\LRf{C^{\LRt{n}{\bul}k}\LRf{X, Y, {\CL}}} & p=0, 1, \ldots, d \\
0 & {\rm otherwise} \\
} 
$$
The latter one is said to be the same, with support in $Y$. 
{\it To keep the notations distinguishable from the Chow-Witt groups $\widetilde{CW}$ in Motivic paradigm, we use $\ora{CW}$ for arrowtic Chow-Witt groups}. 
}
\eD
Our most basic invariants are, as follows.
\bExr\label{milJinh}{\rm 
Let ${\BF}$ be a field with $\frac{1}{2} \in {\BF}$, and $X=\spec{{\BF}}$. In this case, for 
$k\in {\BZ}$ we have 
$$
C^{\LRt{n}pk}\LRf{X}= \LBrace{{ll}
GW_k^{\LRt{k+n}}\LRf{{\BF}} & if~p=0\\
0 & if ~p\neq 0\\
}
$$
Consequently, 
$$
\ora{CW}^{\LRt{n}pk}\LRf{{\BF}}= \LBrace{{ll}
GW_k^{\LRt{k+n}}\LRf{{\BF}} & if~p=0\\
0 & if ~p\neq 0\\
}
$$
}
\eExr
%%%%%%%%%%%%%%%%%%%%%%%%%%%%%%%%%%%%%%%%%%%%

%\bR{\rm 
%To explain the role of $k$, the degree zero term
%$$
%$$
%}\eR
Augmenting the term $ G{\CW}^{\LRt{n+k}}_k\LRf{C{\BM}^0\LRf{X, {\CL}}}$,
at the tail,  the augmented Gersten complex %$C^{\LRt{n}{\bul}k}\LRf{X}:= 
$C^{\LRt{n}{\bul}k}\LRf{X, {\CL}}$ can be displayed as follows: 
$$
\diagram
0 \ar[d]\\
 G{\CW}^{\LRt{n+k}}_k\LRf{C{\BM}^0\LRf{X, {\CL}}}\ar[d]\\
\bigoplus_{x\in X^{\LRf{0}}}G{\CW}^{\LRt{n+k}}_k\LRf{C{\BM}^0\LRf{X_x}}\ar[d]\\
\bigoplus_{x\in X^{\LRf{1}}}G{\CW}^{\LRt{n+k-1}}_{k-1}\LRf{C{\BM}^1\LRf{X_x}}\ar[d]\\
\cdots\ar[d]\\
\bigoplus_{x\in X^{\LRf{k-1}}}G{\CW}^{\LRt{n+1}}_{1}\LRf{C{\BM}^{k-1}\LRf{X_x}}\ar[d] \\
\bigoplus_{x\in X^{\LRf{k}}}G{\CW}^{\LRt{n}}_{0}\LRf{C{\BM}^k\LRf{X_x}}\ar[d]\\ % 
\bigoplus_{x\in X^{\LRf{k+1}}}{\BG}W_{-1}^{\LRt{n-1}}\LRf{C{\BM}^{k+1}\LRf{X_x}} \ar[d]\\
\cdots\ar[d]\\
\bigoplus_{x\in X^{\LRf{d}}}{\BG}W_{-(d-k)}^{\LRt{n-(d-k)}}\LRf{C{\BM}^d\LRf{X_x}} \ar[d] \\
0\\
\enddiagram
$$
% 
%The above was written down, with the  formula \cite[pp 1799]{S17},\cite[pp 466]{M23} 
For non-negative groups \cite[pp 1799]{S17},
\cite[pp 466]{M23}, we have 
$$
 {\BG}W_q^{\LRt{q+n}}\LRf{C{\BM}^{k-q}\LRf{X_x, {\CL}_x}}= {G}{\CW}_q^{\LRt{q+n}}\LRf{C{\BM}^{k-q}\LRf{X_x, {\CL}_x}}
\quad \forall q\geq 0
$$
However, if $X$ is regular, then by Corollary \ref{consSpan}, \cite[pp 1785]{S17} and d\'{e}vissage \cite{BW2}, 
$\forall q\geq 1$, $x\in X^{\LRf{k+q}}$, we have 
$$
 {\BG}W_{-q}^{\LRt{-q+n}}\LRf{C{\BM}^{k+q}\LRf{X_x, {\CL}_x}}=
 W^{\LRt{n}}\LRf{C{\BM}^{k+q}\LRf{X_x, {\CL}_x}} =
 \LBrace{{ll}
 W\LRf{\kappa\LRf{x}} & if~n=0 ~mod~4\\
 0 & {\rm otherwise}\\
 }
$$
Likewise, by d\'{e}vissage \cite[Thm 6.1]{M25}, when $X$ is regular, we have
$$
{\BG}W_q^{\LRt{q+n}}\LRf{C{\BM}^{k-q}\LRf{X_x, {\CL}_x}} = GW_q^{\LRt{q+n}}\LRf{\kappa(x), {\CL}\otimes\kappa(x)}\qquad \forall q
$$
Since shifted duality has periodicity 4, when $X$ is regular,  the augmented Grersten complexes, corresponding to shift $n=0, 1, 2, 3$ can be written down as follows:% ({\it again, augmenting the term at the tail}): 
\begin{equation}\label{TzwiGerK}
\diagram
n=0\quad C^{\LRt{0}{\bul}k}\LRf{X, {\CL}}\\
0 \ar[d]\\
 G{\CW}^{\LRt{k}}_k\LRf{C{\BM}^0\LRf{X, {\CL}}}\ar[d]\\
\bigoplus_{x\in X^{\LRf{0}}}G{\CW}^{\LRt{k}}_k\LRf{C{\BM}^0\LRf{X_x}, {\CL}_x}\ar[d]\\
\bigoplus_{x\in X^{\LRf{1}}}G{\CW}^{\LRt{k-1}}_{k-1}\LRf{C{\BM}^1\LRf{X_x}, {\CL}_x}\ar[d]\\
\cdots\ar[d]\\
\bigoplus_{x\in X^{\LRf{k-1}}}G{\CW}^{\LRt{1}}_{1}\LRf{C{\BM}^{k-1}\LRf{X_x}, {\CL}_x}\ar[d] \\
\bigoplus_{x\in X^{\LRf{k}}}G{\CW}^{\LRt{0}}_{0}\LRf{C{\BM}^k\LRf{X_x}, {\CL}_x}\ar[d]\\ % 
\bigoplus_{x\in X^{\LRf{k+1}}}W^{\LRt{0}}\LRf{C{\BM}^{k+1}\LRf{X_x}, {\CL}_x} \ar[d]\\
\cdots\ar[d]\\
\bigoplus_{x\in X^{\LRf{d}}}W^{\LRt{0}}\LRf{C{\BM}^d\LRf{X_x},  {\CL}_x} \ar[d] \\
0\\
\enddiagram
\quad 
\diagram
n=1, 2, 3 \quad C^{\LRt{n}{\bul}k}\LRf{X, {\CL}}\\
0 \ar[d]\\
 G{\CW}^{\LRt{n+k}}_k\LRf{C{\BM}^0\LRf{X, {\CL}}}\ar[d]\\
\bigoplus_{x\in X^{\LRf{0}}}G{\CW}^{\LRt{n+k}}_k\LRf{C{\BM}^0\LRf{X_x}, {\CL}_x}\ar[d]\\
\bigoplus_{x\in X^{\LRf{1}}}G{\CW}^{\LRt{n+k-1}}_{k-1}\LRf{C{\BM}^1\LRf{X_x}, {\CL}_x}\ar[d]\\
\cdots\ar[d]\\
\bigoplus_{x\in X^{\LRf{k-1}}}G{\CW}^{\LRt{n+1}}_{1}\LRf{C{\BM}^{k-1}\LRf{X_x}, {\CL}_x}\ar[d] \\
\bigoplus_{x\in X^{\LRf{k}}}G{\CW}^{\LRt{n}}_{0}\LRf{C{\BM}^k\LRf{X_x}, {\CL}_x}\ar[d]\\ % 
0\\
\enddiagram
\end{equation}
Now, we have the following simple consequence:
\bP{\rm 
Suppose $X$ is a quasi projective scheme, over $\spec{A_0}$, and ${\CL}$ be a line bundle on $X$. Assume $\frac{1}{2}\in A_0$. Then for $ n=1, 2, 3$. we have the Chow-Witt groups
$$
\ora{CW}^{\LRt{n}pk}\LRf{X, {\CL}}=0 \quad \LBrace{{ll}
 k\leq -1, \forall p &n=1, 2,3 \\
k\geq 0~and~\forall p\geq k+1 &n=1, 2,3 \\
unless~0\leq p \leq d & \forall n \\
}
$$
Further, for $n=0~mod~4$ and  $k\leq -1$, the Gersten complex reduces to the Spectral sequence of 
the Witt groups, as in \cite{BW2}.
}
\eP
\pf Assume $n=1, 2, 3$ modulo $(4)$. When $k\geq 0, p\geq k+1$, it follows from the right hand complex above (\ref{TzwiGerK}).
If $k\leq -1$ and $n=1, 2, 3~mod~(4)$, the Gersten complex is trivial, by the same argument. Same argument works for the last statement. \pic $\eop$

%%%%%%%%%%%%%%%%%%%%%%%%%%%%%%%%%%%%%%%%%%%%%%%%%%%%%
\vspace{3mm} 
%\TCP{13 June Goodnight}\\
We record the following on periodicity, which is part of the standard arguments in the literature \cite{W3}.
\bP\label{Lsquare}{\rm 
Let $X$ be a quasi projective scheme over an $\spec{A_0}$ with $\frac{1}{2}\in A_0$. 
Let ${\CL}_0$ and ${\CL}$ be two invertible bundles over $X$. Then there are homotopy equivalences
$$
\LBrace{{l} 
G{\CW}^{\LRt{n}}\LRf{C{\BM}^k\LRf{X, {\CL}_0}} \cong G{\CW}^{\LRt{n}}\LRf{C{\BM}^k\LRf{X, {\CL}_0{\CL}^2}}   \\
{\BG}{W}^{\LRt{n}}\LRf{C{\BM}^k\LRf{X, {\CL}_0}} \cong {\BG}{W}^{\LRt{n}}\LRf{C{\BM}^k\LRf{X, {\CL}_0{\CL}^2}}  \\
}\quad \forall n, k
$$
of Homotopy spectra (bispectra).
Consequently, the induced maps
$$
\ora{CW}^{\LRt{n}pk}\LRf{X, {\CL}_0} \cong \ora{CW}^{\LRt{n}k}_p\LRf{X, {\CL}_0{\CL}^2} \qquad  \forall n, k, p.
$$
are isomorphisms. 
}
\eP
\pf As usual, let ${\SV}\LRf{X, {\CL}}$ denote the category of vector bundles over $X$, with duality $P \mapsto {\CH}om\LRf{P, {\CL}}$. Consider the functor 
$$
\diagram
{\SV}\LRf{X, {\CL}_0} \ar[rr]^{F} \ar[rr] && {\SV}\LRf{X, {\CL}_0{\CL}^2}\\
\enddiagram
\quad with \quad F_{{\CL}}\LRf{P} =P {{\CL}}
$$
For $P\in {\SV}\LRf{X}$, we have  
$$
F\LRf{P}^{\vee}=
{\CH}om\LRf{P{\CL}, {\CL}_0{\CL}^2}={\CH}om\LRf{P, {\CL}_0{\CL}}=  {\CH}om\LRf{P, {\CL}_0}{\CL}= F\LRf{P^{\vee}}
$$
So, $F$ is a duality preserving equivalence. Rest follows routinely.
\pic $\eop$

%%%%%%%%%%%%%%%%%%%%%%%%%%%%%%%%%%%%%%%%%%%%%%%%%%%%%
%%%%%%%%%%%%%%%%%%%%%%%%%%%%%%%%%%%%%%%%%%%%%%%%%%%%%
%%%%%%%%%%%%%%%%%%%%%%%%%%%%%%%%%%%%%%%%%%%%%%%%%%%%%
%%%%%%%%%%%%%%%%%%%%%%%%%%%%%%%%%%%%%%%%%%%%%%%%%%%%%
%%%%%%%%%%%%%%%%%%%%%%%%%%%%%%%%%%%%%%%%%%%%%%%%%%%%%
%%%%%%%%%%%%%%%%%%%%%%%%%%%%%%%%%%%%%%%%%%%%%%%%%%%%%
%%%%%%%%%%%%%%%%%%%%%%%%%%%%%%%%%%%%%%%%%%%%%%%%%%%%%
%%%%%%%%%%%%%%%%%%%%%%%%%%%%%%%%%%%%%%%%%%%%%%%%%%%%%
%%%%%%%%%%%%%%%%%%%%%%%%%%%%%%%%%%%%%%%%%%%%%%%%%%%%%
%%%%%%%%%%%%%%%%%%%%%%%%%%%%%%%%%%%%%%%%%%%%%%%%%%%%%

\vspace{3mm}
The following is reminiscent  of \cite[Cor 10.4.12]{F8}. 
\bL\label{1111Chwdm}{\rm 
Consider the setup (\ref{notaSptrCx}(\ref{labelNOta}), \ref{NotaArrCHWDef}). So, $X$ is quasi projective scheme over $\spec{A_0}$, with $\frac{1}{2}\in A_0$, $d=\dim X$, and ${\CL}$ is an invertible sheaf. 
Let $Y\subseteq X$ is a closed subset and $U=X-Y$.
Then   there are exact sequences 
\begin{equation}\label{YUXextSeq}
\diagram 
0 \ar[r] & C^{\LRt{n}{\bul}k}\LRf{X, Y, {\CL}} \ar[r] & C^{\LRt{n}{\bul}k}\LRf{X, {\CL}} \ar[r] &  C^{\LRt{n}{\bul}k}\LRf{U, {\CL}} \ar[r] & 0\\
\enddiagram
\qquad \forall n, k\in {\BZ}.
\end{equation} 
of complexes. 
Consequently, there are  long exact sequences of the arrowtic $\ora{CW}$-groups. 
}
\eL
\pf Obvious! $\eop$ 
%%%%%%%%%%%%%%%%%%%%%%%%%%%%%%%%%%%%%%%%%%%%

\subsection{Regular embeddings}
In this section, we expand on the exact sequence (\ref{YUXextSeq}), and the corresponding 
exact sequences of $\ora{CW}$-groups, for regular embeddings $Y\subseteq X$. 
\bL\label{obvPreIsogw}{\rm 
Consider the setup (\ref{labelNOta}). So, $Y\subseteq X$ is a closed subset. Further assume that $X$ and $Y$ are regular, and $Y$ is regularly embedded in $X$, by the ideal sheaf ${\SI}\subseteq {\CO}_X$. Let  
$r=co\dim(Y, X)$. Write 
$$
\LBrace{{l}
{\CO}_Y=\frac{{\CO}_X}{{\SI}}\\
N_{X/Y}=\Lambda^r\frac{{\SI}}{{\SI}^2}\\
}
$$
Then  $\forall r\leq \d \leq d$ and a codimention $\d-r$ point 
$y\in Y^{\LRf{\d-r}}$ there is a commutative diagram natural isomorphisms, as follows
\begin{equation}\label{WGikajahy}
\diagram 
&G{\CW}_p^{\LRt{n}}\LRf{\kappa(y)}\ar@/^/[dr]^{\sim} \ar@/_/[dl]_{\sim}&\\
G{\CW}_p^{\LRt{n}}\LRf{C{\BM}^{\d-r}(Y_y) ,{\CL}_yN_{X/Y}^{-1}}\ar[rr]^{\sim}
&& G{\CW}_p^{\LRt{n}}\LRf{C{\BM}^{\d}(X_y), {\CL}_y}\\
\enddiagram\quad \forall p, n
\end{equation}
where $\kappa(y)$ denotes the residue field of $y$. 
}
\eL
\pf %Here $r=\codim\LRf{X, Y}$. 
We apply (\ref{subsech}) to the inclusion $\LRs{y} \subseteq Y_y \subseteq X_y$. Thus we obtain a commutative diagram of duality 
preserving  functors:
$$
\diagram 
C{\BM}^{\LRs{y}}\LRf{Y_y, \LRf{{\CL}N_{X/Y}^{-1}}_y}\ar@{=}[d]&{\SV}\LRf{\kappa(y)}\ar@/^/[dr] \ar@/_/[dl]&
C{\BM}^{\LRs{y}}(X_y, {\CL}_y)\ar@{=}[d]\\
%&
C{\BM}^{\d-r}\LRf{Y_y, \LRf{{\CL}N_{X/Y}^{-1}}_y} \ar[rr] && C{\BM}^{\d}(X_y, {\CL}_y)\\
\enddiagram 
$$
Operating $G{\CW}_p^{\LRt{n}}$ to this diagram, we obtain the commutative diagram of maps, as in   (\ref{WGikajahy}).
By Lemma \ref{finfpdLem},  two diagonal maps of (\ref{WGikajahy}) are isomorphisms. 
\pic $\eop$

%%%%%%%%%%%%%%%%%%%%%%%%%%%%%%%%%%%%%%%%%%%%
\bL\label{refhKwwdm}{\rm 
Consider the setup (\ref{labelNOta}, \ref{obvPreIsogw}). Then $\forall k\in {\BZ}$,
there is an exact sequence 
\begin{equation}\label{CBulSeqEqn}
\diagram 
0 \ar[r] & C^{\LRt{n}(k-r)}_{\bul}\LRf{Y, {\CL}N_{X/Y}^{-1}} \ar[r] & C^{\LRt{n}k}_{\bul}\LRf{X, {\CL}} \ar[r] &  C^{\LRt{n}k}_{\bul}\LRf{U, {\CL}} \ar[r] & 0\\
\enddiagram
\end{equation}  
of complexes. 
Consequently, there is a long exact sequence of $\ora{CW}$-groups. %homologies. 
}
\eL
\pf Immediate from (\ref{obvPreIsogw}, \ref{1111Chwdm}).  $\eop$ 
%%%%%%%%%%%%%%%%%%%%%%%%%%%%%%%%%%%%%%%%%%%%

%%%%%%%%%%%%%%%%%%%%%%%%%%%%%%%%%%%%%%%%%%%%
%%%%%%%%%%%%%%%%%%%%%%%%%%%%%%%%%%%%%%%%%%%%

\vspace{3mm}
We record the following  exact sequence, which is  
 reminiscent  of \cite[Rem 10.4.12]{F8}.
 %while I am not sure if the \TCP{\bf degree issue was taken care of in  \cite[Cor 10.4.12]{F8}}.
\bP\label{LongSiSediSp}{\rm 
Consider the setup (\ref{labelNOta}, \ref{NotaArrCHWDef}). So, $X$ is a quasi projective scheme over $\spec{A_0}$, with $\frac{1}{2}\in A_0$ and $d=\dim X$. Let
  $Y\subseteq X$ be a closed subset, with $\codim(Y, X)=r\geq 1$, $U=X-Y$  and 
${\CL}$ be a line bundled on $X$. %Further, $U=X-Y$. 
We continue to abuse notations ${\CL}:= {\CL}_{|U}$, as well.
Assume that both $X$ and $Y$ are regular, and $Y$ is regularly embedded in $X$. 
Corresponding to the exact sequence (\ref{CBulSeqEqn}) of complexes, 
 $\forall n, k\in {\BZ}$, we have a long exact sequence:
 \begin{equation}\label{LonSeqEqn1183}
 \diagram
0 \ar[r] &
 \ora{CW}^{\LRt{n}0k}\LRf{X, {\CL}} \ar[r] &  \ora{CW}^{\LRt{n}0k}\LRf{U, {\CL}} \ar@/^/[lld] \\
  0 \ar[r] &
 \ora{CW}^{\LRt{n}1k}\LRf{X, {\CL}} \ar[r] &  \ora{CW}^{\LRt{n}1k}\LRf{U, {\CL}} \ar@/^/[lld] \\
   \cdots \ar[r]& \cdots \ar[r]&  \ora{CW}^{\LRt{n}(r-1)k}\LRf{U, {\CL}}  \ar[dll]\\
   \ora{CW}^{\LRt{n}0(k-r)}\LRf{Y, {\CL}N_{X/Y}^{-1}} \ar[r] &
 \ora{CW}^{\LRt{n}rk}\LRf{X, {\CL}} \ar[r] &  \ora{CW}^{\LRt{n}rk}\LRf{U, {\CL}} \ar@/^/[lld] \\
 \cdots \ar[r]& \cdots \ar[r]& \cdots \ar@/^/[dll]\\
 \ora{CW}^{\LRt{n}(k-r-1)(k-r)}\LRf{Y, {\CL}N_{X/Y}^{-1}} \ar[r] &
 \ora{CW}^{\LRt{n}(k-1)k}\LRf{X, {\CL}} \ar[r] &  \ora{CW}^{\LRt{n}(k-1)k}\LRf{U, {\CL}} \ar@/^/[lld] \\
  \ora{CW}^{\LRt{n}(k-r)(k-r)}\LRf{Y, {\CL}N_{X/Y}^{-1}} \ar[r] &
 \ora{CW}^{\LRt{n}kk}\LRf{X, {\CL}} \ar[r] &  \ora{CW}^{\LRt{n}kk}\LRf{U, {\CL}} \ar@/^/[lld] \\
   \ora{CW}^{\LRt{n}\LRf{k-r+1}(k-r)}\LRf{Y, {\CL}N_{X/Y}^{-1}} \ar[r] &
 \ora{CW}^{\LRt{n}\LRf{k+1}k}\LRf{X, {\CL}} \ar[r] &  \ora{CW}^{\LRt{n}\LRf{k+1}k}\LRf{U, {\CL}} \ar@/^/[lld] \\
  \cdots \ar[r]& \cdots \ar[r]& \cdots \ar[dll]\\
  \ora{CW}^{\LRt{n}\LRf{d-r}(k-r)}\LRf{Y, {\CL}N_{X/Y}^{-1}} \ar[r] &
 \ora{CW}^{\LRt{n}dk}\LRf{X, {\CL}} \ar[r] &  \ora{CW}^{\LRt{n}dk}\LRf{U, {\CL}} \ar@/^/[lld] \\
0 &&\\
 \enddiagram
 \end{equation} 
 Note, if $n=1, 2, 3 ~mod~(4)$ then the sequence terminates at the degree $k$ line. 
 }
\eP
\pf Immediate from (\ref{refhKwwdm}).
$\eop$

%%%%%%%%%%%%%%%%%%%%%%%%%%%%%%%%%%%%%%%%%%%%
%%%%%%%%%%%%%%%%%%%%%%%%%%%%%%%%%%%%%%%%%%%%
%%%%%%%%%%%%%%%%%%%%%%%%%%%%%%%%%%%%%%%%%%%%
%%%%%%%%%%%%%%%%%%%%%%%%%%%%%%%%%%%%%%%%%%%%
%%%%%%%%%%%%%%%%%%%%%%%%%%%%%%%%%%%%%%%%%%%%
%%%%%%%%%%%%%%%%%%%%%%%%%%%%%%%%%%%%%%%%%%%%
%%%%%%%%%%%%%%%%%%%%%%%%%%%%%%%%%%%%%%%%%%%%
%%%%%%%%%%%%%%%%%%%%%%%%%%%%%%%%%%%%%%%%%%%%
%%%%%%%%%%%%%%%%%%%%%%%%%%%%%%%%%%%%%%%%%%%%
%%%%%%%%%%%%%%%%%%%%%%%%%%%%%%%%%%%%%%%%%%%%
%%%%%%%%%%%%%%%%%%%%%%%%%%%%%%%%%%%%%%%%%%%%

%%%%%%%%%%%%%%%%%%%%%%%%%%%%%%%%%%%%%%%%%%%%%%%%%%%%%%%%%
%%%%%%%%%%%%%%%%%%%%%%%%%%%%%%%%%%%%%%%%%%%%%%%%%%%%%%%%%
%%%%%%%%%%%%%%%%%%%%%%%%%%%%%%%%%%%%%%%%%%%%%%%%%%%%%%%%%
%%%%%%%%%%%%%%%%%%%%%%%%%%%%%%%%%%%%%%%%%%%%%%%%%%%%%%%%%
%%%%%%%%%%%%%%%%%%%%%%%%%%%%%%%%%%%%%%%%%%%%%%%%%%%%%%%%%
%%%%%%%%%%%%%%%%%%%%%%%%%%%%%%%%%%%%%%%%%%%%%%%%%%%%%%%%%
%%%%%%%%%%%%%%%%%%%%%%%%%%%%%%%%%%%%%%%%%%%%%%%%%%%%%%%%%

%%%%%%%%%%%%%%%%%%%%%%%%%%%%%%%%%%%%%%%%%%%%
%%%%%%%%%%%%%%%%%%%%%%%%%%%%%%%%%%%%%%%%%%%%
%%%%%%%%%%%%%%%%%%%%%%%%%%%%%%%%%%%%%%%%%%%%

%%%%%%%%%%%%%%%%%%%%%%%%%%%%%%%%%%%%%%%%%%%%
%%%%%%%%%%%%%%%%%%%%%%%%%%%%%%%%%%%%%%%%%%%%
%%%%%%%%%%%%%%%%%%%%%%%%%%%%%%%%%%%%%%%%%%%%
%%%%%%%%%%%%%%%%%%%%%%%%%%%%%%%%%%%%%%%%%%%%%%%%%
\section{The Residue homomorphism}\label{secResidue} 
In this section, we develop our residue homomorphism for 
local integral domains $\LRf{R, \m, \kappa}$, with $\dim R=1$ and associated exact sequences, 
in this arrowtic paradigm. Such residue maps are the key staring point of the theory,
 in all the alternate paradigms. Most notably,  
the Milnor $K$-theory  and the Witt theory versions were given in   \cite[Lem 2.2, Thm 2.3, Cor 5.1, Thm 5.3]{M69}. 
The Motivic version is given in \cite[pp 55, Thm 3.24]{M12}. Recall, Motivic $K$-theory is  also known as Milnor-Witt $K$-Theory. We define the residue homomorphism.

\bD\label{defResidue}{\rm 
Let $\LRf{R, \m, \kappa}$ be a local integral domain, with $\dim R=1$. Let $n, p\in {\BZ}$.
$$
{\rm Denote} 
\quad 
\LBrace{{l}
Y =\spec{R}, \quad
Z=V\LRf{\m}=\LRs{\m}\\
U =Y-Z=\LRs{{\bf 0}}=D(\pi)~\forall \pi \in \m-\LRs{0}\\
So, ~U=\spec{K\LRf{R}}\quad {\rm where}~ K\LRf{R}~{\rm is~the~field~of~fractions.}\\
} 
$$
 By \cite[Thm 5.6]{M25l} there is a homotopy fibration 
$$
\diagram 
{\BG}W^{\LRt{n}}\LRf{C{\BM}^Z\LRf{Y}} \ar[r] & {\BG}W^{\LRt{n+1}}\LRf{{\SV}(Y)} \ar[r] & {\BG}W^{\LRt{n+1}}\LRf{{\SV}(U)} \\
%&&{\BG}W^{\LRt{n+1}}\LRf{K(R)}\ar[u]_{\wr}\\
\enddiagram 
$$
of the ${\BG}W$ bi-spectra. However, ${\BG}W^{\LRt{n+1}}\LRf{{\SV}(U)}\cong {\BG}W^{\LRt{n+1}}\LRf{K(R)}$.
Consequently, we obtain a connecting homomorphisms, to be called the {\bf Residue  map} 
\begin{equation}\label{transfu}
\diagram 
{\BG}W_{p+1}^{\LRt{n+1}}\LRf{K(R)} \ar[rr]^{\partial^{\LRt{n+1}}_{p+1}} && {\BG}W^{\LRt{n}}_p\LRf{C{\BM}^Z\LRf{Y}} 
\enddiagram 
\end{equation}
Sometimes, we write $\partial:=\partial^{\LRt{n+1}}_{p+1}$.
}
\eD
  Clearly, the residue map $\partial$ is part of a longer exact sequence of ${\BG}W$-groups:
\begin{equation}\label{TaruFaru}
\diagram 
\cdots \ar[r] & 
GW^{\LRt{n+1}}_{p+1}\LRf{K(R)} \ar[r]^{\partial} & GW^{\LRt{n}}_{p}\LRf{C{\BM}^Z\LRf{Y}} \ar[dll]\\
GW^{\LRt{n+1}}_{p}\LRf{{\SV}\LRf{Y}} \ar[r] & 
GW^{\LRt{n+1}}_{p}\LRf{K(R)} \ar[r]^{\partial} & \cdots \\ %GW^{\LRt{n}}_{p-1}\LRf{C{\BM}^Z\LRf{Y}} \\
\enddiagram 
\end{equation}
Analogous to other versions alluded above, we propose the following theorem. 
\bT\label{arrowticExactAll} {\rm 
Let ${\BF}$ be a field, with $\frac{1}{2}\in {\BF}$. Write $A={\BF}\LRt{T}$, and $Y:={\BA}^1=\spec{A}$.  For each $y\in Y^{\LRf{1}}$, there is a residue map $\partial_y$  (\ref{transfu}), as follows:
$$
\diagram
{\BG}W^{\LRt{n+1}}_{p+1}\LRf{{\BF}(T)} \ar[rr]^{\partial_y}  &&  
{\BG}W^{\LRt{n}}_{p}\LRf{C{\BM}^y\LRf{Y_y}} \\
\enddiagram
$$ 
 Taking direct sum, with $Y_y=\spec{A_y}$, we obtain the following sequence
\begin{equation}\label{primeSeq}
\diagram
0\ar[r] &
{\BG}W^{\LRt{n+1}}_{p+1}\LRf{{\BF}} \ar[r]^{\p\quad} &
{\BG}W^{\LRt{n+1}}_{p+1}\LRf{{\BF}(T)} \ar[r]^{\bigoplus \partial_y\qquad }  & \bigoplus_{y\in Y^{\LRf{1}}}
{\BG}W^{\LRt{n}}_{p}\LRf{C{\BM}^y\LRf{Y_y}}\ar[r] &0 \\
\enddiagram 
\end{equation}  
of maps. 
The map $\p$ is the usual pullback map. Note that, by d\'{e}vissage \cite[Thm 6.1]{M25}, 
\begin{equation}\label{alsoDfme}
\forall y\in Y^{\LRf{1}} \quad {\BG}W^{\LRt{n}}_{p}\LRf{C{\BM}^y\LRf{Y_y}} \cong 
{\BG}W^{\LRt{n}}_{p}\LRf{\kappa(y)}
\end{equation}
Then the sequence (\ref{primeSeq}) is a split exact sequence.
}
\eT
\pf First, we need the prove that the sequence is a complex. 
Fix $\forall y\in Y^{\LRf{1}}$, we prove that $\partial_y \p=0$. 
Going back to the exact sequence (\ref{TaruFaru}), we have the exact sequence
\begin{equation}\label{Goingkal}
\diagram 
GW^{\LRt{n+1}}_{p+1}\LRf{{\BF}}\ar[d]\ar@/^/[dr]^{\p}&&\\
GW_{p+1}^{\LRt{n+1}}\LRf{{\SV}(Y_y)}\ar[r]_{\iota_y} & GW_{p+1}^{\LRt{n+1}}\LRf{{\BF}(T)}
\ar[r]^{\partial_y} & GW_p^{\LRt{n}}\LRf{C{\BM}^y\LRf{Y_y}}\\
\enddiagram 
\end{equation}
together with the commutative triangle on the left. 
Since $\p$ factors, as above, we have $\partial_y\p=0$. 
 Therefore, the sequence (\ref{primeSeq}) is a complex.  We complete the proof by Lemmas \ref{splitLem} and  \ref{MorExact}. $\eop$

%\item 
\bL\label{splitLem}{\rm 
The map $\p$ is a split monomorphism. 
}
\eL 
\pf  %Refer to %Lemma \ref{KstartopiOne}, and
 By (\ref{unitHike}),  corresponding to $T\in {\BF}\LRt{T}$,  there is a map $T\cup -$, as in the  
 diagram:
\begin{equation}\label{alphaDia}
\diagram 
G{\CW}^{\LRt{n+1}}_{p+1}\LRf{{\BF}(T)} \ar[rr]^{T \cup -} && G{\CW}^{\LRt{n+2}}_{p+2}\LRf{{\BF}(T)}
\ar@{-->}[d]^{\delta_T} \ar@/^/[dr]^{\partial_T}&\\
G{\CW}^{\LRt{n+1}}_{p+1}\LRf{{\BF}} \ar[u]^{\p}\ar[rr]_1&& G{\CW}^{\LRt{n+1}}_{p+1}\LRf{{\BF}} 
\ar[r]_{\d\qquad\qquad}^{\sim\qquad\qquad }&
G{\CW}^{\LRt{n+1}}_{p+1}\LRf{C{\BM}^{(T)}\LRf{{\BF}\LRt{T}_{(T)}}}\\
\enddiagram 
\end{equation} 
such that the outer diagram commutes (\ref{sliceIt}). 
The  arrow  $\d$ is the  d\'{e}vissage isomorphism (\ref{regLocal}). We define  $\delta_T= \d^{-1}\partial_T$. 
%\TCP{The outer diagram commutes \S \ref{secHomCom}}. 
Hence the rectangle commutes.  So, $\p$ splits. $\eop$

\bL\label{MorExact}{\rm 
%Refer to   (\ref{primeSeq}).
The sequence (\ref{primeSeq}) 
 is a split exact sequence,  $\forall n, p$.
}
\eL
\pf By Lemma \ref{splitLem}, the map $\p$ is a split monomorphism. 
Let $Y={\BA}^1=\spec{{\BF}\LRt{T}}$. 
By \cite[pp 473]{M23} the sequences 
$$
\diagram 
{\BG}W^{\LRt{n}}\LRf{C{\BM}^1\LRf{Y}} \ar[r] & {\BG}W^{\LRt{n+1}}\LRf{C{\BM}^0\LRf{Y}} \ar[r]&
{\BG}W^{\LRt{n+1}}\LRf{{\BF}\LRf{T}}\\
\enddiagram 
$$
is a homotopy fibration, of bispectra. By homotopy invariance \cite[Thm 9.8]{S17}, 
$$
 {\BG}W^{\LRt{n+1}}\LRf{C{\BM}^0\LRf{Y}}=  {\BG}W^{\LRt{n+1}}\LRf{{\BF}\LRt{T}}
=  {\BG}W^{\LRt{n+1}}\LRf{{\BF}}
$$
Further, 
$$
{\BG}W^{\LRt{n}}\LRf{C{\BM}^1\LRf{Y}}  =\coprod_{y\in Y^{\LRf{1}}} {\BG}W^{\LRt{n}}\LRf{C{\BM}^y\LRf{Y_y}} 
$$
So, the above homotopy fibfration is homotopic to the sequence
$$
\diagram 
\coprod_{y\in Y^{\LRf{1}}} {\BG}W^{\LRt{n}}\LRf{C{\BM}^y\LRf{Y_y}}  \ar[r] & 
{\BG}W^{\LRt{n+1}}\LRf{{\BF}} \ar[r]&
{\BG}W^{\LRt{n+1}}\LRf{{\BF}\LRf{T}}\\
\enddiagram 
$$
Consequently, we obtain the long exact sequence, of ${\BG}W$ groups
$$
\diagram 
\coprod_{y\in Y^{\LRf{1}}} {\BG}W_{p+1}^{\LRt{n}}\LRf{C{\BM}^y\LRf{Y_y}}  \ar[r] & 
{\BG}W_{p+1}^{\LRt{n+1}}\LRf{{\BF}} \ar[r]^{\p}&
{\BG}W_{p+1}^{\LRt{n+1}}\LRf{{\BF}\LRf{T}}\ar[dll]_{\partial:=\bigoplus \partial_y}\\
\coprod_{y\in Y^{\LRf{1}}} {\BG}W_{p}^{\LRt{n}}\LRf{C{\BM}^y\LRf{Y_y}}  \ar[r] & 
{\BG}W_{p}^{\LRt{n+1}}\LRf{{\BF}} \ar[r]_{\p}&
{\BG}W_{p}^{\LRt{n+1}}\LRf{{\BF}\LRf{T}}\\
\enddiagram 
$$
Both $\p$'s are split injective, by Lemma \ref{splitLem}. Therefore, $\partial$ is surjective.  Hence
$$
\diagram 
0  \ar[r] & 
{\BG}W_{p+1}^{\LRt{n+1}}\LRf{{\BF}} \ar[r]^{\p}&
{\BG}W_{p+1}^{\LRt{n+1}}\LRf{{\BF}\LRf{T}}\ar[rr]^{\partial\qquad} &&
\coprod_{y\in Y^{\LRf{1}}} {\BG}W_{p}^{\LRt{n}}\LRf{C{\BM}^y\LRf{Y_y}}  \ar[r] & 
0\\
\enddiagram 
$$
is exact. $\eop$ 

\vspace{3mm}
This also completes the proof of Theorem \ref{arrowticExactAll}. $\eop$

%%%%%%%%%%%%%%%%%%%%%%%%%%%%%%%%%%%%%%%%%%%%%%
We define the transfer map as follows.
\bD{\rm 
Let ${\BF}$ be a field, with $\frac{1}{2}\in {\BF}$. Let $Y=\spec{{\BF}\LRt{T}}$. Refer to the split exact sequence 
(\ref{primeSeq}). Let $\s$ be a split of $\partial$. Consider the diagram of maps:
$$
\diagram 
{\BG}W_{p+1}^{\LRt{n+1}}\LRf{{\BF}\LRf{T}}\ar[d]_{\partial} \ar@/^/[dr]^1
&\\
\bigoplus_{y\in Y^{\LRf{1}}} {\BG}W_{p}^{\LRt{n}}\LRf{C{\BM}^y\LRf{Y_y}}\ar[r]^{\qquad\s} &  
{\BG}W_{p+1}^{\LRt{n+1}}\LRf{{\BF}\LRf{T}}\ar[d]^{-\partial^{\frac{1}{T}}}\\
{\BG}W_{p}^{\LRt{n}}\LRf{C{\BM}^y\LRf{Y_y}}\ar@{^(->}[u]\ar@{-->}[r]_{\tau^{n, y}_p} & {\BG}W_{p}^{\LRt{n}}\LRf{{\BF}} \\
\enddiagram 
$$
Here $\partial=\bigoplus_{y\in Y^{\LRf{1}}} \partial_y$.
The map $\tau^{n, y}_p$ is defined by composition. This map is 
called the {\bf geometric transfer map}. In fact, $\tau^{n, y}_p$ is 
independent of the choice of the split $\s$.

}
\eD
\pf Consider the composition
$$
\diagram 
{\BG}W^{\LRt{n}}_p\LRf{{\BF}} \ar[r] & {\BG}W^{\LRt{n+1}}_{p+1}\LRf{{\BF}\LRf{T}} \ar[r]^{\quad-\partial^{\frac{1}{T}}} & {\BG}W^{\LRt{n}}_{p}\LRf{{\BF}}\\
\enddiagram  
$$
Repeating the diagram (\ref{Goingkal}), this composition is trivial. Let $\s_1$ be another split of $\partial$. It follows from the sequence (\ref{primeSeq}), $\s_1-\s$ lands in ${\BG}W^{n+1}_{p+1}\LRf{{\BF}}$. 
$$
\forall x\in \bigoplus_{y\in Y^{\LRf{1}}} {\BG}W_{p}^{\LRt{n}}\LRf{C{\BM}^y\LRf{Y_y}} \quad \partial^{\frac{1}{T}} \LRf{\s_1-\s_2}(x)=0
$$
\pic $\eop$

%%%%%%%%%%%%%%%%%%%%%%%%%%%%%%%%%%%%%%%%%%%%%%%%%
\vspace{3mm}
The following is a counter part of similar results in the literature, in the alternate paradigms.
% in the aan analogues of \cite[Lem 1.15]{F20}, in the . 
\bL\label{arrow115}{\rm 
Let ${\BF}$ be a field, with $\frac{1}{2}\in {\BF}$.
%Char\LRf{{\BF}}\neq 2$. 
Let $Y=\spec{{\BF}\LRt{T}}$. Refer to the split exact sequence 
(\ref{primeSeq}).  Then the geometric transfers $\tau^{n, y}_p$ are well defined. Denote
$$
\tau^n:=\bigoplus_{y\in Y^{\LRf{1}}} \tau^{n, y}_p= -\partial^{\frac{1}{T}}\s
$$
 Then the diagram 
$$
\diagram 
0  \ar[r] & 
{\BG}W_{p+1}^{\LRt{n+1}}\LRf{{\BF}} \ar[r]^{\p}&
{\BG}W_{p+1}^{\LRt{n+1}}\LRf{{\BF}\LRf{T}}\ar[r]^{\partial\qquad}  \ar[d]_{-\partial^{\frac{1}{T}}}&
\coprod_{y\in Y^{\LRf{1}}} {\BG}W_{p}^{\LRt{n}}\LRf{C{\BM}^y\LRf{Y_y}}  \ar[r] \ar@/^/[dl]^{\tau^n}& 0\\
&&{\BG}W_p^{\LRt{n}}\LRf{{\BF}}&&\\
\enddiagram 
$$
commutes.
}
\eL
\pf We have 
$$
\tau \partial=  -\partial^{\frac{1}{T}}\s \partial=  -\partial^{\frac{1}{T}}
$$
\pic $\eop$ 
%%%%%%%%%%%%%%%%%%%%%%%%%%%%%%%%%%%%%%%%%%%%%%%%%

\subsection{The local case}
In this section we extend the above (\ref{arrowticExactAll}) to regular local rings.
Before we proceed further, we establish the following setup. 
\begin{setup}\label{chutuApp}{\rm 
Let $\LRf{A, \m, {\BF}}$ be a regular local ring, with $\dim A=d$ and $\frac{1}{2}\in A$.
%Write ${\BF}=\frac{R}{\m}$. 
 Let $R=A\LRt{T}$ be the polynomial ring and $\tilde{\m}={\m}R$. Write $X=\spec{A}$ and $Y=\spec{R}$. Recall,
$C{\BM}^d\LRf{X}$ is the category of $A$-modules with finite length. The duality on
$C{\BM}^d\LRf{X}$  is given by $M \mapsto Ext^d\LRf{M, A}$. Given an  $n\in {\BZ}$, by \cite[pp 473]{M23}, there is a sequence of homotopy fibrations of ${\BG}W$-spectra, together with two natural vertical maps, as follows: 
$$
\diagram 
%& {\BG}W^{\LRt{n+1}}\LRf{{\BF}}\ar[d]^{\wr}&\\
\coprod _{y\in Y^{\LRf{d+1}}}{\BG}W^{\LRt{n}}\LRf{C{\BM}^{d+1}\LRf{Y_y}} &
 {\BG}W^{\LRt{n+1}}\LRf{C{\BM}^{d}\LRf{X}}\ar[d]&\\
{\BG}W^{\LRt{n}}\LRf{C{\BM}^{d+1}\LRf{Y}}\ar[r]\ar[u]_{\wr} & {\BG}W^{\LRt{n+1}}\LRf{C{\BM}^{d}\LRf{Y}}\ar[r] & 
\coprod_{y\in Y^{\LRf{d}}}{\BG}W^{\LRt{n+1}}\LRf{C{\BM}^{d}\LRf{Y_y}}\\
\enddiagram
$$
Here,  the left  vertical map is a homotopy equivalence, and the middle vertical map is the pullback map. 
Consequently, $\forall n, p\in {\BZ}$, there are long exact sequences and maps of $GW$-groups, as follows
\begin{equation}\label{23MayDegDelt}
\diagram 
& {\BG}W^{\LRt{n+1}}_{p+1}\LRf{C{\BM}^{d}\LRf{X}}\ar[d] \ar@/^/@{-->}[dr]^{\i}&\\
\bigoplus _{z\in Y^{\LRf{d+1}}}{\BG}W_{p+1}^{\LRt{n}}\LRf{C{\BM}^{d+1}\LRf{Y_z}}\ar[r] & {\BG}W_{p+1}^{\LRt{n+1}}\LRf{C{\BM}^{d}\LRf{Y}}\ar[r] & 
\bigoplus_{y\in Y^{\LRf{d}}}{\BG}W_{p+1}^{\LRt{n+1}}\LRf{C{\BM}^{d}\LRf{Y_y}}\ar[lld]_{\partial}\\
\bigoplus _{z\in Y^{\LRf{d+1}}}{\BG}W_{p}^{\LRt{n}}\LRf{C{\BM}^{d+1}\LRf{Y_z}}\ar[r] & {\BG}W_{p}^{\LRt{n+1}}\LRf{C{\BM}^{d}\LRf{Y}}\ar[r] & 
\bigoplus_{y\in Y^{\LRf{d}}}{\BG}W_{p}^{\LRt{n+1}}\LRf{C{\BM}^{d}\LRf{Y_y}}\\
\enddiagram
\end{equation}
Write the connecting homomorphism
$$
\partial=\bigoplus_{y\in Y^{\LRf{d}}, z\in Y^{\LRf{d+1}}} \partial_{y2z}\quad {\rm where} \quad 
\diagram  {\BG}W_{p+1}^{\LRt{n+1}}\LRf{C{\BM}^{d}\LRf{Y_y}}\ar[r]^{\partial_{y2z}} & 
{\BG}W_{p}^{\LRt{n}}\LRf{C{\BM}^{d+1}\LRf{Y_z}} \\ \enddiagram
$$
Thus, we obtain the following short complex:
\begin{equation}\label{22Maymtilde2y}
\diagram 
0 \ar[r] &{\BG}W^{\LRt{n+1}}_{p+1}\LRf{C{\BM}^{d}\LRf{X}}
\ar[r]^{\iota} & {\BG}W_{p+1}^{\LRt{n+1}}\LRf{C{\BM}^{d}\LRf{R_{\tilde{m}}}}
\ar[r]^{\partial_{\tilde{\m}}\qquad} & \bigoplus _{z\in Y^{\LRf{d+1}}}{\BG}W_{p}^{\LRt{n}}\LRf{C{\BM}^{d+1}\LRf{Y_z}}\ar[r] & 0\\
\enddiagram
\end{equation}
where 
$$
\partial_{\tilde{\m}}= \bigoplus_{y\in Y^{\LRf{d+1}}}\partial_{\tilde{\m}2y}
$$
}
\end{setup}

%%%%%%%%%%%%%%%%%%%%%%%%%%%%%%%%%%%%%%%%%%%%%%%%%
%%%%%%%%%%%%%%%%%%%%%%%%%%%%%%%%%%%%%%%%%%%%%%%%%
%%%%%%%%%%%%%%%%%%%%%%%%%%%%%%%%%%%%%%%%%%%%%%%%%
%%%%%%%%%%%%%%%%%%%%%%%%%%%%%%%%%%%%%%%%%%%%%%%%%
%%%%%%%%%%%%%%%%%%%%%%%%%%%%%%%%%%%%%%%%%%%%%%%%%
%%%%%%%%%%%%%%%%%%%%%%%%%%%%%%%%%%%%%%%%%%%%%%%%%
%%%%%%%%%%%%%%%%%%%%%%%%%%%%%%%%%%%%%%%%%%%%%%%%%
%%%%%%%%%%%%%%%%%%%%%%%%%%%%%%%%%%%%%%%%%%%%%%%%%

\vspace{3mm} 
Combining with d\'{e}vissage, the following is a consequence of Lemma \ref{MorExact}.
\bT\label{spexactLocal}{\rm 
Consider the setup (\ref{chutuApp}). Then the sequence (\ref{22Maymtilde2y}) is a split exact sequence.
 Further note that 
$$
Y^{\LRf{d+1}}=\LRs{z\in \max Y: z\cap A=\m}.
$$
}
\eT
\pf Write  ${\BA}^1=\spec{{\BF}\LRt{T}}$. Consider the diagram 
$$% 
\diagram 
0  \ar[r] & 
{\BG}W_{p+1}^{\LRt{n+1}}\LRf{{\BF}} \ar[r]^{\p}\ar[d]^{\wr} &
{\BG}W_{p+1}^{\LRt{n+1}}\LRf{{\BF}\LRf{T}}\ar[r]^{\partial\qquad}\ar[d]^{\wr} &
\bigoplus_{z\in \LRf{{\BA}^1}^{\LRf{1}}} {\BG}W_{p}^{\LRt{n}}\LRf{C{\BM}^z\LRf{\LRf{{\BA}^1}_z}}  \ar[r] \ar[d]^{\wr}& 
0\\
0\ar[r]& {\BG}W_{p+1}^{\LRt{n+1}}\LRf{C{\BM}^d\LRf{X}} \ar[r]_{\i} & 
{\BG}W_{p+1}^{\LRt{n+1}}\LRf{C{\BM}^d\LRf{R_{\tilde{\m}}}} \ar[r]_{\partial_{\tilde{\m}}\quad}& \bigoplus_{z\in Y^{\LRf{d+1}}}
{\BG}W^{\LRt{n}}_p\LRf{C{\BM}^z\LRf{{Y}_z}} \ar[r] &0 \\
\enddiagram 
$$ %\end{equation}
 By (\ref{MorExact}) the upper row is a split exact sequence, and by d\'{e}vissage \cite[Thm 6.1]{M25} the vertical maps are isomorphism (see (\ref{alsoDfme}) {\it for the $3^{rd}$ vertical arrow}). It is clear that the left hand rectangle commutes. It remains to show that the right hand rectangle commutes. Note
 $$
 \coprod_{z\in Y^{\LRf{d+1}}}
{\BG}W^{\LRt{n}}\LRf{C{\BM}^z\LRf{{Y}_z}} \cong {\BG}W^{\LRt{n}}\LRf{C{\BM}^{d+1}\LRf{{Y}_z}}
$$ Consider the diagram of arrows of bispectra:
 $$
 \diagram
  \coprod_{z\in \LRf{{\BA}^1}^{\LRf{1}}} {\BG}W^{\LRt{n}}\LRf{C{\BM}^{z}\LRf{{\BA}^1_z}}\ar[r] \ar[d]_{\wr}& {\BG}W^{\LRt{n+1}}\LRf{C{\BM}^{0}\LRf{{\BA}^1}}\ar[r] \ar[d]& 
{\BG}W^{\LRt{n+1}}\LRf{{\BF}\LRf{T}}\ar[d]\\
%&{\BG}W^{\LRt{n+1}}\LRf{C{\BM}^{V\LRf{\tilde{\m}}}\LRf{Y}}\ar[d] \ar[r]^{loc}&{\BG}W^{\LRt{n+1}}\LRf{C{\BM}^{d}\LRf{Y_{\tilde{\m}}}}\ar[d]\\
%
 \coprod_{z\in Y^{\LRf{d+1}}} {\BG}W^{\LRt{n}}\LRf{C{\BM}^{d+1}\LRf{Y_z}}\ar[r] & {\BG}W^{\LRt{n+1}}
 \LRf{C{\BM}^{d}\LRf{Y}}\ar[r] & 
\coprod_{y\in Y^{\LRf{d}}}{\BG}W^{\LRt{n+1}}\LRf{C{\BM}^{d}\LRf{Y_y}}\\
\enddiagram
$$
The %outer 
diagram commutes, and hence this is a map of homotopy fibrations of spectra. So, the long exact sequences of the homotopy groups commutes:
$$% 
\diagram 
%{\BG}W_{p+1}^{\LRt{n+1}}\LRf{{\BF}} \ar[r]^{\p}\ar[d]^{\wr} &
{\BG}W_{p+1}^{\LRt{n+1}}\LRf{{\BF}\LRf{T}}\ar[r]^{\partial\qquad}\ar[d] &
\bigoplus_{y\in \LRf{{\BA}^1}^{\LRf{1}}} {\BG}W_{p}^{\LRt{n}}\LRf{C{\BM}^y\LRf{\LRf{{\BA}^1}_y}}   \ar[d]^{\wr}\\
% {\BG}W_{p+1}^{\LRt{n+1}}\LRf{C{\BM}^d\LRf{X}} \ar[r]_{\i} & 
\bigoplus_{y\in Y^{\LRf{d}}}{\BG}W_{p+1}^{\LRt{n+1}}\LRf{C{\BM}^{d}\LRf{Y_y}}%{\BG}W_{p+1}^{\LRt{n+1}}\LRf{C{\BM}^d\LRf{R_{\tilde{\m}}}} 
\ar[r]_{\partial}& \bigoplus_{z\in Y^{\LRf{d+1}}}
{\BG}W^{\LRt{n}}_p\LRf{C{\BM}^z\LRf{{Y}_z}}  \\
\enddiagram 
$$
where the second $\partial$ is, as in (\ref{23MayDegDelt}). However, the vertical isomorphisms factor, as follows:
$$% 
\diagram 
{\BG}W_{p+1}^{\LRt{n+1}}\LRf{{\BF}\LRf{T}}\ar[r]^{\partial\qquad}\ar[d]^{\wr} &
\bigoplus_{y\in \LRf{{\BA}^1}^{\LRf{1}}} {\BG}W_{p}^{\LRt{n}}\LRf{C{\BM}^y\LRf{\LRf{{\BA}^1}_y}}   \ar[dd]^{\wr}\\
{\BG}W_{p+1}^{\LRt{n+1}}\LRf{C{\BM}^d\LRf{R_{\tilde{\m}}}}\ar@{^(->}[d] \ar[dr]^{\partial_{\tilde{\m}}} & \\
\bigoplus_{y\in Y^{\LRf{d}}}{\BG}W_{p+1}^{\LRt{n+1}}\LRf{C{\BM}^{d}\LRf{Y_y}}% 
\ar[r]_{\partial}& \bigoplus_{z\in Y^{\LRf{d+1}}}
{\BG}W^{\LRt{n}}_p\LRf{C{\BM}^z\LRf{{Y}_z}}  \\
\enddiagram 
$$ 
\pic $\eop$

\section{\bf Homotopy Invariance} \label{SecHomo}%23 May 26: Surjectivity} 

We proceed to establish the
 homotopy invariance of Chow-Witt groups, in the arrowtic paradigm. This is in analogy to 
 \cite[pp 370]{R96} in  Milnor $K$-theory, and \cite[pp 105]{F8} in Milnor-Witt $K$-theory. 
 Having established one of the main tools (\ref{spexactLocal}), first we deal with the following surjectivity result.
\bP%[Thm 11.2.4 pp 112]
\label{Thm1124Sudh}{\rm 
Let $A$ be a regular commutative ring, with $\frac{1}{2}\in A$.  Let $R=A\LRt{T}$ be the polynomial ring. Write $X=\spec{A}$ and 
$Y:=X\times {\BA}^1=\spec{R}$. Let ${\CL}\in {\SV}\LRf{X}$ be an invertible sheaf. 
Let $\pi: Y \lra X$ denote the structure map. Then the pullback map
$$
\diagram 
\ora{CW}^{\LRt{n}pk} \LRf{X, {\CL}} \ar[r]^{\pi^*} \ar[r]& \ora{CW}^{\LRt{n}pk}\LRf{Y, \pi^*{\CL}} \\
\enddiagram 
\quad {\rm is ~surjective}~\forall n, p, k
$$
}
\eP
\pf  Fix $n$ and $k$.  Consider part of the diagrams below \cite[pp 473]{M23}, to construct the 
Gersten complexes. We  ignore the ${\CL}$-coordinate, for typographical reasons.
{\scalefont{.7}
\begin{equation}\label{473GerXuu}
\diagram 
\bigoplus_{X^{\LRf{p-1}}}G{\CW}_{k-p+1}^{\LRt{n+k-p+1}}\LRf{C{\BM}^{p-1}\LRf{X_x}}\ar[d] 
\ar@{-->}[dr]^{\partial_{p-1}^X}&& G{\CW}_{k-p-1}^{\LRt{n+k-p-2}}\LRf{C{\BM}^{p+2}\LRf{X}}\ar[d]\\
 G{\CW}_{k-p}^{\LRt{n+k-p}}\LRf{C{\BM}^{p}\LRf{X}} \ar[r]\ar[d]&
 \bigoplus_{X^{\LRf{p}}}G{\CW}_{k-p}^{\LRt{n+k-p}}\LRf{C{\BM}^{p}\LRf{X_x}} \ar[r]\ar@{-->}[dr]^{\partial_{p}^X}
 &G{\CW}_{k-p-1}^{\LRt{n+k-p-1}}\LRf{C{\BM}^{p+1}\LRf{X}}\ar[d] \\
 G{\CW}_{k-p}^{\LRt{n+k-p+1}}\LRf{C{\BM}^{p-1}\LRf{X}}  &&  
 \bigoplus_{X^{\LRf{p+1}}}G{\CW}_{k-p-1}^{\LRt{n+k-p-1}}\LRf{C{\BM}^{p+1}\LRf{X_x}}\\
\enddiagram 
\end{equation}
}
Similarly, for $Y=X\times {\BA}^1$, we have
{\scalefont{.7}
\begin{equation}\label{YMaroGerXuu}
\diagram 
\bigoplus_{Y^{\LRf{p-1}}}G{\CW}_{k-p+1}^{\LRt{n+k-p+1}}\LRf{C{\BM}^{p-1}\LRf{Y_y}}\ar[d] 
\ar@{-->}[dr]^{\partial_{p-1}^Y}&& G{\CW}_{k-p-1}^{\LRt{n+k-p-2}}\LRf{C{\BM}^{p+2}\LRf{Y}}\ar[d]\\
 G{\CW}_{k-p}^{\LRt{n+k-p}}\LRf{C{\BM}^{p}\LRf{Y}} \ar[r]\ar[d]&
 \bigoplus_{Y^{\LRf{p}}}G{\CW}_{k-p}^{\LRt{n+k-p}}\LRf{C{\BM}^{p}\LRf{Y_y}} \ar[r]\ar@{-->}[dr]^{\partial_{p}^Y}
 &G{\CW}_{k-p-1}^{\LRt{n+k-p-1}}\LRf{C{\BM}^{p+1}\LRf{Y}}\ar[d] \\
 G{\CW}_{k-p}^{\LRt{n+k-p+1}}\LRf{C{\BM}^{p-1}\LRf{Y}}  &&  
 \bigoplus_{Y^{\LRf{p+1}}}G{\CW}_{k-p-1}^{\LRt{n+k-p-1}}\LRf{C{\BM}^{p+1}\LRf{Y_y}}\\
\enddiagram 
\end{equation}
}
%%%
%}
Denote 
$$
\LBrace{{l}
Y^{\LRf{p}}_{e}= \LRs{y\in Y^{\LRf{p}}: y\cap A\in X^{\LRf{p}} }\\
Y^{\LRf{p}}_{ne}= \LRs{y\in Y^{\LRf{p}}: y\cap A\in X^{\LRf{p-1}} }\\
}
$$ 
More specifically,  
$$
\begin{array}{ll}
 \forall x\in X^{\LRf{p}} ~{\rm denote}&Y^{\LRf{p}}_{e, x}= \LRs{y\in Y^{\LRf{p}}_e: y\cap A=x }=\LRs{xA}\\ %\hline 
 \forall x\in X^{\LRf{p-1}}~{\rm denote} &Y^{\LRf{p}}_{ne, x}= \LRs{y\in Y^{\LRf{p}}_{ne}: y\cap A=x }\\
\end{array}
$$ 
Let $\omega \in \ker\LRf{\partial_p^Y}$.  We can write 
$$
\omega=\alpha +\beta\quad 
\LBrace{{l}
\alpha \in  \bigoplus_{Y^{\LRf{p}}_e}G{\CW}_{k-p}^{\LRt{n+k-p}}\LRf{C{\BM}^{p}\LRf{Y_y}, \pi^*{\CL}_y}\\
\beta \in  \bigoplus_{Y^{\LRf{p}}_{ne}}G{\CW}_{k-p}^{\LRt{n+k-p}}\LRf{C{\BM}^{p}\LRf{Y_y}, \pi^*{\CL}_y}\\
}
$$
First, we show that 
\begin{equation}\label{betaBunryty}
\forall \beta\in \bigoplus_{Y^{\LRf{p}}_{ne}}
G{\CW}_{k-p}^{\LRt{n+k-p}}\LRf{C{\BM}^{p}\LRf{Y_y}, \pi^*{\CL}_y}, \quad
\LBrace{{l} \beta ~{\rm is ~in~the~boundary, ~and~ hence}\\
\LRt{\beta}=0 \in \ora{CW}^{\LRt{n}pk}\LRf{Y, \pi^*{\CL}}
}
\end{equation} 
 %and hence  $\LRt{\beta}=0$ in $\ora{CH}_p^{\LRt\LRf{Y, \pi^*{\CL}}$. 
We can write 
$$
\beta=\beta_{x_1}+\beta_{x_2}+ \cdots + \beta_{x_r} \quad {\rm with} % \quad 
\LBrace{{l}
x_s\in X^{\LRf{p-1}}\quad \forall~ s=1, 2, \ldots, r, \quad {\rm and}\\
\beta_{x_s}\in \bigoplus_{y\in Y^{\LRf{p}}_{ne, x_s}}G{\CW}_{k-p}^{\LRt{n+k-p}}\LRf{C{\BM}^{p}\LRf{Y_y}, \pi^*{\CL}_y}\\
}
$$
It is enough to show  $\beta_{x_s}\in Image\LRf{\partial^Y_{p-1}}$. So, we write $x=x_s$ and $\beta=\beta_x$, with 
$x\in X^{\LRf{p-1}}$.  Denote the extended ideal $\tilde{x}=xR$.
%It is enough to prove that each  \beta_{x_s} \in Image\LRf{\partial^Y_{p+1}}$. 
It follows from  (\ref{spexactLocal}) that the sequence 
%{\scalefont{.7}
$$
\diagram 
0 \ar[r] & {G}{\CW}^{\LRt{n+k-p+1}}_{k-p+1}\LRf{C{\BM}^{p-1}\LRf{X_{x}}}
\ar[r]^{\iota} & {G}{\CW}_{k-p+1}^{\LRt{n+k-p+1}}\LRf{C{\BM}^{p-1}\LRf{R_{\tilde{x}}}}
\ar[ld]^{\partial_{\tilde{x}}\qquad} \\
& \bigoplus _{z\in Y^{\LRf{p}}_{ne, x}}{G}{\CW}_{k-p}^{\LRt{n+k-p}}\LRf{C{\BM}^{p}\LRf{Y_z}}\ar[r] & 0\\
\enddiagram
$$
%}
is exact.  %  
From surjectivity of $\partial_{\tilde{x}}$ it follows $\partial_{p+1}^Y\LRf{\gamma} =\beta$ for some $\gamma$. Therefore $\beta= \partial_{p+1}^Y\LRf{\gamma}$ for some $\gamma$ and $\beta$ is a boundary. So, (\ref{betaBunryty}) is established. 
Therefore, we  assume that $\beta=0$ and $\omega=\alpha$.
We can write 
$$
\omega= \omega_1 + \cdots + \omega_r\quad {\rm where}  \quad \omega_i\in G{\CW}^{\LRt{n+k-p}}_{k-p}\LRf{C{\BM}^{p}\LRf{Y_{y_i}}, 
\pi^*{\CL}_{y_i}}, ~{\rm with} ~y_i\in Y_e^{\LRf{p}}. 
$$
We assume $y_i\neq y_j~\forall i\neq j$. 
Write $x_i=y_i\cap A$ and hence $y_i=x_iA$. Fix $i$. By 
 (\ref{spexactLocal}), the sequence 
%{\scalefont{.7}
$$
\diagram 
0 \ar[r] & {G}{\CW}^{\LRt{n+k-p}}_{k-p}\LRf{C{\BM}^{p}\LRf{X_{x_i}}}
\ar[r]^{\pi_i^*} & {G}{\CW}_{k-p}^{\LRt{n+k-p}}\LRf{C{\BM}^{p}\LRf{R_{\tilde{y}_i}}}
\ar[dl]^{\partial_{\tilde{x}_i}} \\
& \bigoplus _{z\in Y^{\LRf{p+1}}_{e, x_i}}{G}{\CW}_{k-p-1}^{\LRt{n+k-p-1}}\LRf{C{\BM}^{p+1}\LRf{Y_z}}\ar[r] & 0\\
\enddiagram
$$
%}
is exact. % with $\m_i= A_{x_i}$, $\tilde{\m}_i=\m_iA_{\m_i} \subseteq A_{\m_i}$.
Note $Y^{\LRf{k-p+1}}_{x_i}\cap Y^{\LRf{k-p+1}}_{x_j}=\phi$, whenever $i\neq j$.
Here, the maps $\pi_i^*$ is induced by $\pi: Y \lra X$.
Since $\partial^Y_p(\omega)=0$ we have $\partial_{\tilde{x}_i}\LRf{\omega_i}=0$. 
Therefore, $\omega_i= \pi^*_i\LRf{\gamma_i}$, for some $\gamma_i\in {\BG}W^{\LRt{n+k-p}}_{k-p}\LRf{C{\BM}^{p}\LRf{X_{x_i}}}$. 
To complete the proof, we need to prove that such an element $\gamma_i$ is a cycle. Consider the commutative diagram 
\begin{equation}\label{preuveSur}
\diagram 
{G}{\CW}^{\LRt{n+k-p}}_{k-p}\LRf{C{\BM}^{p}\LRf{X_{x_i}}}\ar[d]_{\partial^X_p}
\ar[rr]^{\pi_i^*} && {G}{\CW}_{k-p}^{\LRt{n+k-p}}\LRf{C{\BM}^{p}\LRf{R_{\tilde{y}_i}}}
\ar[d]^{\partial_p^Y } \\
% 
%%%
 \bigoplus _{x\in X^{\LRf{p+1}}} {G}{\CW}_{k-p-1}^{\LRt{n+k-p-1}}\LRf{C{\BM}^{p+1}\LRf{X_x}}\ar[rr]_{\f}&&
 \bigoplus _{z\in Y^{\LRf{p+1}}}{G}{\CW}_{k-p-1}^{\LRt{n+k-p-1}}\LRf{C{\BM}^{p+1}\LRf{Y_z}}\\ %&&\\
\enddiagram 
\end{equation}
where $\f$ is the direct sum of the maps, induced by $\pi$.
%$$
%  
Therefore, $\f$ is injective, by (\ref{spexactLocal}). Therefore, $\gamma_i$ are cycles, and $\LRt{\omega_i}$ is in the
image on $\pi^*$.  
Since this works for each $i$, $\omega$ is in the image of $\pi^*$. 
\pic $\eop$ 

\vspace{3mm} 
The following, established in the proof of  (\ref{Thm1124Sudh}), will be of our future interest. 
\bC\label{Amnii}{\rm 
Consider the setup of (\ref{Thm1124Sudh}).
Suppose 
%$$
%$$
$$
\beta\in \bigoplus_{Y^{\LRf{p}}_{ne}}G{\CW}_{k-p}^{\LRt{n}}\LRf{C{\BM}^{p}\LRf{Y_y}, \pi^*{\CL}_y}
$$
 Then $\beta
\in image\LRf{\partial^Y_{p-1}}$, is a boundary.
}
\eC
\pf See (\ref{betaBunryty}) in the  proof of  (\ref{Thm1124Sudh}). $\eop$ 

\vspace{3mm}
\bT%[Cor 11.2.6]
\label{CorFus1126}
{\rm 
Consider the setup of (\ref{Thm1124Sudh}). Then the pullback map
$$
\diagram 
\ora{CW}^{\LRt{n}pk} \LRf{X, {\CL}} \ar[r]^{\pi^*}_{\sim}& \ora{CW}^{\LRt{n}pk}\LRf{Y, \pi^*{\CL}} \\
\enddiagram \quad
{\rm is ~an~ isomorphism}\quad \forall p, n, k\in {\BZ}, 
$$
%is an isomorphism, $\forall p, n, k\in {\BZ}, k\geq 0$. 
}
\eT

\pf We continue to  use some of the notations from  the proof of  (\ref{Thm1124Sudh}). We combine the diagrams 
(\ref{473GerXuu}) and (\ref{YMaroGerXuu}), as follows:
{\scalefont{.7}
$$
\diagram 
\bigoplus_{X^{\LRf{p-1}}}G{\CW}_{k-p+1}^{\LRt{n+k-p+1}}\LRf{C{\BM}^{p-1}\LRf{X_x}}
\ar[r]^{\qquad\partial_{p-1}^X} \ar@{^(->}[d]_{\pi^*_{p-1}}&
 \bigoplus_{X^{\LRf{p}}}G{\CW}_{k-p}^{\LRt{n+k-p}}\LRf{C{\BM}^{p}\LRf{X_x}} \ar[r]^{\partial_{p}^X\quad}\ar@{^(->}[d]_{\pi^*_{p}}
& \bigoplus_{X^{\LRf{p+1}}}G{\CW}_{k-p-1}^{\LRt{n+k-p-1}}\LRf{C{\BM}^{p+1}\LRf{X_x}}\ar@{^(->}[d]_{\pi^*_{p+1}}\\
%%%
\bigoplus_{Y^{\LRf{p-1}}}G{\CW}_{k-p+1}^{\LRt{n+k-p+1}}\LRf{C{\BM}^{p-1}\LRf{Y_y}}
\ar[r]^{\qquad \partial_{p-1}^Y} 
&
 \bigoplus_{Y^{\LRf{p}}}G{\CW}_{k-p}^{\LRt{n+k-p}}\LRf{C{\BM}^{p}\LRf{Y_y}} \ar[r]^{\partial_{p}^Y\quad}
&  \bigoplus_{Y^{\LRf{p+1}}}G{\CW}_{k-p-1}^{\LRt{n+k-p-1}}\LRf{C{\BM}^{p+1}\LRf{Y_y}}\\
\enddiagram 
$$
}
The vertical maps $\pi^*_{j}$,  are induced by the map $\pi: Y \lra X$. In fact, $\pi^*_j$ is injective, $\forall j$, by (\ref{spexactLocal}).  
Let 
$$
\omega\in \bigoplus_{x\in X^{\LRf{p}}}G{\CW}_{k-p}^{\LRt{n+k-p}}\LRf{C{\BM}^{p}\LRf{X_x}, {\CL}_x} 
\quad \ni~~
\LBrace{{l}
\partial_p^X\LRf{\omega}=0, ~ and \\
\pi^*_p\LRf{\LRt{\omega}}=0\in \ora{CW}^{\LRt{n}pk}\LRf{Y, \pi^*{\CL}}\\
}
$$
Therefore, we have 
$$
\partial_{p-1}^Y\LRf{\gamma}= \pi^*_p\LRf{\omega} \quad {\rm for ~some} \quad\gamma\in
 \bigoplus_{y\in Y^{\LRf{p-1}}} 
G{\CW}^{\LRt{n+k-p+1}}_{k-p+1}
\LRf{C{\BM}^{p-1}\LRf{Y_y}}
$$
Write 
$$
\gamma=\gamma_e +\gamma_{ne} \quad {\rm where}~ 
\LBrace{{l}
\gamma_e\in \bigoplus_{y\in Y^{\LRf{p-1}}_e} 
G{\CW}^{\LRt{n+k-p+1}}_{k-p+1}\LRf{C{\BM}^{p-1}\LRf{Y_y}}\\
\gamma_{ne}\in \bigoplus_{y\in Y^{\LRf{p-1}}_{ne}} 
G{\CW}^{\LRt{n+k-p+1}}_{k-p+1}\LRf{C{\BM}^{p-1}\LRf{Y_y}}\\
}
$$
By (\ref{Amnii}), $\gamma_{ne}$ is a boundary. So, we can assume that $\gamma_{ne}=0$ and 
$$
\gamma=\gamma_e\in  \bigoplus_{y\in Y^{\LRf{p-1}}_e} 
G{\CW}^{\LRt{n+k-p+1}}_{k-p+1}\LRf{C{\BM}^{p-1}\LRf{Y_y}}
$$
We can write 
$$
\gamma=\sum_{i=1}^r\gamma_i \qquad {\rm with} \quad 
\LBrace{{l} \gamma_i\in G{\CW}^{\LRt{n+k-p+1}}_{k-p+1}\LRf{C{\BM}^{p-1}\LRf{Y_{\tilde{x}_i}}} \\
{\rm where}~x_i\in X^{\LRf{p-1}},  \tilde{x}_i=x_iA\LRt{T} \\
}
$$
By (\ref{spexactLocal}), we have the following split exact sequences:
%{\scalefont{.7}
$$ 
\diagram 
0 \ar[r] &{G}{\CW}^{\LRt{n+k-p+1}}_{k-p+1}\LRf{C{\BM}^{p-1}\LRf{X_{x_i}}}
\ar[r]^{\pi^*_{x_i}} & {G}{\CW}_{k-p+1}^{\LRt{n+k-p+1}}\LRf{C{\BM}^{p-1}\LRf{Y_{\tilde{x}_i}}}
\ar[dl]^{\partial_{\tilde{x}_i}} \\
& \bigoplus _{z\in Y_{ne, x_i}^{\LRf{p}}}{G}{\CW}_{k-p}^{\LRt{n+k-p}}\LRf{C{\BM}^{p}\LRf{Y_z}}\ar[r] & 0\\
\enddiagram
$$
%}
%\TCP
%
%}
We have 
$$
\partial^Y_{p-1}\LRf{\gamma} =\sum_{i=1}^r\partial^Y_{p-1}\LRf{\gamma_i} =\pi^*_p\LRf{\omega} \in 
\bigoplus_{y\in Y_e^{\LRf{p}}}G{\CW}^{\LRt{n+k-p}}_{k-p}\LRf{C{\BM}^{p}\LRf{Y_y}}
$$
Since $\partial_{p-1}^Y\LRf{\gamma}$ is supported on $Y^{\LRf{p}}_{ne}$,
it follows   $\partial_{\tilde{x}_i}\LRf{\gamma_i}=0$. Hence $\gamma_i=\pi_{x_i}^*\LRf{\beta_i}$ for some 
$\beta_i$. Adding, we have 
$$
\gamma=\pi^*_{p-1}\LRf{\beta}\quad {\rm for ~some}\quad \beta\in 
\bigoplus_{x\in X^{\LRf{p-1}}}G{\CW}_{k-p+1}^{\LRt{n+k-p+1}}\LRf{C{\BM}^{p-1}\LRf{X_x}}
$$
Therefore,
$$
\pi^*_p\LRf{\omega}= \partial^Y_{p-1}\LRf{\gamma} = \partial^Y_{p-1}\LRf{\pi^*_{p-1}\LRf{\beta}} 
=\pi_p^*\LRf{\partial_{p-1}^X\LRf{\beta}}
$$
By injectivity of $\pi_p^*$, we have $\omega=\partial_{p+1}^X\LRf{\beta}$. So, $\LRt{\omega}=0$.
\pic $\eop$

%%%%%%%%%%%%%%%%%%%%%%%%%%%%%%%%%%%%%%%%%%%%%%%%%%
\vspace{3mm}
Now we proceed to give various other versions of the Homotopy theorem (\ref{CorFus1126}), which follow from  one another.
%%%%%%%%%%%%%%%%%%%%%%%%%%%%%%%%%%%%%%%%%%%%%%%%%%
\bT\label{Thm1124WithSupport}{\rm 
Let $A$ be a regular commutative ring and $R=A\LRt{T}$ be the polynomial ring, with $\frac{1}{2}\in A$. Write $X=\spec{A}$ and 
$Y:=X\times {\BA}^1=\spec{R}$. Also let $B\subseteq X$ be a closed subset of $X$. Let ${\CL}\in {\SV}\LRf{X}$
be an invertible sheaf. 
Let $\pi: Y \lra X$ be the structure map. Then $\forall p, k, n$,  the induced map
$$
\diagram 
\ora{CW}^{\LRt{n}pk} \LRf{X, B, {\CL}} \ar[rr]^{\pi^*\qquad}_{\sim\qquad} && \ora{CW}^{\LRt{n}pk}\LRf{Y, \pi^{-1}B, \pi^*{\CL}} \\
\enddiagram 
$$
is  an isomorphism.  
}
\eT
\pf The arguments in the proofs of (\ref{Thm1124Sudh}, \ref{CorFus1126}) are  local on the base scheme $X$. The proofs go through line for line, in this case also.
$\eop$  
%%%%%%%%%%%%%%%%%%%%%%%%%%%%%%%%%%%%%%%%%%%%%%%%%%

%%%%%%%%%%%%%%%%%%%%%%%%%%%%%%%%%%%%%%%%%%%%%%%%%%
\bT\label{ThAffnSpaces}{\rm 
Let $A$ be a regular commutative ring with $\frac{1}{2}\in A$. Let  $R=A\LRt{T_1, \ldots, T_r}$ be the polynomial ring, in $n$ variables. Write $X=\spec{A}$ and 
$Y=X\times {\BA}^r=\spec{R}$. Also let $B\subseteq X$ be a closed subset of $X$. Let ${\CL}\in {\SV}\LRf{X}$
be a line bundle on $X$.  
Let $\pi: Y \lra X$ be the structure map. Then $\forall p, k, n$, the induced map
$$
\diagram 
\ora{CW}^{\LRt{n}pk} \LRf{X, B, {\CL}} \ar[rr]^{\pi^*\qquad}_{\sim\qquad} && \ora{CW}^{\LRt{n}pk}\LRf{X\times {\BA}^r, \pi^{-1}B, \pi^*{\CL}} \\
\enddiagram 
$$
is  an isomorphism.  
}
\eT
\pf Follows from (\ref{Thm1124WithSupport}), by induction. $\eop$
%%%%%%%%%%%%%%%%%%%%%%%%%%%%%%%%%%%%%%%%%%%%%%%%%%

%%%%%%%%%%%%%%%%%%%%%%%%%%%%%%%%%%%%%%%%%%%%%%%%%%
\bT\label{VectorBunCW}{\rm 
Let $X$ be quasi projective scheme over a noetherian affine scheme $\spec{A_0}$. Assume   $X$ is regular,
with $\frac{1}{2}\in A_0$.  Let 
$\pi: Y \lra X$ be a vector bundle, meaning, $Y=\spec{\Sym{{\SE}}}$, where ${\SE}\in {\SV}\LRf{X}$. 
Let ${\CL}\in {\SV}\LRf{X}$ be a line bundle. All the dualities are induced by the duality $P \mapsto {\CH}om\LRf{-, {\CL}}$ on ${\SV}\LRf{X}$.
 Then $\forall p, k, n$, the pullback map
$$
\diagram 
\ora{CW}^{\LRt{n}pk} \LRf{X,  {\CL}} \ar[rr]^{\pi^*\quad}_{\sim\quad} && \ora{CW}^{\LRt{n}pk}\LRf{Y,  \pi^*{\CL}} \\
\enddiagram 
$$
is  an isomorphism.  
}
\eT
\pf Let $U_0=\spec{R_0} \subseteq X$ be an affine open subset, such that ${\SE}_{|U_0}\cong R_0^r$ is trivial. 
Then $\pi^{-1}U_0\cong \spec{R_0\LRt{T_1, \ldots, T_r}}=X\times {\BA}^r$. By (\ref{ThAffnSpaces}), the 
pullback maps 
$$
\diagram 
\ora{CW}^{\LRt{n}pk} \LRf{U_0,  {\CL}_{|{U_0}}} \ar[rr]^{\pi^*\quad}_{\sim\quad} && \ora{CW}^{\LRt{n}pk}\LRf{\pi^{-1}U_0,  \pi^*{\CL}_{|U_0}} \\
\enddiagram 
$$
are isomorphisms, $\forall n, k, p$. Let $B_0=X-U_0$ be the closed complement. Then there is a natural map
$$
\diagram 
C^{\LRt{n}{\bul}k}\LRf{X, B_0, {\CL}} \ar[r] \ar[d]& C^{\LRt{n}{\bul}k}\LRf{X, {\CL}} \ar[r] \ar[d]^{\pi^*}& C^{\LRt{n}{\bul}k}\LRf{U_0, {\CL}_{|U_0}}\ar[d]\\
C^{\LRt{n}{\bul}k}\LRf{Y, \pi^{-1}B_0, \pi^*{\CL}} \ar[r] & C^{\LRt{n}{\bul}k}\LRf{Y, {\CL}} \ar[r] & C^{\LRt{n}{\bul}k}\LRf{\pi^{-1}U_0, \pi^{-1}{\CL}_{|U_0}}\\
\enddiagram 
$$
exact sequences of complexes. 
It is established that the third vertical map induces isomorphisms of homologies.
So,  we only need to prove that the first vertical maps are isomorphisms.  Note that $X$ has a finite affine open cover
$\LRs{U_i: i=1, \ldots, r}$, such that ${\SE}_{|U_i}$ are trivial.
%which  that trivializes ${\SE}$. 
So, there are affine open sets $U_1=\spec{R_1}\subseteq X$ such that  ${\SE}_{|U_1}\cong R_1^r$ and $B_0\cap U_1\neq \phi$. Write $B_1=B_0-U_1$. So, $B_1\neq B_0$. As before, there are exact sequences of complexes:
$$
\diagram 
0 \ar[r] & C^{\LRt{n}{\bul}k}\LRf{X, B_1} \ar[r]\ar[d] & C^{\LRt{n}{\bul}k}\LRf{X, B_0} \ar[r]\ar[d] & 
C^{\LRt{n}{\bul}k}\LRf{U_1, U_1\cap B_0} \ar[r]\ar[d] & 0\\
0 \ar[r] & C^{\LRt{n}{\bul}k}\LRf{Y, \pi^{-1}B_1} \ar[r] & C^{\LRt{n}{\bul}k}\LRf{Y, \pi^{-1}B_0} \ar[r] & 
C^{\LRt{n}{\bul}k}\LRf{\pi^{-1}U_1, \pi^{-1}\LRf{U_1\cap B_0}} \ar[r] & 0\\
\enddiagram 
$$
The $3^{rd}$ vertical arrow induces isomorphisms of cohomologies, by (\ref{ThAffnSpaces}). To establish that the middle vertical arrow  induces isomorphisms of cohomologies, it is enough to prove that the first vertical arrow does the same.
Repeating the process,
%The last term is fine, because it corresponds to $U_1\times {\BA}^1$. So, it is enough to prove for $B_1$.
we obtain a strictly decreasing sequence of closed sets $\cdots \subseteq B_n \subseteq \cdots \subseteq B_1\subseteq B_0$. By noetherian property, the process stops when $B_n=\LRs{z}$ is a closed point in $X$, and it is enough to prove that the map 
$$
\diagram 
C^{\LRt{n}{\bul}k}\LRf{X, \LRs{z}, {\CL}} \ar[r] &
C^{\LRt{n}{\bul}k}\LRf{Y, \pi^{-1}\LRs{z}, \pi^*{\CL}} \\
\enddiagram
$$
induce isomorphism of cohomologies. Recall $d=\dim X$. In this case, Only nonzero term in
$C^{\LRt{n}{\bul}k}\LRf{X, \LRs{z}, {\CL}}$ is at $\deg=p=d$, as follows
$$
\begin{array}{l}
C^{\LRt{n}dk}\LRf{X, \LRs{z}, {\CL}} =G{\CW}_{k-d}^{\LRt{n+k-d}}\LRf{C{\BM}^d\LRf{X_{z}}}
\\
= \LBrace{{ll} 
G{\CW}_0^{\LRt{n}}\LRf{\kappa\LRf{z}}& if~k=d\\
W\LRf{\kappa\LRf{z}} & if~ k\leq d-1,~n=0~mod~4\\
0 &if~  k\leq d-1,~ n=1, 2, 3~mod~4  }\\
\end{array}
$$
Therefore
$$
\ora{CW}^{\LRt{n}pk}\LRf{X, \LRs{z}} =\LBrace{{ll}
0 & if~p\neq d\\
G{\CW}_{k-d}^{\LRt{n+k-d}}\LRf{C{\BM}^d\LRf{X_{z}}}& p=d\\
}
$$
To compute $\ora{CW}^{\LRt{n}pk}\LRf{Y, \pi^{-1}\LRs{z}}$ we can assume $X=\spec{B}$ affine,
with $z\in U$, and ${\SE}=A^r$ is trivial. Further, we can assume that $r=1$. So, $Y= X\times {\BA}^1=\spec{R}$, with $R=A\LRt{T}$. 
We have
$$
C^{\LRt{n}pk}\LRf{Y, \pi^{-1}\LRs{z}}= \bigoplus_{y\in Y^{\LRf{p}}\cap \pi^{-1}z} G{\CW}_{k-p}^{\LRt{n+k-p}}\LRf{C{\BM}^{p}\LRf{Y_y}}
$$
\bE
\item % 
For $p\neq d$
 we have $C^{\LRt{n}{p}k}\LRf{Y, \pi^{-1}z}=0$. Hence
$$
\ora{CW}^{\LRt{n}pk}\LRf{X, \LRs{z}} = \ora{CW}^{\LRt{n}pk}\LRf{Y, \pi^{-1}\LRs{z}} 
$$
\item Consider the case $p=d$.  
By (\ref{spexactLocal}), we have the exact sequence
%{\scalefont{.7}
$$
\diagram 
0 \ar[r] &G{\CW}^{\LRt{n+k-d}}_{k-d}\LRf{C{\BM}^{d}\LRf{X_z}}
\ar[dl]^{\iota} &\\
 G{\CW}_{k-d}^{\LRt{n+k-d}}\LRf{C{\BM}^{d}\LRf{R_{\tilde{z}}}}
\ar[r]^{\partial_{\tilde{z}}\qquad\qquad} & \bigoplus _{y\in Y_z^{\LRf{d+1}}}G{\CW}_{k-d-1}^{\LRt{n+k-d-1}}\LRf{C{\BM}^{d+1}\LRf{Y_y}}\ar[r] & 0\\
0\ar[r] &C^{\LRt{n}dk}\LRf{Y, \pi^{-1}z}\ar@{=}[u]\ar[r]&C^{\LRt{n}(d+1)k}\LRf{Y, \pi^{-1}z}\ar@{=}[u]&\\
\enddiagram
$$ 
%}
It follows
$$
 \ora{CW}^{\LRt{n}dk}\LRf{Y, \pi^{-1}\LRs{z}} =
G{\CW}^{\LRt{n+k-d}}_{k-d}\LRf{C{\BM}^{d}\LRf{X_z}}=\ora{CW}^{\LRt{n}k}_p\LRf{X, \LRs{z}} 
$$
\eE 
So, in all the cases, we have
$$
 \ora{CW}^{\LRt{n}dk}\LRf{Y, \pi^{-1}\LRs{z}} =
\ora{CW}^{\LRt{n}pk}\LRf{X, \LRs{z}} 
$$
\pic $\eop$ 
%%%%%%%%%%%%%%%%%%%%%%%%%%%%%%%%%%%%%%%%%%%%%%%%%%

\vspace{3mm}
The following is one of our basic result, on homotopy invariance.
\bC\label{ArKoChgp}{\rm 
Let ${\BF}$ be a field with $\frac{1}{2} \in {\BF}$ and ${\BA}^r=\spec{{\BF}\LRt{T_1, T_2, \ldots, T_r}}$. Then 
$$
\ora{CW}^{\LRt{n}pk}\LRf{{\BA}^r} = \ora{CW}^{\LRt{n}pk}\LRf{{\BF}}= \LBrace{{ll}
G{\CW}_k^{\LRt{n+k}}\LRf{{\BF}} & if~p=0\\
0 & if ~p\neq 0\\
}%\qquad \forall n, k, p
$$
}
\eC
\pf Use (\ref{VectorBunCW}, \ref{milJinh}). 
$\eop$
%%%%%%%%%%%%%%%%%%%%%%%%%%%%%%%%%%%%%%%%%%%%%%%%%%

%%%%%%%%%%%%%%%%%%%%%%%%%%%%%%%%%%%%%%%%%%%%%%%%%%

%%%%%%%%%%%%%%%%%%%%%%%%%%%%%%%%%%%%%%%%%%%%%%%%%
%%%%%%%%%%%%%%%%%%%%%%%%%%%%%%%%%%%%%%%%%%%%%%%%%
%%%%%%%%%%%%%%%%%%%%%%%%%%%%%%%%%%%%%%%%%%%%%%%%%
%%%%%%%%%%%%%%%%%%%%%%%%%%%%%%%%%%%%%%%%%%%%%%%%%
%%%%%%%%%%%%%%%%%%%%%%%%%%%%%%%%%%%%%%%%%%%%
%%%%%%%%%%%%%%%%%%%%%%%%%%%%%%%%%%%%%%%%%%%%
%%%%%%%%%%%%%%%%%%%%%%%%%%%%%%%%%%%%%%%%%%%%
%%%%%%%%%%%%%%%%%%%%%%%%%%%%%%%%%%%%%%%%%%%%
%%%%%%%%%%%%%%%%%%%%%%%%%%%%%%%%%%%%%%%%%%%%
%%%%%%%%%%%%%%%%%%%%%%%%%%%%%%%%%%%%%%%%%%%%
%%%%%%%%%%%%%%%%%%%%%%%%%%%%%%%%%%%%%%%%%%%%
\section{Chow-Witt groups of projective spaces}\label{SecPdDesc}
%%%%%%%%%%%%%%%%%%%%%%%%%%%%%%%%%%%%%%%%%%%%
%%%%%%%%%%%%%%%%%%%%%%%%%%%%%%%%%%%%%%%%%%%%
%%%%%%%%%%%%%%%%%%%%%%%%%%%%%%%%%%%%%%%%%%%%
%%%%%%%%%%%%%%%%%%%%%%%%%%%%%%%%%%%%%%%%%%%%
As usual, ${\BF}$ will denote a field, with $\frac{1}{2}\in {\BF}$ and ${\BP}^d=\proj{{\BF}\LRt{X_0, X_1, 
\ldots, X_d}}$ will denote the projective space. In this section, we compute the arrowtic Chow-Witt groups
$\ora{CW}^{\LRt{n}pk}\LRf{{\BP}, {\CO}\LRf{\l}}$ for $\l\in {\BZ}$. However, because of periodicity (\ref{Lsquare}), we only  consider $\l=0, \pm 1$.  
In this section, we will often substitute the twist ${\CO}\LRf{\l}$, by the integer $\l$ itself. So, 
  $$
  \LBrace{{l}
  C^{\LRt{n}k}_{\bul}\LRf{{\BP}^d, \l}:= C^{\LRt{n}k}_{\bul}\LRf{{\BP}^d, {\CO}\LRf{\l}}\\
  \ora{CW}^{\LRt{n}pk}\LRf{{\BP}^d, \l}:= \ora{CW}^{\LRt{n}pk}\LRf{{\BP}^d, {\CO}\LRf{\l}} \\
  G{\CW}^{\LRt{n}}_p\LRf{{\BP}^d, \l}=  G{\CW}_p^{\LRt{n}}\LRf{{\BP}^d, {\CO}\LRf{\l}}\\
  }
  \quad {\rm and, ~likewise}.
  $$
As usual, if $\l=0$, we ignore the second coordinate. 
The following  is an immediate consequence of (\ref{subsech}). %ZXYcontexdt}).
\bL\label{CorPnnn}{\rm 
Let ${\CL}$ be line bundle on ${\BP}^d$. 
Write ${\BP}^{d-k}=V\LRf{X_{d-k+1}, \ldots, X_d}\subseteq {\BP}^d$. Then the following
%there is a commutative diagram of duality preserving functors:
$$
\diagram 
&{\SV}\LRf{{\BP}^{d-k}, {\CL}_{|{\BP}^{d-k}}} \ar@/^/[dr]\ar@/_/[dl]& \\
C{\BM}^{{\BP}^{d-k}}\LRf{{\BP}^{d-1}, {\CL}_{|{\BP}^{d-1}}(k-1)}\ar[rr]^{\qquad X_d=0}_{\varphi} && C{\BM}^{{\BP}^{d-k}}\LRf{{\BP}^d, {\CL}(k)}\\
\enddiagram \\
$$ 
is a commutative diagram of exact categories with duality. 
}
\eL
\pf Follows from  (\ref{subsech}). $\eop$ % 

%%%%%%%%%%
\vspace{3mm}
The following is a particular case of the exact sequence (\ref{CBulSeqEqn}), for this section.
\bL\label{PdkoPromecx}{\rm 
Consider ${\BA}^d=\LRf{X_d\neq 0}\subseteq {\BP}^d$.
Then the following is an exact sequence of complexes
\begin{equation}\label{CPnuSeqEqn}
\diagram 
0 \ar[r] & C^{\LRt{n}(k-1)}_{\bul}\LRf{{\BP}^{d-1}, \l-1} \ar[r] & C^{\LRt{n}k}_{\bul}\LRf{{\BP}^d, \l} \ar[r] &  C^{\LRt{n}k}_{\bul}\LRf{{\BA}^{d}, {\CO}\LRf{\l}_{|{\BA}^{d}}} \ar[r] & 0\\
\enddiagram
\end{equation}
}
\eL
\pf Consider $Y={\BP}^{d-1}=\LRf{X_d=0}\subseteq {\BP}^d$. The rest follows from Lemma \ref{refhKwwdm}. 
$\eop$ 

\vspace{3mm}
Further, the following is a particular case of  (\ref{LongSiSediSp}).
%%%%%%%%%%%%%%%%%%%%%%%%%%%%%%%%%%%%%%%%%%%%
\bT\label{longsaPnre}{\rm
Let ${\BF}$ be a field, with $\frac{1}{2}\in {\BF}$. Let ${\BP}^d:={\BP}^d_{{\BF}}$ denote the projective
$d$ space and $\l\in {\BZ}$.  Then  
\bE
 \item\label{oneOfPd} First,    $\ora{CW}^{\LRt{n}pk}\LRf{{\BP}^d, \l} =0$ unless $p=0, 1, \ldots, d$ 
 \item\label{twoOfPd}   The maps
 \begin{equation}\label{ortherTerhInd}
   \diagram
  \ora{CW}^{\LRt{n}(p-1)(k-1)}\LRf{{\BP}^{d-1}, \l-1} \ar[r]^{\qquad\quad \sim} &
 \ora{CW}^{\LRt{n}pk}\LRf{{\BP}^d, \l} \\
 \enddiagram
 \end{equation}
 are isomorphisms $ \forall p= 2, 3, \ldots, d$. 
 \item\label{for4Exact}
 We have a four term exact sequence 
  \begin{equation}\label{khubJoruri}
  \diagram
 0 \ar[r] & %\ora{CH}^{\LRt{n}0(k-1)}\LRf{{\BP}^{r-1}, {\CO}{(\l-1)}} \ar[r] &
 \ora{CW}^{\LRt{n}0k}\LRf{{\BP}^d, \l} \ar[r] &  G{\CW}^{\LRt{n+k}}_k\LRf{{\BF}} \ar[lld]_{\delta_k} \\
  \ora{CW}^{\LRt{n}0(k-1)}\LRf{{\BP}^{d-1}, \l-1} \ar[r] &
 \ora{CW}^{\LRt{n}1k}\LRf{{\BP}^d, \l} \ar[r] & 0 &\\
 \enddiagram
  \end{equation} 
 
 \eE

 }
\eT

 %
% \TCP
\pf Proof of (\ref{oneOfPd}) is obvious. %With $Y=\LRf(X_d=0)$ and $U={\BP}^d-Y\cong {\BA}^{d}$,
 From  (\ref{PdkoPromecx}), we obtain the following exact sequence:
  \begin{equation}\label{LonPanTaEqn1011}
 \diagram
 0 \ar[r] %& \ora{CH}^{\LRt{n}(p-1)(k-1)}\LRf{{\BP}^{r-1}, {\CO}{(\l-1)}} \ar[r] 
 &
 \ora{CW}^{\LRt{n}0k}\LRf{{\BP}^d, \l} \ar[r] &  \ora{CW}^{\LRt{n}0k}\LRf{{\BA}^d} \ar[dll] \\
  \ora{CW}^{\LRt{n}0(k-1)}\LRf{{\BP}^{d-1}, \l-1} \ar[r] &
 \ora{CW}^{\LRt{n}1k}\LRf{{\BP}^d, \l} \ar[r] &  \ora{CW}^{\LRt{n}0k}\LRf{{\BA}^d} \ar[dll] \\
   \ora{CW}^{\LRt{n}1(k-1)}\LRf{{\BP}^{d-1}, \l-1} \ar[r] &
 \ora{CW}^{\LRt{n}2k}\LRf{{\BP}^r, \l} \ar[r] & \cdots \\% \ora{CH}^{\LRt{n}(p+1)k}_{-1}\LRf{{\BA}^r} \ar[r] &\cdots\\
 \enddiagram
 \end{equation} 
 By homotopy invariance (\ref{ArKoChgp}), on the $3^{rd}$ column, we have
 $$
 \ora{CW}^{\LRt{n}pk}\LRf{{\BA}^d}=  \ora{CW}^{\LRt{n}pk}\LRf{\spec{{\BF}}}=
 \LBrace{{ll}
 G{\CW}_k^{\LRt{n+k}}\LRf{{\BF}} & if ~p=0\\
 0 & {\rm otherwise}\\
 }
 $$
  Now, (\ref{twoOfPd}) and (\ref{for4Exact}) follow immediately. \pic $\eop$

%  
 %Consequently, from Proposition \ref{LongSiSediSp}, we have the following:

 %%%%%%%%%%%%%%%%%%%%%%%%%%%%%%%%%%%%%%%%%%%%
 \vspace{3mm}
Now we proceed to compute Chow-Witt groups of ${\BP}^d$. 
%Because of periodicity (\ref{Lsquare}), we only consider two dualities. 
The case of $d=0$ is immediate, as follows.
\begin{equation}\label{disze0case}
\ora{CW}^{\LRt{n}pk}\LRf{{\BP}^0}= \ora{CW}^{\LRt{n}pk}\LRf{\spec{{\BF}}}=
 \LBrace{{ll}
 GW^{\LRt{n+k}}_k\LRf{{\BF}} & if ~p=0\\
 0 & {\rm otherwise}\\
 }
\end{equation}
%%%

\vspace{3mm}
%We start with the following case of $p=0$. 
%\subsection{The degree zero lemma} 
The case of degree $p=0$ and the twist $\l=0$ proves to be a key ingredient, for the subsequent 
arguments. We refer to it zero-zero lemma, as follows. 
\bL[Zero-zero]\label{degreeZero}{\rm 
Let ${\BF}$ be a field, with $\frac{1}{2}\in {\BF}$. Let ${\BP}^d=\proj{{\BF}\LRt{X_0, X_1, \ldots, X_d}}$. 
Then 
$$
\ora{CW}^{\LRt{n}0k}\LRf{{\BP}^d, 0} = G{\CW}^{\LRt{n+k}}_k\LRf{{\BF}} \qquad \forall k, n
$$ 
}
\eL
\pf % 
Let  $U={\BA}^d=(X_d\neq 0)\subseteq {\BP}^d$.  We have 
${\BA}^d=\spec{{\BF}\LRt{T_0, T_1, \ldots, T_{d-1}}}$, with $T_i=\frac{X_i}{X_d}$. 
Let ${\BF}\LRf{{\BP}^d}\cong {\BF}\LRf{{\BA}^d}$ denote the respective field of rational functions. 
Consider the commutative diagram of maps (with $\l=0$ for all the dualities)
 \begin{equation}\label{resketch4P1}
  \diagram 
\ker\LRf{\partial_{{\BP}}}\ar@{^(->}[d] \ar@{^(->}[r] &G{\CW}_k^{\LRt{n+k}}\LRf{{\BF}}\ar@{^(->}[d]^{\pi^*}\\
  G{\CW}_k^{\LRt{n+k}}\LRf{{\BF}\LRf{{\BP}^d}} \ar[r]^{\iota}_{\sim} \ar[d]_{\partial_{{\BP}}}& G{\CW}_k^{\LRt{n+k}}\LRf{{\BF}\LRf{{\BA}^d}}\ar[d]^{\partial_{{\BA}}} \\
 \bigoplus_{x\in \LRf{{\BP}^d}^{\LRf{1}}}   G{\CW}_{k-1}^{\LRt{n+k-1}}\LRf{C{\BM}^1\LRf{{\BP}^d_x}} \ar[r]_{\q} & 
  \bigoplus_{x\in \LRf{{\BA}^d}^{\LRf{1}}} G{\CW}_{k-1}^{\LRt{n+k-1}}\LRf{C{\BM}^1\LRf{{\BA}^d_x}} \\
  \enddiagram
  \end{equation}
  The vertical maps $\partial_{{\BP}}$ and $\partial_{{\BA}}$ denote the tail of the respective Gersten ${G}{\CW}$-complexes
  \cite[pp 473]{M23}. By definition,
  $\ora{CW}^{\LRt{n}0k}\LRf{{\BP}^d, 0}=\ker\LRf{\partial_{{\BP}}}$. 
  By homotopy invariance (\ref{ArKoChgp}), we have 
  $\ker\LRf{\partial_{{\BA}}}=\ora{CW}^{\LRt{n}0k}\LRf{{\BA}^d}= G{\CW}_k^{\LRt{n+k}}\LRf{{\BF}}$. Since middle horizontal map is an isomorphism, the top arrow is also injective. 
  %\TCP{ Alternately, }
% 
  Consider the commutative diagram %({\it we ignore the duality related coordinate ${\CO}\LRf{0}$})
  $$
  \diagram 
  &G{\CW}^{\LRt{n+k}}_k\LRf{{\BF}}\ar[d]^{\iota^{-1}\pi_*}\ar@/_/[dl]_{\theta}&\\
  G{\CW}_k^{\LRt{n+k}}\LRf{C{\BM}^0\LRf{{\BP}^d}} \ar[r] &  
  G{\CW}_k^{\LRt{n+k}}\LRf{{\BF}\LRf{{\BP}^d}} \ar[r]^{\partial_{\BP}\qquad\qquad} & 
   \bigoplus_{x\in \LRf{{\BP}^d}^{\LRf{1}}}   G{\CW}_{k-1}^{\LRt{n+k-1}}\LRf{C{\BM}^1\LRf{{\BP}^d_x}} \\
  \enddiagram
  $$
  where $\theta$ is the pullback, which exists because duality corresponds to $\l=0$.
  The horizontal line is exact \cite[pp 473]{M23}.
  The diagram shows that the vertical natural pullback map factors through $\theta$. Therefore
  $$
  G{\CW}^{\LRt{n+k}}_k\LRf{{\BF}} \subseteq   \ker\LRf{\partial_{{\BP}}}. \quad {\rm Hence}
  \quad   G{\CW}^{\LRt{n+k}}_k\LRf{{\BF}} =   \ker\LRf{\partial_{{\BP}}}=
  \ora{CW}^{\LRt{n}0k}\LRf{{\BP}^d, 0}.
  $$
\pic $\eop$

\vspace{3mm}
Combining the zero-zero lemma (\ref{degreeZero}), the four term exact sequence (\ref{khubJoruri}) takes the following form.
%{\scalefont{.8}
\begin{equation}\label{zeroHero10}
  \diagram
 0 \ar[r] & 
 \ora{CW}^{\LRt{n}0k}\LRf{{\BP}^d, 1} \ar[r] &  G{\CW}^{\LRt{n+k}}_k\LRf{{\BF}} \ar[r]^{\b^{\LRt{n+k}d}_{k}} & 
G{\CW}^{\LRt{n+k-1}}_{k-1}\LRf{{\BF}} \ar[r] &
 \ora{CW}^{\LRt{n}1k}\LRf{{\BP}^d, 1} \ar[r] & 0 &\\
 \enddiagram
  \end{equation} 
 % }
  where we relabeled $\b^{\LRt{n+k}d}_{k}:=\delta_k$. 
  
  %%%%%%%%%%%%%%%%%%%%%%%%%%%%%%%%%%%%%%%%%%%%%%
  %%%%%%%%%%%%%%%%%%%%%%%%%%%%%%%%%%%%%%%%%%%%%%
  %%%%%%%%%%%%%%%%%%%%%%%%%%%%%%%%%%%%%%%%%%%%%%
  %%%%%%%%%%%%%%%%%%%%%%%%%%%%%%%%%%%%%%%%%%%%%%
  %%%%%%%%%%%%%%%%%%%%%%%%%%%%%%%%%%%%%%%%%%%%%%
  %%%%%%%%%%%%%%%%%%%%%%%%%%%%%%%%%%%%%%%%%%%%%%
  
  \subsection{Consistency} 
  The exact sequence (\ref{zeroHero10}) will be instrumental for the purpose descriptions of 
  the  Chow-Wittt groups $\ora{CW}^{\LRt{n}pk}\LRf{{\BP}^d, \l}$. 
  We proceed to establish that all these maps 
  $\b^{\LRt{n+k}d}_{k}$  are independent of the dimension $d$. 
  %%%%
  \bN\label{notaVott}{\rm 
  Let ${\BF}$ be a field, with $\frac{1}{2}\in {\BF}$. Let 
  ${\BP}^d=\proj{{\BF}\LRt{X_0, X_1, \ldots, X_d}}$ and $m, k\in {\BZ}$. 
  The field of rational functions on ${\BP}^d$ is denoted by 
  ${\BF}\LRf{{\BP}^d}$ ({\it likewise, for ${\BA}^d$ or any integral scheme $X$}). We consider ${\BA}^d=\LRf{X_d\neq 0}\subseteq {\BP}^d$. 
     }
  \eN
 % \TCP{XXX}
By reindexing $m:=n+k$,  we rewrite  the exact sequence (\ref{CPnuSeqEqn}), as follows
\begin{equation}\label{ReinIndSeqEqn}
\diagram 
0 \ar[r] & C^{\LRt{m-k}(k-1)}_{\bul}\LRf{{\BP}^{d-1},  \l-1} \ar[r] & C^{\LRt{m-k}k}_{\bul}\LRf{{\BP}^d, \l} \ar[r] &  C^{\LRt{m-k}k}_{\bul}\LRf{{\BA}^{d}} \ar[r] & 0\\
\enddiagram
\end{equation}
With $\l=1$ the  four term exact sequence (\ref{khubJoruri}) takes the following form:
\begin{equation}\label{minx4Term}
\diagram 
0  \ar[r] & \ora{CW}^{\LRt{m-k}0k}\LRf{{\BP}^d, 1} \ar[r] &  \ora{CW}^{\LRt{m-k}0k}\LRf{{\BA}^d}\ar[dll] \\
\ora{CW}^{\LRt{m-k}0k-1}\LRf{{\BP}^{d-1}, 0} \ar[r] & \ora{CW}^{\LRt{m-k}1k}\LRf{{\BP}^d, 1} \ar[r]&0\\
\enddiagram
\end{equation}
By Zero-zero lemma (\ref{degreeZero})  and homotopy invariance, we have the exact sequence 
\begin{equation}\label{minx42ndTerm}
\diagram 
0  \ar[r] & \ora{CW}^{\LRt{m-k}0k}\LRf{{\BP}^d, 1} \ar[r] & G{\CW}^{\LRt{m}}_k\LRf{{\BF}}\ar[r]^{\b^{\LRt{m} d}_k} &
G{\CW}^{\LRt{m-1}}_{k-1}\LRf{{\BF}} %, \TCP{0}} 
\ar[r] & \ora{CW}^{\LRt{m-k}1k}\LRf{{\BP}^d, 1} \ar[r]&0\\
\enddiagram
\end{equation}
 The construction of the connecting homomorphism $\b^{\LRt{m}, d}_k$ is presented by the diagram: 
\begin{equation}\label{VottNtaDef}
\diagram 
& G{\CW}^{m}_k\LRf{{\BF}}\ar@{_(->}[d]^{\pi}  \ar@/_/@{-->}[ld]_{\b^{\LRt{m}d}_k}
 \\
 G{\CW}^{m-1}_{k-1}\LRf{{\BF}}\ar@{^(->}[d]_{\pi}&G{\CW}^{m}_k\LRf{{\BF}\LRf{{\BP}^d}}% , \TCP{1}} 
   \ar[d]^{\partial_{X_d}}&\\
G{\CW}^{m-1}_{k-1}\LRf{{\BF}\LRf{{\BP}^{d-1}}}%, \TCP{0}} %
\ar[r]_{\D_d\quad}  
&G{\CW}_{k-1}^{\LRt{m-1}}\LRf{C{\BM}^1\LRf{{\BP}_{\LRf{X_d}}^d}}&
  \\
\enddiagram
\end{equation} 
The vertical map $\partial_{X_d}$ is the residue map, corresponding to the point $(X_d)\in \LRf{{\BP}^d}^{\LRf{1}}$, and 
two $\pi$'s denote the pullback maps. The map $\D_d$ is the d\'{e}vissage isomorphism. 
The map $ b_k^{\LRt{m}d}$ is constructed by following the usual method of constructing connecting homomorphisms. 
 
We prove that the maps $\b^{\LRt{m}d}_k$ are independent of $\dim {\BP}^d$. 
\bL\label{Agreebott}{\rm 
Consider the notations from (\ref{VottNtaDef}). %and diagram (\ref{VottNtaDef}),
 Then
$\b^{\LRt{m}1}_k:=\b^{\LRt{m}d}_k$ for all $d\geq 1$. % are independent of $d$. 
}
\eL
\pf% \TCP{It follows, 
%$$
%\pi\b^{\LRt{m}}_k= \pi\LRf{\D_1^{-1}\partial_1\pi}= \LRf{\D^{-1}\iota^*}\partial_1\pi
%=\D^{-1}\LRf{\partial \iota^*}\pi= {\D}^{-1}\partial \pi
%$$
%Therefore, $\b^{\LRt{m}}_k=\b^{\LRt{m}d}_k$.
%}
%\\ \TCP{Seems working.} For the general case, 
We consider ${\BP}^1=\proj{{\BF}\LRt{X_0, X_d}}\subseteq {\BP}^d$. Recall the map $\b^{\LRt{m}1}_k$,
is given by the diagram:
$$
\diagram
&G{\CW}^{\LRt{m}}_k\LRf{{\BF}}\ar@/^/@{_(->}[dr]^{\pi}\ar@{-->}[d]_{\p_1}\ar@/_/[ddl]_{\b^{\LRt{m}1}_k}&\\
&G{\CW}^{\LRt{m}}_k\LRf{{\BF}\LRf{{\BP}^1}, 1} \ar[r]_{\sim} \ar[d]_{\partial_{X_d}} & 
G{\CW}^{\LRt{m}}_k\LRf{{\BF}\LRf{{\BA}^1}}\\
G{\CW}^{\LRt{m-1}}_{k-1}\LRf{{\BF}}\ar[r]_{\D_1\qquad}^{\sim\qquad}&
G{\CW}^{\LRt{m-1}}_{k-1}C{\BM}^1\LRf{{\BP}^1_{\LRf{X_d}}, 1}
& \\
\enddiagram 
%\quad \b=\D_1^{-1}\partial_{X_1} \p_1
$$ 
So, $\b^{\LRt{m}1}_k= \D_1^{-1}\partial{X_d}\p_1$, where the broken arrow $\p_1$ is given by composition, via 
the pullback $\pi$.
We 
construct the diagram 
$$
\diagram 
&&G{\CW}_k^{\LRt{m}}{\BF} \ar@/^/[rd]^{\pi} \ar@/_/@{-->}[ddll]_{\b^{\LRt{m}1}_k}^{\b^{\LRt{m}d}_k}& \\
&&G{\CW}_k^{\LRt{m}}\LRf{{\BF}{\BP}^1, 1}\ar@/_/[dl]_{\partial_{\infty}}\ar[r]^{\sim}\ar[d]^{\otimes}&G{\CW}_k^{\LRt{m}}\LRf{{\BF}{\BA}^1, 1}\ar[d]\\
 G{\CW}_{k-1}^{\LRt{m-1}}\LRf{{\BF}} \ar[d]_{\otimes}\ar[r]^{\D_1} 
 &G{\CW}_{k-1}^{\LRt{m-1}}\LRf{{\BP}^1_{\LRf{X_d}}}\ar@{-->}[d]^{\otimes}
& G{\CW}_k^{\LRt{m}}\LRf{{\BF}{\BP}^d, 1} 
\ar@/^/[dl]^{\partial_{\infty}}\ar[r]^{\sim} \ar[d] & G{\CW}_k^{\LRt{m}}\LRf{{\BF}{\BA}^d, 1}\ar[d]\\
G{\CW}_{k-1}^{\LRt{m-1}}\LRf{{\BF}{\BP}^1, 0}\ar[r]_{\D_d}&G{\CW}_{k-1}^{\LRt{m-1}}\LRf{{\BP}^d_{\LRf{X_d}}}\ar[r] &\bigoplus G{\CW}_{k-1}^{\LRt{m-1}}\LRf{C{\BM}^1{\BP}^d_x, 1} \ar[r]&
\bigoplus G{\CW}_{k-1}^{\LRt{m-1}}\LRf{C{\BM}^1{\BA}^d_x, 1} \\
\enddiagram
$$
The d\'{e}vissage rectangle commutes $\otimes \D_1= \D_d\otimes$, because of the simple commutativity at the dg category level.  To check the commutativity of the  $\partial_{\infty}$-rectangle, write $T_r=\frac{X_r}{X_0}$. The 
diagram
$$
\diagram
{\bf dg}^{\LRt{m-1}} C{\BM}^1{\BF}\LRt{T_d}_{\LRf{T_2}} \ar[r] \ar[d]& {\bf dg}^{\LRt{m}} C{\BM}^0{\BF}\LRt{T_d}_{\LRf{T_d}} \ar[r]\ar[d] & {\bf dg}^{\LRt{m}} C{\BM}^0{\BF}\LRf{T_d} \ar[d]\\
{\bf dg}^{\LRt{m-1}} C{\BM}^1{\BF}\LRt{T_1, \cdots, T_d}_{\LRf{T_2}} \ar[r] & {\bf dg}^{\LRt{m}} C{\BM}^0{\BF}\LRt{T_1, 
\cdots, T_d}_{\LRf{T_d}} \ar[r] & {\bf dg}^{\LRt{m}} C{\BM}^0{\BF}\LRf{T_1, \cdots, T_d}\\
\enddiagram
$$
commutes. 
This induces a map of corresponding homotopy fibration of $G{\CW}$-spectra. Therefore, the $\partial_{\infty}$ rectangle commutes. 
It follows from 
the methods of constructing such connecting homomorphism
$$
{\b^{\LRt{m}1}_k}={\b^{\LRt{m}d}_k}
$$
 \pic $\eop$

\vspace{3mm} 
So, henceforth, we  write $\b^{\LRt{m}}_k=\b^{\LRt{m}d}_k$. 
%%%%%%%%%%%%%%%%%%%%%%%%%%%%%%%%%%%%%%%%%%%%
  %%%%%%%%%%%%%%%%%%%%%%%%%%%%%%%%%%%%%%%%%%%%%%%%%%
  %%%%%%%%%%%%%%%%%%%%%%%%%%%%%%%%%%%%%%%%%%%%%%%%%%
  %%%%%%%%%%%%%%%%%%%%%%%%%%%%%%%%%%%%%%%%%%%%%%%%%%
  %%%%%%%%%%%%%%%%%%%%%%%%%%%%%%%%%%%%%%%%%%%%%%%%%%
  %%%%%%%%%%%%%%%%%%%%%%%%%%%%%%%%%%%%%%%%%%%%%%%%%%

 % \noindent\TCP{endBott}
 % \subsection{Agreement of connecting homomorphisms} 
 % Given $m, k\in {\BZ}$ most canonical connecting homomorphism 
  %$\diagram 
 % G{\CW}^{\LRt{m}}_k\LRf{{\BF}} \ar[r] & 
%G{\CW}^{\LRt{m-1}}_{k-1}\LRf{{\BF}}\\
%\enddiagram$
%is given by the Algebraic Bott sequence \cite[pp 1874]{S17}. In this section, we shall establish that all these maps 
 %$\delta^{\LRt{n+k}}_{d,k}$ agrees with the one we get from the Bott sequence. \\
% (\TCP{We postpone talk about Bott sequence till later}.)  
  %%%%%%%%%%%%%%%%%%%%%%%%%%%%%%%%%%%%%%%%%%%%%%%%%%
  %%%%%%%%%%%%%%%%%%%%%%%%%%%%%%%%%%%%%%%%%%%%%%%%%%
  %%%%%%%%%%%%%%%%%%%%%%%%%%%%%%%%%%%%%%%%%%%%%%%%%%
  %%%%%%%%%%%%%%%%%%%%%%%%%%%%%%%%%%%%%%%%%%%%%%%%%%
  %%%%%%%%%%%%%%%%%%%%%%%%%%%%%%%%%%%%%%%%%%%%%%%%%%
  %%%%%%%%%%%%%%%%%%%%%%%%%%%%%%%%%%%%%%%%%%%%%%%%%%
  %%%%%%%%%%%%%%%%%%%%%%%%%%%%%%%%%%%%%%%%%%%%%%%%%%

%%%%%%%%%%%
\subsection{Final Computations} %Case of ${\BP}^1$}
First, we establish the Chow-Witt groups of ${\BP}^1$.
\bL\label{POneAllOn}{\rm 
Let ${\BF}$ be a field, with $\frac{1}{2}\in {\BF}$. Let ${\BP}^1=\proj{{\BF}\LRt{X_0, X_1}}$. 
Then 
$$
\ora{CW}^{\LRt{n}pk}\LRf{{\BP}^1, 0} = 
\LBrace{{ll}
G{\CW}^{\LRt{n+k}}_k\LRf{{\BF}}&p=0 \\
G{\CW}^{\LRt{n+k-1}}_{k-1}\LRf{{\BF}}&p=1 \\
0 & {\rm otherwise}\\
}
$$ 
}
\eL
\pf The case $p=0$ follows from (\ref{degreeZero}). We give a proof for  $p=1$. Re-sketch the  complex
(\ref{CPnuSeqEqn}), with  $d=1$.
$$
  \diagram 
&\ker\LRf{\partial_{{\BP}}}\ar@{^(->}[d] \ar@{^(->}[r] &G{\CW}_k^{\LRt{n+k}}\LRf{{\BF}}\ar@{^(->}[d]^{\pi^*}\\
0\ar[d]&  G{\CW}_k^{\LRt{n+k}}\LRf{{\BF}\LRf{{\BP}^1}} \ar[r]^{\iota}_{\sim} \ar[d]_{\partial_{{\BP}}}& G{\CW}_k^{\LRt{n+k}}\LRf{{\BF}\LRf{{\BA}^1}}\ar[d]^{\partial_{{\BA}}} \\
G{\CW}^{\LRt{n+k-1}}_{k-1}\LRf{{\BF}}\ar[r]\ar[d]& \bigoplus_{x\in \LRf{{\BP}^1}^{\LRf{1}}}   G{\CW}_{k-1}^{\LRt{n+k-1}}\LRf{C{\BM}^1\LRf{{\BP}^1_x}} \ar[r]_{\q} \ar[d]& 
  \bigoplus_{x\in \LRf{{\BA}^1}^{\LRf{1}}} G{\CW}_{k-1}^{\LRt{n+k-1}}\LRf{C{\BM}^1\LRf{{\BA}^1_x}} \ar[d]\\
0  &0&0\\
  \enddiagram
$$
 The right vertical map is also exact, in the sense of Theorem \ref{arrowticExactAll} ({\it meaning, $\partial_{{\BA}}$ is surjective, as well}).  The map $\q$ is split surjective,
with 
$$
\ker\LRf{\q} = G{\CW}_{k-1}^{\LRt{n+k-1}}\LRf{C{\BM}^1\LRf{{\BP}^d_{\infty}}, 0} \cong G{\CW}^{\LRt{n+k-1}}_{k-1}\LRf{{\BF}}
$$
where $\infty:=\LRf{X_1}\in {\BP}^1$, and the latter isomorphism is given by d\'{e}vissage (\ref{regLocal}).
It is obvious now, $\ora{CW}^{\LRt{n}1k}\LRf{{\BP}^1, 0}=G{\CW}^{\LRt{n+k-1}}_{k-1}\LRf{{\BF}}$.
%({\it one can use Snake Lemma, as well}). 
\pic $\eop$

%\TCP{XXXold}
%\pf The case $p=0$ follows from (\ref{degreeZero}). Only other case that needs a proof is the case of $p=1$. Re-sketch the diagram (\ref{resketch4P1}), with  $d=1$. The right vertical map is also exact, in the sense of Theorem \ref{arrowticExactAll} ({\it meaning, $\partial_{{\BA}}$ is surjective, as well}).  The map $\q$ is split surjective,with 
%$$
%\ker\LRf{\q} = G{\CW}_{k-1}^{\LRt{n+k-1}}\LRf{C{\BM}^1\LRf{{\BP}^d_{\infty}}, 0} \cong G{\CW}^{\LRt{n+k-1}}_{k-1}\LRf{{\BF}}
%$$
%where $\infty:=\LRf{X_1}\in {\BP}^1$, and the latter isomorphism is given by d\'{e}vissage (\ref{regLocal}).It is obvious now, $\ora{CW}^{\LRt{n}1k}\LRf{{\BP}^1, 0}=G{\CW}^{\LRt{n+k-1}}_{k-1}\LRf{{\BF}}$ ({\it one can use Snake Lemma, as well}). \pic $\eop$

\vspace{3mm}
Next, we compute $\ora{CW}^{\LRt{n}pk}\LRf{{\BP}^1, 1}$. 
\bL\label{oraCHP1O1}{\rm 
Let ${\BF}$ be a field, with $\frac{1}{2}\in {\BF}$. Let ${\BP}^1=\proj{{\BF}\LRt{X_0, X_1}}$. 
Then 
$$
\ora{CW}^{\LRt{n}pk}\LRf{{\BP}^1,  \pm 1} =
\LBrace{{ll}
\ker\LRf{\b_k^{n+k}}   & p=0\\
co\ker\LRf{\b_k^{n+k}}  & p=1 \\
0 & {\rm otherwise} \\
}
$$
}
\eL
\pf In this proof, we follow the diagram (\ref{zeroHero10}), %(\ref{bottAndBF}), 
with $n:=n+k$. It is immediate,
$$
\LBrace{{l}
\ora{CW}^{\LRt{n}0k}\LRf{{\BP}^1, 1} %\ker\LRf{\partial_k^{\LRt{n+k}}}
= \ker\LRf{\b_k^{\LRt{n+k}}}\\
\ora{CW}^{\LRt{n}1k}\LRf{{\BP}^1, 1} %\ker\LRf{\partial_k^{\LRt{n+k}}}
= co\ker\LRf{\b_k^{\LRt{n+k}}}\\
}
$$
Since $\dim {\BP}^1=1$, the proof is complete. $\eop$

%%%%%%%%%%%%%%%%%%%
%\subsection{Case of ${\BP}^2$}

\vspace{3mm}
Next, we compute the Chow-Witt groups $\ora{CW}^{\LRt{n}pk}\LRf{{\BP}^2, \l}$, as in the following lemma.
\bL\label{P2e00done}{\rm 
Let ${\BF}$ be a field, with $\frac{1}{2}\in {\BF}$. Let ${\BP}^2=\proj{{\BF}\LRt{X_0, X_1, X_2}}$. Then 
 $$
 \ora{CW}^{\LRt{n}pk}\LRf{{\BP}^2, 0}=  
 \LBrace{{ll}
GW^{\LRt{n+k}}_k\LRf{{\BF}} & p=0\\
%\ora{CW}^{\LRt{n}0\LRf{k-1}}\LRf{{\BP}^1,-1}=\
ker\LRf{\b_{k-1}^{\LRt{n+k-1}}}  & p=1\\
%\ora{CW}^{\LRt{n}1\LRf{k-1}}\LRf{{\BP}^1, -1}=
co\ker\LRf{\b_{k-1}^{\LRt{n+k-1}}} & p=2\\
  0 & {\rm otherwise}\\
 }
 $$
 Further,
 $$
 \ora{CW}^{\LRt{n}pk}\LRf{{\BP}^2, \pm 1} =
 \LBrace{{ll}
\ker\LRf{\b^{\LRt{n+k}}_k} &p=0\\
 co\ker\LRf{\b^{\LRt{n+k}}_k}& p=1\\
GW_{k-2}^{\LRt{n+k-2}}\LRf{{\BF}}& p=2\\
  0 & {\rm otherwise}\\
 }
 $$
 }
 \eL 
 \pf First, we compute $ \ora{CW}_p^{\LRt{n}k}\LRf{{\BP}^2, 0}$. The case $p=0$  follows 
 directly from Lemma \ref{degreeZero}. By (\ref{khubJoruri}), %\ref{for4Exact}, \ref{CBulSeqEqn}), 
 we have a four term exact sequence:
 $$
   \diagram
 %0 \ar[r] &  
 \ora{CW}^{\LRt{n}0k}\LRf{{\BP}^2, 0} \ar@{^(->}[r] &  GW^{\LRt{n+k}}_k\LRf{{\BF}} \ar[r] &
  \ora{CW}^{\LRt{n}0\LRf{k-1}}\LRf{{\BP}^{1}, -1} \ar@{->>}[r] &
 \ora{CW}^{\LRt{n}1k}\LRf{{\BP}^1, 0} \\
 \enddiagram
 $$
% \TCP{ZZZ}
 First injective arrow is an isomorphism, by the case $p=0$. Therefore, the $2^{nd}$ surjective arrow
 is also an isomorphism. Consequently,
 $p=1$ case follows form Lemma \ref{oraCHP1O1}. The case of $p=2$ follows from 
 Lemma \ref{oraCHP1O1}, in conjunction with the isomorphisms (\ref{ortherTerhInd}). For  $p\neq 0, 1, 2$
 it follows from the dimension consideration. 
 
 Now, we establish the formula for $\ora{CW}^{\LRt{n}pk}\LRf{{\BP}^2, \pm 1}$. It follows from 
 (\ref{ortherTerhInd}) that
 $$
  \diagram
  \ora{CW}^{\LRt{n}(p-1)(k-1)}\LRf{{\BP}^{1}, 0} \ar[r]^{\qquad\quad \sim\qquad} &
 \ora{CW}^{\LRt{n}pk}\LRf{{\BP}^2, 1} \\
 \enddiagram
 \quad \forall p\neq 0, 1
 $$
 Therefore,
 $$
  \ora{CW}^{\LRt{n}pk}\LRf{{\BP}^2, 1}= \ora{CW}^{\LRt{n}(p-1)(k-1)}\LRf{{\BP}^{1}, 0}
  =\LBrace{{ll}
  G{\CW}^{\LRt{n+k-2}}_{k-2}\LRf{{\BF}}& p=2\\
  0 & p\neq 0, 1, 2\\
  }
 $$
 To deal with the cases of $p=0, 1$ refer to the four term exact sequence (\ref{khubJoruri}), with $\l=1$. We have the diagram of exact sequence:
$$
  \diagram
 % &GW^{\LRt{n+k}}_k\LRf{{\BF}} \ar[r]^{\d_k}\ar@{=}[d]
 && GW^{\LRt{n+k-1}}_{k-1}\LRf{{\BF}}\ar[d]_{\wr}&  \\
 \ora{CW}^{\LRt{n}0(k-1)}\LRf{{\BP}^{2}, 1}  \ar@{^(->}[r]&   GW^{\LRt{n+k}}_k\LRf{{\BF}}  \ar[r]_{\delta_k\qquad}\ar@/^/[ru]^{\b^{\LRt{n+k}}_k}& \ora{CW}^{\LRt{n}0(k-1)}\LRf{{\BP}^{1}, 0} \ar@{->>}[r]&   \ora{CW}^{\LRt{n}1(k-1)}\LRf{{\BP}^{2}, 1}  \\
 %&G{\CW}_k^{n+k}\LRf{{\BA}^2}\ar@{=}[u]&&\\
 \enddiagram
 $$
The formula follows. \pic $\eop$

%%%%%%%%%%%%%%%%%%%%%%%%%%%%%%%%%%%%%%%%%%%%
\bT\label{mainPSpace}{\rm 
Let ${\BF}$ be a field, with $\frac{1}{2}\in {\BF}$. Let $k, n\in {\BZ}$.
Then 
$$
\ora{CW}^{\LRt{n}pk}\LRf{{\BP}^d, 0} =
\LBrace{{ll}
G{\CW}^{n+k}_k\LRf{{\BF}} & if ~p=0\\
\ora{CW}^{\LRt{n}(p-1)(k-1)}\LRf{{\BP}^{d-1}, -1} & p=1, \ldots, d\\
0 & {\rm otherwise}\\
} 
$$
If $d$ is even, then 
$$
\ora{CW}^{\LRt{n}pk}\LRf{{\BP}^d, 1} =
\LBrace{{ll}
\ker\LRf{\b_k^{\LRt{n+k}}} & if ~p=0\\
co\ker\LRf{\b_k^{\LRt{n+k}}} &  if ~p=1\\
\ker\LRf{\b_{k-2}^{\LRt{n+k-2}}} & if ~p=2\\
co\ker\LRf{\b_{k-2}^{\LRt{n+k-2}}} & if ~ p=3\\
%\ora{CH}^{\LRt{n}\TCP{k-4}}_p\LRf{{\BP}^{\TCP{r-4}}, {\CO}\LRf{\TCP{1}}} &p=k-4, \ldots, k-r\\
\cdots\\
%\ora{CH}^{\LRt{n}\TCP{k-(r-2)}}_p\LRf{{\BP}^{\TCP{2}}, {\CO}\LRf{\TCP{1}}} &p=k-(r-2), k-(r-1), k-r\\
\ker\LRf{\b_{k-(d-2)}^{\LRt{n+k-(d-2)}}} & if ~p=d-2\\
co\ker\LRf{\b_{k-(d-2)}^{\LRt{n+k-(d-2)}}} & if ~p= d-1\\
G{\CW}_{k-d}^{\LRt{n+k-d}}\LRf{{\BF}}& if ~p= d\\
0 & {\rm otherwise}\\
} 
$$
If $d$ is odd,  then
$$
\ora{CW}^{\LRt{n}pk}\LRf{{\BP}^d, 1} =
\LBrace{{ll}
\ker\LRf{\b_k^{\LRt{n+k}}} & if ~p=0\\
co\ker\LRf{\b_k^{\LRt{n+k}}} & if ~ p=1\\
\ker\LRf{\b_{k-2}^{\LRt{n+k-2}}} & if ~p=2\\
co\ker\LRf{\b_{k-2}^{\LRt{n+k-2}}} & if ~ p=3\\
%\ora{CH}^{\LRt{n}\TCP{k-4}}_p\LRf{{\BP}^{\TCP{r-4}}, {\CO}\LRf{\TCP{1}}} &p=k-4, \ldots, k-r\\
\cdots\\
%\ora{CH}^{\LRt{n}\TCP{k-(r-3)}}_p\LRf{{\BP}^{\TCP{3}}, {\CO}\LRf{\TCP{1}}} &p=k-(r-3),\ldots, k-r\\
%\ker\LRf{\partial_{k-(r-3)}^{\LRt{n+k-(r-3)}}} &p=k-(r-3)\\
%co\ker\LRf{\partial_{k-(r-3)}^{\LRt{n+k-(r-3)}}}&p=k-(r-2)\\
\ker\LRf{\b_{k-(d-1)}^{\LRt{n+k-(d-1)}}}& if ~p=d-1\\
co\ker\LRf{\b_{k-(d-1)}^{\LRt{n+k-(d-1)}}}& if ~p= d\\
0 & {\rm otherwise}\\
} 
$$

}
\eT
\pf First, we establish the formula for $\ora{CW}^{\LRt{n}pk}\LRf{{\BP}^d, 0}$.
The case $p=0$  follows directly from the Zero-zero Lemma \ref{degreeZero}.
From Theorem \ref{longsaPnre}, we have isomorphisms
$$
 \forall p\neq 0, 1\quad % {\rm the ~maps}\quad
   \diagram
  \ora{CW}^{\LRt{n}p(k-1)}\LRf{{\BP}^{d-1},\l-1} \ar[r]^{\quad\sim} &
 \ora{CW}^{\LRt{n}pk}\LRf{{\BP}^d, \l} \quad \forall \l\\
 \enddiagram
 $$
Further, consider the four term exact sequence (\ref{khubJoruri})
 {\scalefont{.8}
 $$ 
  \diagram
  0 \ar[r] &
 \ora{CW}^{\LRt{n}0k}\LRf{{\BP}^d, \l} \ar[r] &  GW^{\LRt{n+k}}_k\LRf{{\BF}} \ar[r] 
  &\ora{CW}^{\LRt{n}0(k-1)}\LRf{{\BP}^{d-1}, \l-1} \ar[r] &
 \ora{CW}^{\LRt{n}1k}\LRf{{\BP}^r, \l} \ar[r] & 0 &\\
 \enddiagram
$$ 
 }
For $\l=0$,   the $2^{nd}$ arrow is an isomorphism, and hence so is the $3^{rd}$ arrow. This establishes the formula for  $\ora{CW}^{\LRt{n}pk}\LRf{{\BP}^d, 0}$.

Now, we deal with the formulas for  $\ora{CW}^{\LRt{n}pk}\LRf{{\BP}^d, 1}$.
With $\l=\pm 1$, the four term exact sequence takes the following form: 
 %{\scalefont{.9}
 $$ 
  \diagram
  0 \ar[r] &
 \ora{CW}^{\LRt{n}0k}\LRf{{\BP}^d,1} \ar[r] &  GW^{\LRt{n+k}}_k\LRf{{\BF}} \ar[r]^{\b_k^{\LRt{n+k}}}
  &GW^{\LRt{n+k-1}}_{k-1}\LRf{{\BF}} \ar[r]%\ora{CH}^{\LRt{n}(k-1)}_{k-1}\LRf{{\BP}^{d-1}, 0} \ar[r]\ar[d]^{\wr} 
  &
 \ora{CW}^{\LRt{n}1k}\LRf{{\BP}^d, 1} \ar[r] & 0 &\\
 \enddiagram
$$
%}
Therefore, we have
$$
\ora{CW}^{\LRt{n}pk}\LRf{{\BP}^d, 1} =
\LBrace{{ll}
\ker\LRf{\b_k^{\LRt{n+k}}} & if ~p=0\\
co\ker\LRf{\b_k^{\LRt{n+k}}} & p=1\\
\ora{CW}^{\LRt{n}(p-1)(k-1)}\LRf{{\BP}^{d-1}, 0} &p=2, \ldots, d\\
0 & {\rm otherwise}\\
} 
$$
By the formula for $\ora{CW}^{\LRt{n}pk}\LRf{{\BP}^d, 0}$, we have % Which is
$$
\ora{CW}^{\LRt{n}pk}\LRf{{\BP}^d, 1} 
=
\LBrace{{ll}
\ker\LRf{\b_k^{\LRt{n+k}}} & if ~p=0\\
co\ker\LRf{\b_k^{\LRt{n+k}}} & p=1\\
\ora{CW}^{\LRt{n}(p-2)\LRf{k-2}}\LRf{{\BP}^{d-2}, 1} &p=2, \ldots, d\\
0 & {\rm otherwise}\\
} 
$$
Assume $d$ is even. Inductively, then 
$$
\ora{CW}^{\LRt{n}pk}\LRf{{\BP}^d, 1} =
\LBrace{{ll}
\ker\LRf{\b_k^{\LRt{n+k}}} & p=0\\
co\ker\LRf{\b_k^{\LRt{n+k}}} & p=1\\
\ker\LRf{\b_{k-2}^{\LRt{n+k-2}}} & p=2\\
co\ker\LRf{\b_{k-2}^{\LRt{n+k-2}}} & p=3\\
%\ora{CW}^{\LRt{n}(p-4)\LRf{k-4}}\LRf{{\BP}^{d-4}, 1} &p=4, \ldots, d\\
\cdots\\
\ora{CW}^{\LRt{n}(p-d+2)\LRf{k-(d-2)}}_p\LRf{{\BP}^{2}, 1} &p=d-2, d-1, d\\
0 & {\rm otherwise}\\
} 
$$
Therefore
$$
\ora{CW}^{\LRt{n}pk}\LRf{{\BP}^d, 1} =
\LBrace{{ll}
\ker\LRf{\b_k^{\LRt{n+k}}} & p=0\\
co\ker\LRf{\b_k^{\LRt{n+k}}} & p=1\\
\ker\LRf{\b_{k-2}^{\LRt{n+k-2}}} & p=2\\
co\ker\LRf{\b_{k-2}^{\LRt{n+k-2}}} & p=3\\
\cdots\\
%\ker\LRf{\b_{k-(r-2)}^{\LRt{n+k-(r+2)}}} &p=k-(r-2)\\
%co\ker\LRf{\b_{k-(r-2)}^{\LRt{n+k-(r+2)}}} &p= k-(r-1)\\
%G{\CW}_{k-r}^{\LRt{n+k-r}}\LRf{{\BF}}&p= k-r\\
\ker\LRf{\b_{k-(d-2)}^{\LRt{n+k-(d-2)}}} &p=d-2\\
co\ker\LRf{\b_{k-(d-2)}^{\LRt{n+k-(d-2)}}} &p= d-1\\
G{\CW}_{k-d}^{\LRt{n+k-d}}\LRf{{\BF}}&p= d\\
0 & {\rm otherwise}\\
} 
$$
%%%
Assume $d$ is odd. Inductively, then
$$
\ora{CW}^{\LRt{n}pk}\LRf{{\BP}^d, 1} =
\LBrace{{ll}
\ker\LRf{\b_k^{\LRt{n+k}}} & p=0\\
co\ker\LRf{\b_k^{\LRt{n+k}}} & p=1\\
\ker\LRf{\b_{k-2}^{\LRt{n+k-2}}} & p=2\\
co\ker\LRf{\b_{k-2}^{\LRt{n+k-2}}} & p=3\\
%\ora{CH}^{\LRt{n}\TCP{k-4}}_p\LRf{{\BP}^{\TCP{r-4}}, {\CO}\LRf{\TCP{1}}} &p=k-4, \ldots, k-r\\
\cdots\\
%\ora{CH}^{\LRt{n}\TCP{k-(r-3)}}_p\LRf{{\BP}^{\TCP{3}}, {\CO}\LRf{\TCP{1}}} &p=k-(r-3),\ldots, k-r\\
%\ker\LRf{\partial_{k-(r-3)}^{\LRt{n+k-(r-3)}}} &p=k-(r-3)\\
%co\ker\LRf{\partial_{k-(r-3)}^{\LRt{n+k-(r-3)}}}&p=k-(r-2)\\
\ker\LRf{\b_{k-(d-1)}^{\LRt{n+k-(d-1)}}}&p=d-1\\
co\ker\LRf{\b_{k-(d-1)}^{\LRt{n+k-(d-1)}}}&p= d\\
0 & {\rm otherwise}\\
} 
$$
\pic $\eop$ 
%%%%%%%%%%%%%%%%%%%%%%%%%%%%%%%%%%%%%%%%%%%%
%%%%%%%%%%%%%%%%%%%%%%%%%%%%%%%%%%%%%%%%%%%%
%%%%%%%%%%%%%%%%%%%%%%%%%%%%%%%%%%%%%%%%%%%%
%%%%%%%%%%%%%%%%%%%%%%%%%%%%%%%%%%%%%%%%%%%%
%%%%%%%%%%%%%%%%%%%%%%%%%%%%%%%%%%%%%%%%%%%%
%%%%%%%%%%%%%%%%%%%%%%%%%%%%%%%%%%%%%%%%%%%%
%%%%%%%%%%%%%%%%%%%%%%%%%%%%%%%%%%%%%%%%%%%%
%%%%%%%%%%%%%%%%%%%%%%%%%%%%%%%%%%%%%%%%%%%%
%%%%%%%%%%%%%%%%%%%%%%%%%%%%%%%%%%%%%%%%%%%%
 %%%%%%%%%%%%%%%%%%%%%%%%%%%%%%%%%%%%%%%%%%%%%%
%%%%%%%%%%%%%%%%%%%%%%%%%%%%%%%%%%%%%%%%%%%%%% 
%%%%%%%%%%%%%%%%%%%%%%%%%%%%%%%%%%%%%%%%%%%%%% 

\subsection{Other descriptions} \label{secOtheDes} %{Better on $\b^{\LRt{n}}_k$} 
In this section, we provide other descriptions of the Chow-Witt groups $\ora{CW}^{\LRt{n}pk}\LRf{{\BP}^d, \l}$. In fact, we only give other descriptions of the expressions:
$$
\LBrace{{l}
\ker\LRf{\b_k^{\LRt{n}}} \\
co\ker\LRf{\b_k^{\LRt{n}}}\\
}
$$
By definition of $\b^{\LRt{n}}_k$ is given by the following commutative diagram:
$$
\diagram 
 G{\CW}^{\LRt{n-1}}_{k-1}\LRf{{\BF}}\ar[d]_{\D_1}^{\wr} &G{\CW}^{\LRt{n}}_k\LRf{{\BF}} \ar@{^(->}[d]^{\pi} \ar@{-->}[l]_{\b^{\LRt{n}}_k}  \\
G{\CW}^{\LRt{n-1}}_{k-1}\LRf{C{\BM}^1{\BF}\LRt{T}_{\LRf{T}}} &G{\CW}^{\LRt{n}}_k\LRf{{\BF}\LRf{T}} \ar[l]^{\quad\qquad \partial_T} \\
\enddiagram 
$$ 
where $\D_1$ is the d\'{e}vissage isomorphism, $\pi$ is the pullback and $\partial_T$ is the residue homomorphism. 
Subsequently, we
 consider ${\BA}^1=\LRf{X_1\neq 0}\subseteq {\BP}^1$ and $T=\frac{X_1}{X_0}$. 
First, we have the following descriptions.
\bL{\rm 
Let ${\BF}$ be a field with $\frac{1}{2}\in {\BF}$. Then 
$$
\LBrace{{l}
\ker\LRf{\b_k^{\LRt{n}}} 
=\ker\LRf{\partial_{{\BP}}}\\
co\ker\LRf{\b_k^{\LRt{n}}}=co\ker\LRf{\partial_{{\BP}}}\\
}
$$
where $\partial_{{\BP}}$  is  the connecting homomorphism as in the diagram below
(\ref{actuallyPartol}), in the proof. 
}
\eL
\pf We let ${\BA}^1=\LRf{X_1\neq 0}\subseteq {\BP}^1$. Consider the map 
$$
\diagram 
G{\CW}^{\LRt{n-1}}\LRf{C{\BM}^1{\BP}^1, 1} \ar[r] \ar[d]& G{\CW}^{\LRt{n}}\LRf{C{\BM}^0{\BP}^1, 1} \ar[r] \ar[d]& 
G{\CW}^{\LRt{n}}\LRf{{\BF}\LRf{{\BP}^1}, 1} \ar[d]\\
G{\CW}^{\LRt{n-1}}\LRf{C{\BM}^1{\BA}^1} \ar[r] & G{\CW}^{\LRt{n}}\LRf{C{\BM}^0{\BA}^1} \ar[r] & 
G{\CW}^{\LRt{n}}\LRf{{\BF}\LRf{{\BA}^1}} \\
\enddiagram
$$
of the $G{\CW}$-spectra, where the vertical maps are restriction maps.
%By  \cite[Rem 9.11, pp 1812]{S17}, \cite{W3}, there is a  homotopy equivalence ${\BK}\LRf{{\BP}^1} \iso G{\CW}^{\LRt{n}}\LRf{C{\BM}^0{\BP}^1, 1}$ of spectra. 
 Combining, $\forall k\in {\BZ}$,  we obtain the commutative diagram:
\begin{equation}\label{actuallyPartol}
\diagram 
&&
G{\CW}^{\LRt{n}}_k\LRf{{\BF}}\ar[d]^{\pi}\ar@/_/@{-->}[lld]_{\b^{\LRt{n}}_k}\\
G{\CW}^{\LRt{n-1}}_{k-1}\LRf{{\BF}}\ar@{=}[d]_{\D_1}& G{\CW}^{\LRt{n}}_k\LRf{{\BF}\LRf{{\BP}^1}, 1}
\ar@/_/[dl]^{\partial_{X_1}}  \ar[d]^{\partial_{{\BP}}} \ar[r]^{\sim}_{\iota} & G{\CW}^{\LRt{n}}_k\LRf{{\BF}\LRf{{\BA}^1}} \ar[d]^{\partial_{{\BA}}}\\
  G{\CW}^{\LRt{n-1}}_{k-1}\LRf{C{\BM}^1{\BP}^1_{\LRf{X_1}}, 1}\ar[r]& G{\CW}^{\LRt{n-1}}_{k-1}\LRf{C{\BM}^1{\BP}^1, 1}\ar[r]&G{\CW}^{\LRt{n-1}}_{k-1}\LRf{C{\BM}^1{\BA}^1}\\
\enddiagram 
\end{equation}
where %vertical maps are restrictionmmaps, and 
 $\partial_{\BP}$ and $\partial_{{\BA}}$ are the connecting homomorphism. 
The $3^{rd}$-vertical line is a short exact sequence (\ref{arrowticExactAll}), and $\iota$ is the obvious isomorphism. 
By Snake lemma, the sequence
$$
\diagram 
0 \ar[r] & \ker\LRf{\partial_{{\BP}}} \ar[r] & G{\CW}^{\LRt{n}}_k\LRf{{\BF}} \ar[r]^{\b^{\LRt{n}}_k} & 
G{\CW}^{\LRt{n-1}}_{k-1}\LRf{{\BF}} \ar[r] & co\ker\LRf{\partial_{{\BP}}}\ar[r] & 0\\
\enddiagram
$$
is exact. 
\pic $\eop$ 
%%%%%%%%%%%%%%%%%%%%%%%%%%%%%%%%%%%%%%%%%%%%%% 

We provide a further description of $\ker\LRf{\b_k^{\LRt{n}}}$ and $co\ker\LRf{\b_k^{\LRt{n}}}$, in terms of 
Algebraic Bott sequence \cite[PP 1784]{S17}. Recall the following from \cite[Thm 6.1 pp 1784]{S17}.
\bT\label{MarcoPolo}{\rm 
Let ${\SA}$ be a dg category with weak equivalences and duality, with $\frac{1}{2}\in {\SA}$. Let $n\in {\BZ}$.
Then the sequence
\begin{equation}\label{Bott611784}
\diagram
G{\CW}^{\LRt{n-1}}\LRf{{\SA}} \ar[r]^{\quad F} & {\BK}\LRf{{\SA}} \ar[r]^{H\quad} & G{\CW}^{\LRt{n}}\LRf{{\SA}}\\
\enddiagram 
\end{equation}
is a homotopy fibration of $G{\CW}$-spectra, where $F$ is the forgetful and $H$ is the hyperbolic functor.  
This sequence is  referred to as the Algebraic Bott sequence.
}
\eT

%%%%%%%%%%%%%%%%%%%%%%%%%%%%%%%%%%%%%%%%%%%%%% 
\vspace{3mm}
The following particular case of the Algebraic Bott sequence (\ref{Bott611784}) is of our subsequent interest. 
\bC\label{BottBFfile}{\rm 
Let ${\BF}$ be a field, with  $\frac{1}{2}\in {\BF}$. Then $\forall n\in {\BZ}$ there is a sequence 
\begin{equation}\label{FilduVishu}
\diagram
G{\CW}^{\LRt{n-1}}\LRf{{\BF}} \ar[r]^{\quad F} & {\BK}\LRf{{\BF}} \ar[r]^{H\quad} & G{\CW}^{\LRt{n}}\LRf{{\BF}}\\
\enddiagram 
\end{equation}
homotopy fibration of $G{\CW}$-spectra. 
Consequently, $\forall k\in {\BZ}$ the is a natural connecting homomorphism 
$$
\diagram
G{\CW}_{k}^{\LRt{n}}\LRf{{\BF}} \ar[rr]^{\beta_{k}^{\LRt{n}}} && G{\CW}_{k-1}^{\LRt{n-1}}\LRf{{\BF}}\\
\enddiagram
$$
}
\eC
\pf Follows form (\ref{MarcoPolo}), with ${\SA}={\bf dg}\LRf{{\BF}}$. $\eop$
%%%%%%%%%%%%%%%%%%%%%%%%%%%%%%%%%%%%%%%%%%%%%% 

\vspace{3mm}
Now we have the other description.
\bP\label{betaBohu}{\rm 
Let ${\BF}$ be a field, with  $\frac{1}{2}\in {\BF}$. Then $\forall n, k\in {\BZ}$, we have
$$
\LBrace{{l}
\ker\LRf{\b^{\LRt{n}}_k} = \ker\LRf{\beta^{\LRt{n}}_k}\\
co\ker\LRf{\b^{\LRt{n}}_k} = co\ker\LRf{\beta^{\LRt{n}}_k}\\
}
$$
}
\eP
\pf With notations as in (\ref{actuallyPartol}), we prove
$$
\LBrace{{l}
\ker\LRf{\partial_{{\BP}}} = \ker\LRf{\beta^{\LRt{n}}_k}\\
co\ker\LRf{\partial_{{\BP}}} = co\ker\LRf{\beta^{\LRt{n}}_k}\\
}
$$
We consider ${\BP}^1=\proj{{\BF}\LRt{X_0, X_1}}$.
Dualities considered are induced by the duality $P\mapsto {\CH}om\LRf{P, {\CO}\LRf{1}}$
 on $C{\BM}^0{\BP}^1$.  Consider the following diagram of maps of dg categories
 \begin{equation}\label{dgGapMap}
 \diagram
 {\bf dg}^{\LRt{n-1}}\LRf{{\BF}} \ar[r]^{\f} &  {\CH}\LRf{{\bf dg}\LRf{{\BF}}} 
 \ar[r]^{\h} \ar[d]^{\F} &  {\bf dg}^{\LRt{n}}\LRf{{\BF}} \ar[d]^{\tau}\\
  {\bf dg}^{\LRt{n-1}}C{\BM}^1\LRf{{\BP}^1} \ar[r]_{\c} &  {\bf dg}^{\LRt{n}}C{\BM}^0\LRf{{\BP}^1}
 \ar[r]_{\t} &  {\bf dg}^{\LRt{n}}{\BF}\LRf{{\BP}^1} \\
 \enddiagram 
 \end{equation}
 where ${\CH}\LRf{-}$ denotes the hyperbolic category,  $\f$ is the forgetful and $\h$ is the hyperbolic functors. 
 Here 
 $$
 \LBrace{{l}
 \f\LRf{M_{\bul}}=\LRf{M_{\bul}, M_{\bul}^{\#_{n-1}}}\\
 \h\LRf{M_{\bul}, N_{\bul}}= M_{\bul}\oplus N_{\bul}^{\#_{n}}\\
 }
 $$
 where $\LRf{-}^{\#_n}$ denotes the $n$-shifted duality. 
 The map $\F$ is defined by the composition:
 $$
 \diagram 
 %{\bf dg}\LRf{{\BF}} \ar[rr]^{\pi^*} && {\bf dg}C{\BM}^0\LRf{{\BP}^1}\\
 {\CH}\LRf{{\bf dg}\LRf{{\BF}}} \ar[rr]^{\F_0}\ar@/_/[drr]_{\F} && {\CH}\LRf{ {\bf dg}C{\BM}^0\LRf{{\BP}^1}}\ar[d]^{\H}\\
 && {\bf dg}^{\LRt{n}}C{\BM}^0\LRf{{\BP}^1}\\
 \enddiagram 
 \qquad 
 \LBrace{{l}
 {\F}_0\LRf{M_{\bul}, N_{\bul}} =\LRf{M_{\bul}{\CO}, N_{\bul}{\CO}} \\
 {\H}\LRf{P_{\bul}, Q_{\bul}}= P_{\bul}\oplus Q_{\bul}^{\#_{n, 1}}\\
  {\F}\LRf{M_{\bul}, N_{\bul}}=M_{\bul}{\CO}\oplus N_{\bul}{\CO}^{\#_{n, 1}}\\
 }
 $$ 
 where $M_{\bul}{\CO}:= M_{\bul}\otimes {\CO}$, and likewise, and  $\LRf{-}^{\#_{n,1}}$ denotes the $n$-shifted duality, induced by $P\mapsto {\CH}om\LRf{P, {\CO}\LRf{1}}$,
 .

  It follows that
  the righthand square of diagram (\ref{dgGapMap}) commutes. So, (\ref{dgGapMap})
   induces a map 
  \begin{equation}\label{dgGabrSlskMessi}
 \diagram
 G{\CW}^{\LRt{n-1}}\LRf{{\BF}} \ar[r]^{\LR|{\f}} \ar@{-->}[d]_{\theta}&  {\BK}\LRf{{\BF}}
 \ar[r]^{\LR|{\h}} \ar[d]^{\LR|{\F}} &   G{\CW}^{\LRt{n}}\LRf{{\BF}} \ar[d]^{\LR|{\tau}}\\
  G{\CW}^{\LRt{n-1}}C{\BM}^1\LRf{{\BP}^1} \ar[r]_{\LR|{\c}} &  G{\CW}^{\LRt{n}}\LRf{{\BP}^1, 1}
 \ar[r]_{\LR|{\t}} &  G{\CW}^{\LRt{n}}{\BF}\LRf{{\BP}^1} \\
 \enddiagram 
 \end{equation}
 of homotopy $G{\CW}$-spectra.
 The upper sequence in the Algebraic Bott fibration (\ref{FilduVishu}). The lower line is also a homotopy fibration of 
 $G{\CW}$-spectra \cite[pp 473]{M23}. The map $\theta$ ({\it which is not relevant for out subsequent discussions}) is induced by the properties of homotopy fibrations. 
 However, $\LR|{\F}$ is a homotopy equivalence \cite[Rem 9.11 pp 1812]{S17}, \cite{W3}. This induces the following diagram of exact sequences of the homotopy groups and connecting homomorphisms 
 \begin{equation}\label{dataGanesh} 
  \diagram
 {\BK}_{k}\LRf{{\BF}}
 \ar[r]^{\LR|{\h}} \ar[d]^{\LR|{\F}}_{\wr} &   G{\CW}_{k}^{\LRt{n}}\LRf{{\BF}} \ar@{^(->}[d]^{\LR|{\tau}}\ar[r]^{\beta^{\LRt{n}}_k}
 &  G{\CW}^{\LRt{n-1}}_{k-1}\LRf{{\BF}}  %\ar@{-->}[d]_{\LR|{\theta}}
 \ar[r]^{\LR|{\f}}
 & {\BK}_{k-1}\LRf{{\BF}} \ar[d]^{\LR|{\F}}_{\wr} \\
  G{\CW}_{k}^{\LRt{n}}\LRf{{\BP}^1, 1}
 \ar[r]_{\LR|{\t}} &  G{\CW}_{k}^{\LRt{n}}{\BF}\LRf{{\BP}^1} \ar[r]_{\partial_{{\BP}}\quad}&   G{\CW}_{k-1}^{\LRt{n-1}}C{\BM}^1\LRf{{\BP}^1}\ar[r]_{\LR|{\c}}& G{\CW}_{k-1}^{\LRt{n}}\LRf{{\BP}^1, 1}\\
 \enddiagram 
 \end{equation}
 The map $\LR|{\tau}$ is also injective (\ref{arrowticExactAll}). Therefore, $\LR|{\tau}$ induces an injective map
 $$
 \diagram
\ker\LRf{\beta^{{\LRt{n}}}_k} \ar@{^(->}[r]^{\LR|{\tau}} &  \ker\LRf{\partial_{{\BP}}} \\
 \enddiagram
 $$
 We need a proof that $\LR|{\tau}$ maps $\ker\LRf{\beta^{{\LRt{n}}}_k}$ into $\ker\LRf{\partial_{{\BP}}}$. 
 Suppose $x\in \ker\LRf{\beta^{\LR|{n}}_k}= Image\LRf{\LR|{\h}}$. So, $x=\LR|{\h}\LRf{x_0}$
 and $\LR|{\tau}\LRf{x}=\LR|{\tau}\LR|{\h}\LRf{x_0}=\LR|{\t}\LR|{\F}\LRf{x_0}\in \ker\LRf{\partial_{{\BP}}}$.
 
 To prove subjectivity of the map, let  $y\in  \ker\LRf{\partial_{{\BP}}}$. So, 
 $$
 y=\LR|{{\t}}\LR|{{\F}}\LRf{z}, ~~{\rm for~some}~z\in {\BK}_k\LRf{{\BF}}. \quad {\rm So.}\quad  
% $$
 \LR|{\tau}\LR|{\h}\LRf{z}= \LR|{{\t}}\LR|{{\F}}\LRf{z}=y
 $$
 So, the map $ \diagram
\ker\LRf{\beta^{{\LRt{n}}}_k} \ar[r]^{\LR|{\tau}}_{\sim} &  \ker\LRf{\partial_{{\BP}}} \\
 \enddiagram$ is an isomorphism. 
 Similarly, $\LR|{\F}$ induces an injective map
 \begin{equation}\label{July272nd}
 \diagram
 co\ker\LRf{\beta^{{\LRt{n}}}_k} \ar@{^(->}[r]^{\LR|{\F}} &  co\ker\LRf{\partial_{{\BP}}}\\
 \enddiagram 
 \end{equation}
 We prove that $\LR|{\F}$ maps $co\ker\LRf{\beta^{{\LRt{n}}}_k}$
 into  $co\ker\LRf{\partial_{{\BP}}}$. Let $x\in co\ker\LRf{\beta^{{\LRt{n}}}_k}= Image\LRf{\LR|{\f}}$. So, 
 $x=\LR|{\f}\LRf{x_0}$ for some $x_0$. So, $\LR|{\t}\LR|{\F}\LRf{x}= \LR|{\t}\LR|{\F}\LR|{\f}\LRf{x_0}=
 \LR|{\tau}\LR|{\h}\LR|{\f}\LRf{x_0}=0$. So, $\LR|{\F}\LRf{x} \in \ker\LRf{\LR|{\t}}=Image\LRf{\LR|{\c}}=co\ker\LRf{\partial_{{\BP}}}$. So, (\ref{July272nd}) is an injective map, as shown. 
 
 Now let $z\in  co\ker\LRf{\partial_{{\BP}}}=Image\LRf{\LRt{\c}}$. Then $z=\LR|{\c}\LRf{y}$ for some $y$. So, 
 $z=\LR|{\F}\LRf{x}$ for some $x$. Since $\LR|{\tau}$ is injective, it follows 
 $x\in \ker\LRf{\LR|{\h}}=Image\LRf{\LR|{\f}}=co\ker\LRf{\beta_k^{\LRt{n}}}$. Therefore, 
 $
 \diagram
 co\ker\LRf{\beta^{{\LRt{n}}}_k} \ar[r]^{\LR|{\F}}_{\sim} &  co\ker\LRf{\partial_{{\BP}}}\\
 \enddiagram 
 $
 is an isomorphism. $\eop$
%%%%%%%%%%%%%%%%%%%%%%%%%%%%%%%%%%%%%%%%%%%%%%
%%%%%%%%%%%%%%%%%%%%%%%%%%%%%%%%%%%%%%%%%%%%%% 
%%%%%%%%%%%%%%%%%%%%%%%%%%%%%%%%%%%%%%%%%%%%%% 
%%%%%%%%%%%%%%%%%%%%%%%%%%%%%%%%%%%%%%%%%%%%%% 
%%%%%%%%%%%%%%%%%%%%%%%%%%%%%%%%%%%%%%%%%%%%%%%%%%%%%%%%%
%%%%%%%%%%%%%%%%%%%%%%%%%%%%%%%%%%%%%%%%%%%%%%%%%%%%%%%%%

%\noindent\TCP{\bf Vanishing 10 August 2026}

\subsection{The Algebraic Bott Fibration} \label{Secmoidyam} %Refreshing: 10 August 2026} 

We close this section by drawing attention  to some of the important cases and observations. First, we wish the put our results on projective spaces, 
in the mold of Algebraic Bott fibration \cite[pp 1782]{S17}. This aspect establishes the fundamental nature of the same. 
Fix $k, n\in {\BZ}$. Consider the Algebraic Bott fibration
$$
\diagram 
G{\CW}^{\LRt{n+k-1}}\LRf{{\BF}} \ar[r]^{\qquad F} & {\BK}\LRf{{\BF}} \ar[r]^{H \quad} & G{\CW}^{\LRt{n+k}}\LRf{{\BF}} \\
\enddiagram 
$$
We can write down the corresponding long exact sequence
$$
\diagram 
\cdots \ar[r] &G{\CW}_k^{\LRt{n+k-1}}\LRf{{\BF}} \ar[r]^{\qquad\f_k^{n+k-1}} & {\BK}_k\LRf{{\BF}} \ar[r]^{\h_k^{n+k}\quad} & G{\CW}_k^{\LRt{n+k}}\LRf{{\BF}} \ar@/_/[lld]^{\beta^{\LRt{n+k}}_k} &\\
&G{\CW}_{k-1}^{\LRt{n+k-1}}\LRf{{\BF}} \ar[r]^{\qquad\f_{k-1}^{n+k-1}} & {\BK}_{k-1}\LRf{{\BF}} \ar[r]^{\h_{k-1}^{n+k}\quad} & G{\CW}_{k-1}^{\LRt{n+k}}\LRf{{\BF}}  \ar@/_/[lld]^{\beta^{\LRt{n+k}}_{k-1}} &\\
&G{\CW}_{k-2}^{\LRt{n+k-1}}\LRf{{\BF}} \ar[r]^{\qquad\f_{k-2}^{n+k-1}} & {\BK}_{k-2}\LRf{{\BF}} \ar[r]^{\h_{k-2}^{n+k}\quad} & G{\CW}_{k-2}^{\LRt{n+k}}\LRf{{\BF}} \ar@/_/[lld]^{\beta^{\LRt{n+k}}_{k-2}} &\\
&G{\CW}_{k-3}^{\LRt{n+k-1}}\LRf{{\BF}} \ar[r]^{\qquad\f_{k-3}^{n+k-1}} & {\BK}_{k-3}\LRf{{\BF}} \ar[r]^{\h_{k-3}^{n+k}\quad} & G{\CW}_{k-3}^{\LRt{n+k}}\LRf{{\BF}}  \ar@/_/[lld]^{\beta^{\LRt{n+k}}_{k-3}}&\\
&G{\CW}_{k-4}^{\LRt{n+k-1}}\LRf{{\BF}} \ar[r]^{\qquad\f_{k-4}^{n+k-1}} & {\BK}_{k-4}\LRf{{\BF}} \ar[r]^{\h_{k-4}^{n+k}\quad} & G{\CW}_{k-4}^{\LRt{n+k}}\LRf{{\BF}}  \ar[r]^{\qquad \beta^{\LRt{n+k}}_{k-4}}& \cdots\\
\enddiagram 
$$
of $G{\CW}$-groups. 
We have 
$$
\begin{array}{ll}
\ora{CH}^{\LRt{n}pk}\LRf{{\BP}^d, 1}=&
\LBrace{{lr}
\ker\LRf{\beta_k^{{\LRt{n+k}}}} =Im\LRf{\h_k^{n+k}}& p=0\\
co\ker\LRf{\beta_k^{{\LRt{n+k}}}}=\ker\LRf{\f_{k-1}^{n+k}} & p=1\\
}\\
\ora{CH}^{\LRt{n+1}pk}\LRf{{\BP}^d, 0}=&
\LBrace{{lr}
\ker\LRf{\beta_{k-1}^{{\LRt{n+1+k-1}}}} =Im\LRf{\h_{k-1}^{n+k}}& p=1\\
co\ker\LRf{\beta_{k-1}^{{\LRt{n+k}}}}=\ker\LRf{\f_{k-2}^{n+k-1}}& p=2\\
}\\
\ora{CH}^{\LRt{n+2}pk}\LRf{{\BP}^d, 1}=&
\LBrace{{lr}
\ker\LRf{\beta_{k-2}^{{\LRt{n+k}}}} =Im\LRf{\h_{k-2}^{n+k}}& p=2\\
co\ker\LRf{\beta_{k-2}^{{\LRt{n+k}}}}=\ker\LRf{\f_{k-3}^{n+k-1}} & p=3\\
}\\
\ora{CH}^{\LRt{n+3}pk}\LRf{{\BP}^d, 0}=&
\LBrace{{lr}
\ker\LRf{\beta_{k-3}^{{\LRt{n+k}}}} =Im\LRf{\h_{k-3}^{n+k}}& p=3\\
co\ker\LRf{\beta_{k-3}^{{\LRt{n+k}}}}=\ker\LRf{\f_{k-4}^{n+k-1}}& p=4\\
}\\
\end{array}
$$
The shift being four periodic, varying $n=0,1, 2, 3~mod~4$, we obtain all the Chow-Witt groups 
$\ora{CH}^{\LRt{n}pk}\LRf{{\BP}^d, \l}$. We summarize all these, as follows.
%%%
\bC\label{alsBor10Aug}{\rm 
Let ${\BF}$ be a field, with $\frac{1}{2}\in {\BF}$. Fix, $n, k\in {\BZ}$. Then 
\begin{equation}\label{GanesherAug10Ashirvad}
\begin{array}{l}
\ora{CW}^{\LRt{n}pk}\LRf{{\BP}^d,  {\CO}\LRf{1}} =
\LBrace{{ll}
Im\LRf{\h_{k-p}^{\LRt{n+k-p}}} & 0\leq p=2r\leq d-1\\
\ker\LRf{\f_{k-p}^{\LRt{n+k-p}}} & 1\leq p=2r+1\leq d\\
G{\CW}^{\LRt{n+k-d}}_{k-d}\LRf{{\BF}} & p=d~even\\
%\ora{CW}^{\LRt{n}(p-1)(k-1)}\LRf{{\BP}^{d-1}, 0} &p=2, \ldots, d\\
0 & Otherwise\\
} \\
%$$
%and 
%$$
\ora{CW}^{\LRt{n}pk}\LRf{{\BP}^d,  {\CO}\LRf{0}} =
\LBrace{{ll}
G{\CW}^{\LRt{n+k}}_k\LRf{{\BF}} & p=0\\
%\ora{CW}^{\LRt{n}(p-1)(k-1)}\LRf{{\BP}^{d-1}, {\CO}\LRf{1}} & p=1, \ldots, d\\
Im\LRf{\h_{k-p}^{\LRt{n+k-p}}} & 1\leq p=2r+1\leq d-1\\ %&0\leq p-1=2r\leq d-2\\
\ker\LRf{\f_{k-p}^{\LRt{n+k-p}}} & 2\leq p=2r+2\leq d\\ %& 1\leq p-1=2r+1\leq d-1\\
G{\CW}^{\LRt{n+k-d}}_{k-d}\LRf{{\BF}} & p=d~odd\\
0 & Otherwise\\
} \\
\end{array}
\end{equation} 
}
\eC
\pf Follows from the above and Theorem \ref{mainPSpace}. $\eop$

\vspace{3mm}
 The shift 
$n=0$ case is also of special interest, because it corresponds to the motivic Chow-Witt groups \cite{BL8, M12, F8, F13}.
We write it down explicitly, as follows. 

\bC\label{niszeroface}{\rm 
Let ${\BF}$ be a field, with $\frac{1}{2}\in {\BF}$, and $n, k\in {\BZ}$. Then 
$$
\ora{CW}^{\LRt{0}pk}\LRf{{\BP}^d, {\CO}\LRf{1}}=
\LBrace{{ll}
\ker\LRf{\beta_{k-2r}^{\LRt{k-2r}}} & 0\leq p =2r\leq d-1\\
co\ker\LRf{\beta_{k-2r}^{\LRt{k-2r}}} & 1\leq p =2r+1\leq d\\
G{\CW}_{k-d}^{\LRt{k-d}}\LRf{{\BF}} & p=d=2r~even\\
0 & Otherwise\\
}
$$
And,
$$
\ora{CW}^{\LRt{0}pk}\LRf{{\BP}^d, {\CO}\LRf{0}}=
\LBrace{{ll}
G{\CW}^{\LRt{k}}_k\LRf{{\BF}} & p=0\\
%\ora{CW}^{\LRt{0}(p-1)(k-1)}\LRf{{\BP}^d, {\CO}\LRf{1}} & p=1, \ldots, d\\
\ker\LRf{\beta_{k-1-2r}^{\LRt{k-1-2r}}} & 1\leq p =2r+1\leq d-1\\
%&0\leq p-1 =2r\leq d-2\\
co\ker\LRf{\beta_{k-1-2r}^{\LRt{k-1-2r}}} & 2\leq p =2r+2\leq d\\
%& 1\leq p -1=2r+1\leq d-1\\
G{\CW}_{k-d}^{\LRt{k-d}}\LRf{{\BF}} & p=d=2r+1~odd\\
0 & Otherwise\\
}
$$

}
\eC
\pf This is an immediate consequence of (\ref{betaBohu}) and Theorem \ref{mainPSpace}. $\eop$

%\subsection{Vanishing groups}
When $k-p \leq -1$ is negative then most of the Chow-Witt groups $\ora{CW}^{\LRt{n}pk}\LRf{{\BP}^d, \l}=0$ vanish and others reduce to Witt groups. 
 Recall, $\p:=k-p$ is used, as the alternate homological indexing. We give a list of groups that vanish. 
 
 % So, the following is a statement on Chow-Witt groups in the negative homological degrees. 
\bC\label{udhao}{\rm 
Suppose ${\BF}$ is a field and $\frac{1}{2}\in {\BF}$. Then 
\begin{equation}\label{GanesherAshirvadInsert}
\begin{array}{l}
\ora{CW}^{\LRt{n}pk}\LRf{{\BP}^d,  {\CO}\LRf{1}} =
\LBrace{{ll}
{\bf 0} & 0\leq p=2r\leq d-1,~k-2r<0\\
{\bf 0}  & 0\leq p=2r+1\leq d, ~k-2r<0\\
W^{\LRt{n}}\LRf{{\BF}}& p=d~is~even~if~k<d\\
{\bf 0} & unless~p=0, 1, \ldots, d\\
} \\
\ora{CW}^{\LRt{n}pk}\LRf{{\BP}^d,  {\CO}\LRf{0}} =
\LBrace{{ll}
W^{\LRt{n}}\LRf{{\BF}}
 & p=0, k<0\\
{\bf 0} & 0\leq p-1=2r\leq d-2,~k-2r<1\\
{\bf 0} & 0\leq p-1=2r+1\leq d-1,~k-2r<1\\
{W}^{\LRt{n}}\LRf{{\BF}}& p=d~odd, ~if~k<d\\
{\bf 0} & unless~p=0, 1, \ldots, d\\
} \\
\end{array}
\end{equation} 
Further, note that $W^{\LRt{n}}\LRf{{\BF}}= \LBrace{{ll}
{\bf 0} & if~n\neq 0~mod~4\\
W\LRf{{\BF}} & if~n= 0~mod~4}$.
}
\eC 
\pf Follow from the fact that $\beta^{\LRt{n}}_r$ is an isomorphism, when $r\leq -1$. $\eop$

%%%%%%%%%%%%%%%%%%%%%%%%%%%%%%%%%%%%%%%%%%%%%%%%%%%%%%%%%
%%%%%%%%%%%%%%%%%%%%%%%%%%%%%%%%%%%%%%%%%%%%%%%%%%%%%%%%%
%%%%%%%%%%%%%%%%%%%%%%%%%%%%%%%%%%%%%%%%%%%%%%%%%%%%%%%%%
%\noindent\TCP{\bf endVanishing}
%%%%%%%%%%%%%%%%%%%%%%%%%%%%%%%%%%%%%%%%%%%%%%%%%%%%%%%%%
%%%%%%%%%%%%%%%%%%%%%%%%%%%%%%%%%%%%%%%%%%%%%%%%%%%%%%%%%
%%%%%%%%%%%%%%%%%%%%%%%%%%%%%%%%%%%%%%%%%%%%%%%%%%%%%%%%%
%%%%%%%%%%%%%%%%%%%%%%%%%%%%%%%%%%%%%%%%%%%%%% 
%%%%%%%%%%%%%%%%%%%%%%%%%%%%%%%%%%%%%%%%%%%%%% 
%%%%%%%%%%%%%%%%%%%%%%%%%%%%%%%%%%%%%%%%%%%%%% 
%%%%%%%%%%%%%%%%%%%%%%%%%%%%%%%%%%%%%%%%%%%%%% 
%%%%%%%%%%%%%%%%%%%%%%%%%%%%%%%%%%%%%%%%%%%%%%%%%%%%%%%%%
%%%%%%%%%%%%%%%%%%%%%%%%%%%%%%%%%%%%%%%%%%%%%%%%%%%%%%%%%
%%%%%%%%%%%%%%%%%%%%%%%%%%%%%%%%%%%%%%%%%%%%%%%%%%%%%%%%%
%%%%%%%%%%%%%%%%%%%%%%%%%%%%%%%%%%%%%%%%%%%%%%%%%%%%%%%%%
%%%%%%%%%%%%%%%%%%%%%%%%%%%%%%%%%%%%%%%%%%%%%%%%%%%%%%%%%
%%%%%%%%%%%%%%%%%%%%%%%%%%%%%%%%%%%%%%%%%%%%
%%%%%%%%%%%%%%%%%%%%%%%%%%%%%%%%%%%%%%%%%%%%
%%%%%%%%%%%%%%%%%%%%%%%%%%%%%%%%%%%%%%%%%%%%
%%%%%%%%%%%%%%%%%%%%%%%%%%%%%%%%%%%%%%%%%%%%
%%%%%%%%%%%%%%%%%%%%%%%%%%%%%%%%%%%%%%%%%%%%
 
%%%%%%%%%%%%%%%%%%%%%%%%%%%%%%%%%%%%%%%%%%%%%%%%%%%%%%%%%
%%%%%%%%%%%%%%%%%%%%%%%%%%%%%%%%%%%%%%%%%%%%%%%%%%%%%%%%%
%%%%%%%%%%%%%%%%%%%%%%%%%%%%%%%%%%%%%%%%%%%%%%%%%%%%%%%%%
%%%%%%%%%%%%%%%%%%%%%%%%%%%%%%%%%%%%%%%%%%%%%%%%%%%%%%%%%
%%%%%%%%%%%%%%%%%%%%%%%%%%%%%%%%%%%%%%%%%%%%%%%%%%%%%%%%%
\section{Arrowtic Chow groups} \label{secChow}
The $K$-theory aspects of the theory above do not follow directly from that on $GW$-theory, while these need to be worked out by repeating the same line of arguments. We briefly outline the notations, definitions and results without proofs, in this section. 

\bD\label{eaKnotaSptrCx}{\rm 
As above, let  $X$ be a quasi projective scheme over an affine scheme $\spec{A_0}$ with $d=\dim X$. 
%\TC
Then $\forall k \in {\BZ}$, with cohomology $\deg=p=0, 1, 2, \ldots, d$, % $\deg=p=p=k, k-1, \ldots, k-d$, 
define
%$$
%$
%
\begin{equation}\label{eaKlabelNOta}
{\CC}^{pk}\LRf{X}= \bigoplus_{x\in X^{\LRf{p}}}{\BK}_{k-p}\LRf{C{\BM}^{p}\LRf{X_x}} \\
\end{equation}
Write ${\CC}^{pk}\LRf{X}=0$ unless $p=0, 1, \ldots, d$. 
For $k\in {\BZ}$,   ${\CC}^{{\bul}k}\LRf{X}$ is a complex. In this section, we work with the non-connective ${\BK}$-theory only, which makes no difference in the non-singular case. For $p, k\in {\BZ}$, define the arrowtic Chow groups, as
follows
$$
\ora{CW}^{pk}\LRf{X, {\CL}}= H^p\LRf{{\CC}^{{\bul}k}\LRf{X}} 
$$
}
\eD
The following is immediate.
\bExr\label{eaKmilJinh}{\rm 
Let ${\BF}$ be a field, and $X=\spec{{\BF}}$. Then for 
$k\in {\BZ}$ we have 
$$
{\CC}^{pk}\LRf{X}= \LBrace{{ll}
{\BK}_k\LRf{{\BF}} & if~p=0\\
0 & if ~p\neq 0\\
}
%$$
\quad {\rm So, }\quad
%$$
\ora{CH}^{pk}\LRf{{\BF}}= \LBrace{{ll}
{\BK}_k\LRf{{\BF}} & if~p=0\\
0 & if ~p\neq 0\\
}
$$
}
\eExr
%%%%%%%%%%%%%%%%%%%%%%%%%%%%%%%%%%%%%%%%%%%%
Following (\ref{defResidue}) we define the residue homomorphism. 

\bD\label{eaKdefResidue}{\rm 
Let $\LRf{R, \m, \kappa}$ be a local integral domain, with $\dim R=1$, and other notations be as in 
(\ref{defResidue}).
 There is a homotopy fibration 
$$
\diagram 
{\BK}\LRf{C{\BM}^Z\LRf{Y}} \ar[r] & {\BK}\LRf{{\SV}(Y)} \ar[r] & {\BK}\LRf{{\SV}(U)} \\
\enddiagram 
$$
of the ${\BK}$-theory spectra, and ${\BK}\LRf{{\SV}(U)}\cong {\BK}\LRf{K(R)}$.
Consequently, we obtain a connecting homomorphism, to be called the {\bf Residue  map} 
\begin{equation}\label{eaKtransfu}
\diagram 
{\BK}_{p+1}\LRf{K(R)} \ar[rr]^{\partial_{p+1}} && {\BK}_p\LRf{C{\BM}^Z\LRf{Y}} 
\enddiagram \qquad \forall p\in {\BZ}.
\end{equation}
%Sometimes, we write $\partial:=\partial^{\LRt{n+1}}_{p+1}$.
}
\eD
Analogous to (\ref{arrowticExactAll}) we have the following theorem.
\bT\label{eaKarrowticExactAll} {\rm 
Let ${\BF}$ be a field. Write $A={\BF}\LRt{T}$, and $Y:={\BA}^1=\spec{A}$.  For each $y\in Y^{\LRf{1}}$, there is a residue map $\partial_y$  (\ref{transfu}), as follows:
$$
\diagram
{\BK}_{p+1}\LRf{{\BF}(T)} \ar[rr]^{\partial_y}  &&  
{\BK}_{p}\LRf{C{\BM}^y\LRf{Y_y}} \\
\enddiagram
$$ 
 Taking direct sum, with $Y_y=\spec{A_y}$, we obtain the following sequence
\begin{equation}\label{eaKprimeSeq}
\diagram
0\ar[r] &
{\BK}_{p+1}\LRf{{\BF}} \ar[r]^{\p\quad} &
{\BK}_{p+1}\LRf{{\BF}(T)} \ar[r]^{\bigoplus \partial_y\qquad }  & \bigoplus_{y\in Y^{\LRf{1}}}
{\BK}_{p}\LRf{C{\BM}^y\LRf{Y_y}}\ar[r] &0 \\
\enddiagram 
\end{equation}  
of maps. 
The map $\p$ is the usual pullback map. 
Then the sequence is a split exact sequence.
}
\eT
\pf Same as that of (\ref{arrowticExactAll}). $\eop$ 
%%%%%%%%%%%%%%%%%%%%%%%%%%%%%%%%%%%%%%%%%%%%%%%%%%%%%%%%%
%%%%%%%%%%%%%%%%%%%%%%%%%%%%%%%%%%%%%%%%%%%%%%%%%%

Analogous to (\ref{VectorBunCW}), we obtain the homotopy invariance of ${\BK}$-theory, as follows.
\bT\label{eaKVectorBunCW}{\rm 
Let $X$ be quasi projective scheme over a noetherian affine scheme $\spec{A_0}$. Assume   $X$ is regular.  Let 
$\pi: Y \lra X$ be a vector bundle. %, meaning, $Y=\spec{\Sym{{\SE}}}$, where ${\SE}\in {\SV}\LRf{X}$. 
 Then $\forall p, k$, the pullback map
$$
\diagram 
\ora{CH}^{pk} \LRf{X} \ar[rr]^{\pi^*\quad}_{\sim\quad} && \ora{CH}^{pk}\LRf{Y} \\
\enddiagram 
$$
is  an isomorphism.  
}
\eT
\pf Same the proof of (\ref{VectorBunCW}). $\eop$

\vspace{3mm}
%%%%%%%%%%%%%%%%%%%%%%%%%%%%%%%%%%%%%%%%%%%%%%%%%%%%%%%%%
As in the $GW$-theory (\ref{ArKoChgp}), we obtain the following basic result.
\bC\label{eaKArKoChgp}{\rm 
Let ${\BF}$ be a field and ${\BA}^d=\spec{{\BF}\LRt{T_1, T_2, \ldots, T_d}}$. Then 
$$
\ora{CH}^{pk}\LRf{{\BA}^d} = \ora{CH}^{pk}\LRf{{\BF}}= \LBrace{{ll}
{\BK}_k\LRf{{\BF}} & if~p=0\\
0 & if ~p\neq 0\\
}%\qquad \forall n, k, p
$$
}
\eC
Now we describe the same for the projective spaces.
%%%%%%%%%%%%%%%%%%%%%%%%%%%%%%%%%%%%%%%%%%%%
\bT\label{eaKmainPSpace}{\rm 
Let ${\BF}$ be a field, and $k\in {\BZ}$, and ${\BP}^d=\proj{{\BF}\LRt{T_1, T_2, \ldots, T_d}}$
Then 
$$
\ora{CH}^{pk}\LRf{{\BP}^d} =
\LBrace{{ll}
{\BK}_{k-p}\LRf{{\BF}} & if ~p=0, 1, \ldots, d\\
0 & Otherwise\\
} 
$$
}
\eT
\pf Same as (\ref{mainPSpace}). $\eop$ 
%%%%%%%%%%%%%%%%%%%%%%%%%%%%%%%%%%%%%%%%%%%%%%%%%%%%%%%%%

%%%%%%%%%%%%%%%%%%%%%%%%%%%%%%%%%%%%%%%%%%%%
%%%%%%%%%%%%%%%%%%%%%%%%%%%%%%%%%%%%%%%%%%%%
%%%%%%%%%%%%%%%%%%%%%%%%%%%%%%%%%%%%%%%%%%%%
%%%%%%%%%%%%%%%%%%%%%%%%%%%%%%%%%%%%%%%%%%%%
%%%%%%%%%%%%%%%%%%%%%%%%%%%%%%%%%%%%%%%%%%%%
%%%%%%%%%%%%%%%%%%%%%%%%%%%%%%%%%%%%%%%%%%%%
%%%%%%%%%%%%%%%%%%%%%%%%%%%%%%%%%%%%%%%%%%%%
%%%%%%%%%%%%%%%%%%%%%%%%%%%%%%%%%%%%%%%%%%%%

%\TCP{PPP}
%%%%%%%%%%%%%%%%%%%%%%%%%%%%%%%%%%%%%%%%%%%%

%%%%%%%%%%%%%%%%%%%%%%%%%%%%%%%%%%%%%%%%%%%%%%%%%
%%%%%%%%%%%%%%%%%%%%%%%%%%%%%%%%%%%%%%%%%%%%%%%%%
%%%%%%%%%%%%%%%%%%%%%%%%%%%%%%%%%%%%%%%%%%%%%%%%%
%%%%%%%%%%%%%%%%%%%%%%%%%%%%%%%%%%%%%%%%%%%%%%%%%
%%%%%%%%%%%%%%%%%%%%%%%%%%%%%%%%%%%%%%%%%%%%%%%%%
%%%%%%%%%%%%%%%%%%%%%%%%%%%%%%%%%%%%%%%%%%%%%%%%%
%%%%%%%%%%%%%%%%%%%%%%%%%%%%%%%%%%%%%%%%%%%%%%%%%
%%%%%%%%%%%%%%%%%%%%%%%%%%%%%%%%%%%%%%%%%%%%%%%%%
%\newpage
 \appendix
%%%%%%%%%%%%%%%%%%%%%%%%%%%%%%%%%%%%%%%%%%%%%%%%%
\section{Machinery}
%%%%%%%%%%%%%%%%%%%%%%%%%%%%%%%%%%%%%%%%%%%%

%%%%%%%%%%%%%%%%%%%%%%%%%%%%%%%%%%%%%%%%%%%%%%%%%%%%%%
%%%%%%%%%%%%%%%%%%%%%%%%%%%%%%%%%%%%%%%%%%%%%%%%%%%%%%
%%%%%%%%%%%%%%%%%%%%%%%%%%%%%%%%%%%%%%%%%%%%%%%%%%%%%%
%%%%%%%%%%%%%%%%%%%%%%%%%%%%%%%%%%%%%%%%%%%%%%%%%%%%%%
%%%%%%%%%%%%%%%%%%%%%%%%%%%%%%%%%%%%%%%%%%%%%%%%%%%%%%
%%%%%%%%%%%%%%%%%%%%%%%%%%%%%%%%%%%%%%%%%%%%%%%%%%%%%%
%%%%%%%%%%%%%%%%%%%%%%%%%%%%%%%%%%%%%%%%%%%%%%%%%%%%%%
%%%%%%%%%%%%%%%%%%%%%%%%%%%%%%%%%%%%%%%%%%%%%%%%%%%%%%
%%%%%%%%%%%%%%%%%%%%%%%%%%%%%%%%%%%%%%%%%%%%%%%%%%%%%%
%%%%%%%%%%%%%%%%%%%%%%%%%%%%%%%%%%%%%%%%%%%%%%%%%%%%%%
%%%%%%%%%%%%%%%%%%%%%%%%%%%%%%%%%%%%%%%%%%%%%%%%%%%%%%
%%%%%%%%%%%%%%%%%%%%%%%%%%%%%%%%%%%%%%%%%%%%%%%%%%%%%%
%%%%%%%%%%%%%%%%%%%%%%%%%%%%%%%%%%%%%%%%%%%%%%%%%%%%%%

\subsection{The ${\BF}^{\star}$ action} 
We propose the following multiplication map. 
\bP\label{unitHike}{\rm 
Let ${\BF}$ be field, with $\frac{1}{2}\in {\BF}$. Let  $u\in {\BF}^{\star}$ be a unit and $n, p \in {\BZ}$. Let ${\SB}$ be a dg ${\BF}$-category, with weak equivalences and duality. Then there is a natural map
$$
\diagram 
GW_p^{\LRt{n}}\LRf{{\SB}} \ar[r]^{u\cup -\qquad \quad} & GW_{p+1}^{\LRt{n+1}}\LRf{{\SC}h^b\LRf{{\BF}}\otimes {\SB}} \\
\enddiagram 
$$
In particular, there is a  natural map
$$
\diagram 
GW_p^{\LRt{n}}\LRf{{\BF}} \ar[r]^{u\cup -\quad } & GW_{p+1}^{\LRt{n+1}}\LRf{{\BF}} \\
\enddiagram 
$$
}
\eP 
\pf There is a natural map $\lambda: {\BF}^* \lra G{\CW}_1^{\LRt{1}}\LRf{{\BF}}$, by for example \cite[pp 1784]{S17}, the hyperbolic map
$K_1\LRf{{\BF}} \lra G{\CW}_1^{\LRt{1}}\LRf{{\BF}}$. ({\it Alternately, see} \cite[pp 67]{BL8}.)
Now, 
consider the diagram of maps:
$$
\diagram 
GW_p^{\LRt{n}}\LRf{{\SB}}\ar[d]_{\lambda(u)\otimes -} \ar@/^/@{-->}[rd]^{u \cup -}& \\
GW_1^{\LRt{1}}\LRf{{\SC}h^b\LRf{{\BF}}} \otimes GW_p^{\LRt{n}}\LRf{{\SB}}\ar[r] & GW_{p+1}^{\LRt{n+1}}\LRf{{\SC}h^b\LRf{{\BF}}\otimes {\SB}}\\
\enddiagram
$$
The horizontal map is given by tensor product \cite[pp 1782]{S17}. The vertical map sends $x\mapsto \lambda\LRf{u} 
\otimes x$. The required map    is defined as composition. 
\pic $\eop$
%%%%%%%%%%%%%%%%%%%%%%%%%%%%%%%%%%%%%%%%%%%%%%%%%
%%%%%%%%%%%%%%%%%%%%%%%%%%%%%%%%%%%%%%%%%%%%%%%%%
%%%%%%%%%%%%%%%%%%%%%%%%%%%%%%%%%%%%%%%%%%%%%%%%%

%%%%%%%%%%%%%%%%%%%%%%%%%%%%%%%%%%%%%%%%%%%%%%%%%
\subsection{Proof of the splitting lemma} %Naturality of the connecting morphism}
\label{secHomCom}
In this section, we prove that the diagram (\ref{alphaDia}) commutes. 
\bL\label{sliceIt}{\rm 
Let ${\BF}$ be a field, with $\frac{1}{2}\in {\BF}$. Let
$$
\LBrace{{ll}
A={\BF}\LRt{T} &\m=TA,~R=A_{\m} \\
Y=\spec{A} &X=\spec{R}\\
V=Y-\LRs{(T)} & U=X -\LRf{\m}=\spec{{\BF}\LRf{T}}\\
}
$$
In the following diagram, $-\cup T$ is the multiplication map (\ref{unitHike}), $\partial_T$ is the residue map,  $\d$ is the d\'{e}vissage, and $\delta_T=\d^{-1}\partial_T$. 
Then the diagram 
$$
\diagram 
G{\CW}^{\LRt{n}}_{p}\LRf{{\BF}(T)} \ar[rr]^{- \cup T} && G{\CW}^{\LRt{n+1}}_{p+1}\LRf{{\BF}(T)}
\ar@{-->}[d]^{\delta_T} \ar@/^/[dr]^{\partial_T}&\\
G{\CW}^{\LRt{n}}_{p}\LRf{{\BF}} \ar[u]^{\p}\ar[rr]_1&& G{\CW}^{\LRt{n}}_{p}\LRf{{\BF}} 
\ar[r]_{\d\qquad\qquad}^{\sim\qquad\qquad }&
G{\CW}^{\LRt{n}}_{p}\LRf{C{\BM}^{1}\LRf{R}}\\
\enddiagram 
$$
commutes. In particular, the pullback  map $\p$ is a split monomorphism. 
}
\eL 
\pf For $p\leq -2$ all the groups in the diagram are isomorphic to the Witt groups. For $p=-1$,  other than 
$G{\CW}_0^{\LRt{n+1}}\LRf{{\BF}\LRf{T}}$, all groups are isomorphic to the Witt groups. Further, if $n\neq 0~mod ~4$, the Witt groups 
${W}^{\LRt{n}}\LRf{-}=0$, for any field. Consequently, $\forall p\leq -1$ the  lemma follows from the 
classical  result \cite[Cor 5.1, pp 334]{M69} of Milnor, on Witt groups. So, we assume $p\geq 0$ and we work with the 
${\bf GW}$-spaces. 
%However, our methods work $G{\CW}$-spectra and do not distinguish the sign of $p$. So, we treat $p\in {\BZ}$. 
%
%
 We consider the diagram of  pullback (tensor product) maps:
$$
\diagram 
  GW^{\LRt{n}}_p\LRf{{\BF}} \ar@/_/[dr]\ar[r] & GW^{\LRt{n}}_p\LRf{R} \ar[d]\\ 
& GW^{\LRt{n}}_p\LRf{{\BF}\LRf{T}} \\  
  \enddiagram
$$
Fix $\LRt{\sigma}\in  GW_p^{\LRt{n}}\LRf{{\BF}}$.
Let 
$\tilde{\sigma}\in GW_p^{\LRt{n}}\LRf{{\BF}\LRf{T}}$ be the pullback of $\sigma$ and 
 $\Sigma\in GW_p^{\LRt{n}}\LRf{R}$ be the pullback of $\sigma$.
Let 
$$
\LBrace{{l}
\sigma:{\BS}^p \lra GW^{\LRt{n}}\LRf{{\BF}}, ~{\rm represent}~\LRt{\sigma}\in GW^{\LRt{n}}_p\LRf{{\BF}\LRf{T}}\\% 
T: {\BS}^1 \lra GW^{\LRt{1}}\LRf{{\BF}\LRf{T}}, ~{\rm correspond~to~}T. \quad~{\rm Denote} ~
\ell\LRf{T}:=\LRt{T}\in GW_1^{\LRt{1}}\LRf{{\BF}\LRf{T}}.  \\
\tilde{\sigma}\wedge \ell\LRf{T}: {\BS}^p\wedge {\BS}^1 \lra GW^{\LRt{n}}\LRf{{\BF}\LRf{T}}\wedge  GW^{\LRt{1}}\LRf{{\BF}\LRf{T}}\\
}
$$
%We consider a few diagrams. 
First consider the diagram of functors of dg categories: 
$$
\diagram 
&& {\bf dg}C{\BM}^1\LRf{R}\ar[d]^{cone}\\
&&  {\bf dg}^{+1}C{\BM}^0\LRf{R}\ar[d]^{\otimes}\\
{\bf dg}{\BF} \ar[r]_{\otimes\qquad \quad}\ar@/^/[urr]^{\Theta}\ar@/^/[rruu]^{\iota}&{\bf dg}C{\BM}^0\LRf{{\BF}\LRf{T}}&  {\bf dg}^{+1}C{\BM}^0\LRf{{\BF}\LRf{T}}\\
\enddiagram 
$$
The functor $\iota$ is the pushforward map, give by the map $\diagram R \ar[r] & {\BF}\\ \enddiagram$. The functor   $\Theta$ is obtained by the resolving complexes, and taking the double complex.  
In this case, we have a canonical choice
$$
\Theta\LRf{M_{\bul}}=Cone\LRf{\diagram RM_{\bul} \ar[r]^T & RM_{\bul}\\ \enddiagram}
$$
Applying $GW$ to the above diagram we obtain diagram of maps of spectra:
$$
\diagram 
&&  GW^{\LRt{n}}\LRf{C{\BM}^1\LRf{R}}\ar[d]\\
&&  GW^{\LRt{n+1}}\LRf{R}\ar[d]^{\pi}\\
GW^{\LRt{n}}\LRf{{\BF}}\ar[r]_{\p}\ar@/^/[urr]^{\Theta}\ar@/^/[rruu]^{\d}_{\sim}&
GW^{\LRt{n}}\LRf{{\BF}\LRf{T}}&  GW^{\LRt{n+1}}\LRf{{\BF}\LRf{T}}\\
\enddiagram  
$$
%This map is a quasi isomorphism. 
The upper diagonal map is homotopy equivalence, by d\'{e}vissage. 
In ${\bf dg}^{\LRt{1}}C{\BM}^0R$, consider the  symmetric space (quasi-isomorphism)
$$
\Phi:=
\LBrace{{l}
\diagram
\deg= &1&0\\
{C}^{\LRt{1}}_{\bul}\LRf{{R}}:=\ar[d]_{\varphi}& R \ar@{^(->}[r]^T\ar[d]_{-1}& R\ar[d]^{1}\\
{C}^{\LRt{1}}_{\bul}\LRf{{R}}^{\#_1}:=& R^{\vee} \ar@{^(->}[r]_{-T}& R^{\vee}\\
 \enddiagram  % 
 }
 $$
 We consider the homotopy fibration \cite[pp 424]{M23}
 $$
 \diagram 
 {\bf GW}^{\LRt{1}}\LRf{R} \ar[r] & {\BB}{\BQ}\LRf{{R}, h^{\LRt{1}}} \ar[r]^{\f} & {\BB}{\BQ}\LRf{{R}} \\
 \enddiagram 
 $$
 where $h^{\LRt{1}}$ denotes the duality in ${\bf dg}^{\LRt{1}}C{\BM}^0R$, and ${\BQ}$ denotes the Quillen categories. 
 The homotopy fiber of $\f$ is given by \cite[pp 600]{M23}
 $$
 {\SF}\LRf{\f} =\LRs{
 \LRf{x, \gamma}\in {\BB}{\BQ}\LRf{R, h^{\LRt{1}}}  \times {\bfP}\LRf{{\BB}{\BQ}\LRf{R}}: 
 \LBrace{{l}
 \gamma(0)=0\\
  \gamma(1)=\f\LRf{x}\\
  }.
 }
 $$
 where ${\bfP}$ denotes the path space. 
 In fact, ${\SF}\LRf{\f}=:{\bf GW}^{\LRt{1}}\LRf{R}$ is our model of the ${\bf GW}$-space. 
 There is an arrow $\P$ in the hermitian Quillen category ${\BQ}\LRf{R, h^{\LRt{1}}}$ and its image in
 the Quillen category  ${\BQ}\LRf{R}$, as follows
 $$
 \P:=
 \diagram 
 {\bf 0} \ar[d]  & {C}^{\LRt{1}}_{\bul}\LRf{{R}}\ar[d]_{\varphi}^{=\pm 1}\ar@{->>}[l] \ar@{^(->}[r]^1& {C}^{\LRt{1}}_{\bul}\LRf{{R}}\ar[d]^{\varphi}\\
 {\bf 0}\ar@{^(->}[r] & {C}^{\LRt{1}}_{\bul}\LRf{{R}}^{\#_1} & {C}^{\LRt{1}}_{\bul}\LRf{{R}}^{\#_1}\ar@{->>}[l]^1\\
 \enddiagram 
 \mapsto 
  \diagram 
 {\bf 0}   & {C}^{\LRt{1}}_{\bul}\LRf{{R}}\ar@{->>}[l] \ar@{=}[r]& {C}^{\LRt{1}}_{\bul}\LRf{{R}}\\
 \enddiagram 
 $$
 So, there are paths 
 $$
 \diagram
  I \times \P \ar[rr] && {\BB}{\BQ}\LRf{R, h} \ar[d]^{\f}\\
 I \ar[u]_{\wr}\ar[rru]^{\p}\ar[rr]_{\f\p}  &&  {\BB}{\BQ}\LRf{R}\\
 \enddiagram \quad {\rm where}\quad  I \times \P ~{\rm denotes~the}~1~{\rm simplex}.
 $$
 Note 
$$
\LBrace{{l}
\f\p(0)=\f\LRf{0\times \P}={\bf 0}\\
\f\p(1)=\f\LRf{1\times \P}=C_{\bul}\LRf{R}\\
\f\p(t)=\f\LRf{t\times \P}=t\times \f\LRf{\P}\\
}
$$
and 
$$
\LBrace{{ll}
\theta(s):=\f\p(st)=\f\LRf{st\times \P}=st\times \f\LRf{\P} &0\times \f\P \mapsto s\times \f\P\\
\theta(0)=\f\p(0t)=\f\LRf{0\times \P}=0\times \f\LRf{\P}={\bf 0} & {\rm the ~constant~path}\\
\theta(1):=\f\p(1t)=\f\LRf{t\times \P}=t\times \f\LRf{\P}&=\f\p\\
}
$$
 and $\f\p(1)=\f\LRf{\varphi}$ , $\LRf{\varphi, \f\p} \in {\SF}\LRf{\f}$. Further
 $$
\ell_R(s):= \LRf{\LRf{s, \P}, \theta(s)}\in {\SF}\LRf{\f}\qquad 
\LBrace{{l}
\ell_R(0)= \LRf{(0, \P), 0}\\
\ell_R(1)= \LRf{\LRf{1, \P}, \f\p}=\LRf{C_{\bul}\LRf{R}, \varphi}\\
}
 $$
It follows
$$
\ell_R\LRf{T}_{|t=1}=\LRs{\diagram 
R\ar@{^(->}[r]^T\ar[d]_{-1} & R\ar[d]^1\\
R^{\vee} \ar@{^(->}[r]_{-T} & R^{\vee}\\
\enddiagram}=\Theta\LRf{1:{\BF} \lra {\BF}^{\vee}}
$$
 Now $\ell_R\LRf{T}\otimes {\BF}\LRf{T}=\ell\LRf{T}$. This is because the process is same, while the 
 terminal point of the latter is isometric to zero.

 %%%%%%%%%%%%%%%%%%%%%%%%%%%%%%%%%%%%%%%%%%%%%%%%%
Consider the diagram with $I=\LRt{0, 1}$, while we replace ${\BS}^p$ by $I^p$.
$$
\diagram 
&{\bf GW}^{\LRt{n}}\LRf{{\BF}}\ar[rrr]^{\d}_{\sim}\ar[rrrd]^{\Theta}&& &{\bf GW}^{\LRt{n}}\LRf{C{\BM}^{1}\LRf{R}}\ar[d]^{res}\\
I^p\times 1\ar@/^/@{-->}[ru]^{\sigma}\ar@/_/[dr]
&I^p\times 0\ar[d]^{(t, 0)}\ar@{-->}[rr]^{base\qquad \qquad\qquad} 
&& {\bf GW}^{\LRt{n}}C{\BM}^0\LRf{R} \wedge  {\bf GW}^{\LRt{1}}C{\BM}^0\LRf{R}\ar[r] \ar[d]
&{\bf GW}^{\LRt{n+1}}\LRf{R}\ar[d]\\
&I^p\times I\ar[rr]_{\tilde{\sigma}\wedge \ell\LRf{T}\qquad\qquad} \ar[rru]^{\Sigma\wedge \ell_R\LRf{T}}
%\Sigma\wedge \ell_R\LRf{T}}
&&
{\bf GW}^{\LRt{n}}\LRf{{\BF}\LRf{T}} \wedge {\bf GW}^{\LRt{1}}\LRf{{\BF}\LRf{T}} \ar[r]&{\bf GW}^{\LRt{n+1}}\LRf{{\BF}\LRf{T}}\\
\enddiagram
$$
The diagram commutes. 
Consequently,
$$
\delta\LRf{-\cup T}\p\LRf{\LRt{\sigma}}=\LRt{\sigma}
$$
as required. \pic $\eop$

%%%%%%%%%%%%%%%%%%%%%%%%%%%%%%%%%%%%%%%%%%%%%%%%%
%%%%%%%%%%%%%%%%%%%%%%%%%%%%%%%%%%%%%%%%%%%%%%%%%
%%%%%%%%%%%%%%%%%%%%%%%%%%%%%%%%%%%%%%%%%%%%%%%%%
%%%%%%%%%%%%%%%%%%%%%%%%%%%%%%%%%%%%%%%%%%%%%%%%%
%%%%%%%%%%%%%%%%%%%%%%%%%%%%%%%%%%%%%%%%%%%%%%%%%
%%%%%%%%%%%%%%%%%%%%%%%%%%%%%%%%%%%%%%%%%%%%%%%%%

\subsection{D\'{e}vissage}  
\bL\label{regLocal}{\rm 
Let $\LRf{R, \m, \kappa}$ be a regular local ring, $X=\spec{R}$, $Z=V\LRf{\m}$. Assume $\frac{1}{2}\in R$.  Then 
$$
\diagram
{G}W^{\LRt{n}}\LRf{{\bf dg}{\SV}\LRf{\kappa}} \ar[r] & {G}W^{\LRt{n}}\LRf{{\bf dg}C{\BM}^Z\LRf{X}}\\
\enddiagram
\quad {\rm is~a~homotopy~equivalence~} \forall n\in {\BZ}.
$$
}
\eL
\pf It follows from \cite[Thm 5.3]{M25}. 
Note $C{\BM}^Z\LRf{X}$ is the category of finite length $R$-modules, with duality $M \mapsto {E}xt^d\LRf{M, R}$. 
Hence ${\SA}:=Ch^b\LRf{C{\BM}^Z\LRf{X}}$
 is an abelian category. Let ${\SB}:= Ch^b\LRf{{\SV}\LRf{\kappa}}$.
Consider the diagram
$$ 
\diagram 
M_{\bul}^0\ar@{^(->}[d] &&0\ar@{^(->}[d]\ar[r] & M_{n-1}\ar[r] & \cdots \ar[r] & M_0 \ar[r] & 0\\
M_{\bul}^1\ar@{^(->}[d] &0\ar[r]&N^{1}\ar@{^(->}[d] \ar[r] & M_{n-1}\ar[r] & \cdots \ar[r] & M_0 \ar[r] & 0\\
\cdots\ar@{^(->}[d]  &0\ar[r]&\cdots\ar@{^(->}[d]\ar[r] & M_{n-1}\ar[r] & \cdots \ar[r] & M_0 \ar[r] & 0\\
M_{\bul}^{r-1}\ar@{^(->}[d] &0\ar[r] &N^{r-1}\ar@{^(->}[d]\ar[r] & M_{n-1}\ar[r] & \cdots \ar[r] & M_0 \ar[r] & 0\\
M_{\bul}:  &0 \ar[r] & M_n \ar[r] & M_{n-1}\ar[r] & \cdots \ar[r] & M_0 \ar[r] & 0\\
\enddiagram 
$$
Here $M_{\bul}\in \Obj{{\SA}}$ and with $N^r=M_n$ is a Jordan decomposition of $M_n$, with $\frac{N^k}{N^{k-1}}\cong \kappa$. It follows $\frac{M^k_{\bul}}{N^{k-1}_{\bul}}\cong \kappa^{\LRt{n}}\in {\SB}$. Here 
$\kappa^{\LRt{k}}=\LBrace{{ll}\kappa & if ~k=n\\ 0 & otherwise\\}$.  
Inductively, any $M_{\bul}\in \Obj{{\SA}}$ has a decomposition, by objects of ${\SB}$. By a similar argument,
 ${\SA}$ is also noetherian. Therefore, \cite[Thm 5.3]{M25} applies. 
$\eop$ 

%%%%%%%%%%%%%%%%%%%%%%%%%%%%%%%%%%%%%%%%%%%%%%%%%
\vspace{3mm}
We mostly use D\'{e}vissage in the following form. 
\bL\label{finfpdLem}{\rm 
Let $\LRf{A, \m, \kappa}$ be a regular local ring, with $\dim A=d$. Assume $\frac{1}{2}\in R$. Then  
$$
GW_p^{\LRt{n}}\LRf{\kappa} \cong GW_p^{\LRt{n}}\LRf{C{\BM}^d(\spec{A})}\quad \forall p, n
$$
}
\eL
\pf Refer to D\'{e}vissage \cite[Thm 6.1]{M25}. $\eop$
%%%%%%%%%%%%%%%%%%%%%%%%%%%%%%%%%%%%%%%%%%%%%%%%%
%%%%%%%%%%%%%%%%%%%%%%%%%%%%%%%%%%%%%%%%%%%%%%%%%

\subsection{Coincidence of $GW$ spectra and bispectra} 
Let $X$ be a quasi projective  scheme over an affine scheme $\spec{A}$ and ${\CL}$ be a line bundle on $X$.
Assume $\frac{1}{2}\in A$. 
 In this section, we record a proof that, if $X$ is regular, then the map $\diagram G{\CW}^{\LRt{n}}\LRf{C{\BM}^k\LRf{X}} \ar[r] & {\BG}{W}^{\LRt{n}}\LRf{C{\BM}^k\LRf{X}} \\ \enddiagram$ is a homotopy equivalence. Consequently, the induced maps of  the Gersten-Witt
complexes
$$
\diagram
\bigoplus_{x\in X^{\LRf{p}}}G{\CW}^{\LRt{n+k-p}}_{k-p}\LRf{C{\BM}^p\LRf{X_x, {\CL}_x}}
\ar[r] &\bigoplus_{x\in X^{\LRf{p}}}{\BG}{W}^{\LRt{n+k-p}}_{k-p}\LRf{C{\BM}^p\LRf{X_x, {\CL}_x}}
\enddiagram
$$
are  isomorphisms. 
The following is from \cite[Thm 7]{S06}.
\bT\label{NoethAbeMarchu} {\rm 
Let ${\SA}$ be a small noetherian abelian category. Then 
$$
{\BK}_p\LRf{{\SA}}= 0 \qquad \forall p\leq -1
$$

}\eT

%%%%%%%%%%%%%%%%%%%%%%%%%%%%%%%%%%%%%%%%%%%%
\bP\label{negKgrzero}{\rm 
Let $X$ be a regular quasi projective scheme
over an affine scheme $\spec{A}$.
% and ${\CL}$ be a line bundle on $X$.Assume $\frac{1}{2}\in A$.. 
Then 
$$
{\BK}_p\LRf{C{\BM}^k\LRf{X}}=0 \qquad \forall p\leq -1
$$
Consequently, the map
$$
\diagram 
{\bfK}\LRf{C{\BM}^k\LRf{X}}\ar[r] & {\BK}\LRf{C{\BM}^k\LRf{X}}
\enddiagram
$$
is homotopy equivalence of spectra, $\forall k$. 
}
\eP
\pf The proposition holds for $k=0$. The exact sequence 
$$
\diagram 
{\BK}_0\LRf{C{\BM}^{k+1}\LRf{X}} \ar[r] & {\BK}_0\LRf{C{\BM}^{k}\LRf{X}}  \ar[r]^{\p\qquad } & \bigoplus_{x\in X^{\LRf{k}}} 
{\BK}_0\LRf{C{\BM}^{k}\LRf{X_x}} \ar[dll]\\
{\BK}_{-1}\LRf{C{\BM}^{k+1}\LRf{X}} \ar[r] & {\BK}_{-1}\LRf{C{\BM}^{k}\LRf{X}}  \ar[r] & \bigoplus_{x\in X^{\LRf{k}}} 
{\BK}_{-1}\LRf{C{\BM}^{k}\LRf{X_x}} \ar[dll]\\\
{\BK}_{-2}\LRf{C{\BM}^{k+1}\LRf{X}} \ar[r] & {\BK}_{-2}\LRf{C{\BM}^{k}\LRf{X}}  \ar[r] & \bigoplus_{x\in X^{\LRf{k}}} 
{\BK}_{-2}\LRf{C{\BM}^{k}\LRf{X_x}}\ar[dll]\\\
\cdots \ar[r] & \cdots\ar[r] & \cdots \\
\enddiagram
$$
is exact \cite[pp 414]{M23}.
By induction,
$$
{\BK}_{-p}\LRf{C{\BM}^{k}\LRf{X}} =0 \quad \forall p\geq 1
$$
Further, by (\ref{NoethAbeMarchu}) 
$$
 \bigoplus_{x\in X^{\LRf{k}}} 
{\BK}_{-p}\LRf{C{\BM}^{k}\LRf{X_x}}=0 \quad \forall p\geq 1
$$
Remains to prove that
$$
{\BK}_{-1}\LRf{C{\BM}^{k+1}\LRf{X}} =0
$$
However, the map $\p$ is surjective. To see this, fix $x_0\in X^{\LRf{k}}$ and let $M\in C{\BM}^k\LRf{X_{x_0}}$. Then there is a resolution 
$$
0 \subseteq M_1 \subseteq M_2 \subseteq M_3 \subseteq \cdots \subseteq M_n=M
\qquad \ni \quad \LRt{\frac{M_i}{M_{i-1}}}= \kappa\LRf{x_0}
$$
In other words, 
${\BK}_0\LRf{C{\BM}^{k}\LRf{X_{x_0}} }$ is generated by $\LRt{\kappa\LRf{x_0}}$. Let ${\wp}$ be the
homogeneous prime ideal, corresponding to $x_0$. Then 
$$
\p\LRf{\LRt{\frac{{\CO}_X}{\tilde{{\wp}}}}}=\kappa\LRf{x_0}
$$
\pic $\eop$

%%%%%%%%%%%%%%%%%%%%%%%%%%%%%%%%%%%%%%%%%%%%%%%%
\bC\label{consSpan}{\rm 
Let $X$ be a regular quasi projective scheme over an affine scheme $\spec{A}$ and $\frac{1}{2}\in A$. 
Let ${\CL}$ be a line bundle on $X$.
Then the map
$$
\diagram 
G{\CW}^{\LRt{n}}\LRf{C{\BM}^k\LRf{X, {\CL}}}  \ar[r] & {\BG}{W}^{\LRt{n}}\LRf{C{\BM}^k\LRf{X, {\CL}}}  \\
\enddiagram
$$
is a homotopy equivalence of bi spectrum.
}
\eC
\pf Fix $k, n$. By \cite[Prop 8.7]{S17}, \cite[pp 466]{M23}, we have 
$$
{G}{\CW}_p^{\LRt{n}}\LRf{C{\BM}^k\LRf{X, {\CL}}}= {\BG}{W}_p^{\LRt{n}}\LRf{C{\BM}^k\LRf{X, {\CL}}}\quad \forall p\geq 0
$$ 
%For negative groups the proof of \cite[Prop 9.3, 1805]{S17} works (\TCP{\it which is beyond the scope of my patience and concentration span.}) 
We write ${\SA}= {\bf dg}C{\BM}^k\LRf{X, {\CL}}$. By \cite[Prop 8.14, 1803]{S17} the diagram of the bi spectra
$$
\diagram
G{\CW}^{\LRt{n}}\LRf{{\SA}} \ar[r] \ar[d]& {\BG}{W}^{\LRt{n}}\LRf{{\SA}}\ar[d]\\
{\bf K}\LRf{{\SA}}^{hC_2} \ar[r] & {\BK}\LRf{{\SA}}^{hC_2} \\
\enddiagram 
$$
is homotopy cartesian. 
By Prop \ref{negKgrzero}, the map
$$
\diagram
{\bf K}\LRf{{\bf dg}C{\BM}^k\LRf{X}} \ar[r] & {\BK}\LRf{{\bf dg}C{\BM}^k\LRf{X}} \\
\enddiagram 
\quad {\rm is ~a~ stable~ equivalence}.
$$
Therefore the lower larrow of the above rectangle is s stable equivalence. Hence so is the upper arrow.
\pic $\eop$ 

%%%%%%%%%%%%%%%%%%%%%%%%%%%%%%%%%%%%%%%%%%%%%%%%
%%%%%%%%%%%%%%%%%%%%%%%%%%%%%%%%%%%%%%%%%%%%%%%%
%%%%%%%%%%%%%%%%%%%%%%%%%%%%%%%%%%%%%%%%%%%%%%%%

%\section{Parity of homotopy spectra}
%\bL{\rm 
%Let $\LRf{A, \m}$ be a regular local ring, with $\dim A=d$, and $\frac{1}{2}\in A$. Let $R=A[T]$ be a polynomial ring, and $Y=\spec{R}$, and $\tilde{\m}={\m}R$. Fix $y_0\in Y^{\LRf{d+1}}$ and $Y_0=Y_{y_0}$. By \cite[pp 473]{M23} we have a homotopy fibration of spectrum
%$$
%\diagram 
%{\BG}W^{\LRt{n}}\LRf{C{\BM}^{d+1}\LRf{Y_0}} \ar[r] & {\BG}W^{\LRt{n+1}}\LRf{C{\BM}^{d}\LRf{Y_0}} \ar[r] & 
%\coprod_{x\in Y_0^{\LRf{d}}}{\BG}W^{\LRt{n+1}}\LRf{C{\BM}^{d}\LRf{Y_0}_x} \\
%&&\TCP{\tilde{m}~not~alone}\ar@{..>}[u]\\
%\enddiagram 
%$$
%}
%\eL
%%%%%%%%%%%%%%%%%%%%%%%%%%%%%%%%%%%%%%%%%%%%%%%%%
%%%%%%%%%%%%%%%%%%%%%%%%%%%%%%%%%%%%%%%%%%%%%%%%%
%%%%%%%%%%%%%%%%%%%%%%%%%%%%%%%%%%%%%%%%%%%%%%%%%
%%%%%%%%%%%%%%%%%%%%%%%%%%%%%%%%%%%%%%%%%%%%%%%%%
%%%%%%%%%%%%%%%%%%%%%%%%%%%%%%%%%%%%%%%%%%%%%%%%%
%%%%%%%%%%%%%%%%%%%%%%%%%%%%%%%%%%%%%%%%%%%%%%%%%
\subsection{Duality and the normal bundle} %From Altman and Kleiman \cite{AK70}}
We start with the following from  \cite[Chap 1, Thm 4.5]{AK70}.
\bT\label{AK4570}{\rm 
Let $Y$ be a scheme, and $X\subseteq Y$ be a closed subscheme, defined by an ideal sheaf $J \subseteq {\CO}_Y$. Assume that $X$ is regularly embedded in $Y$ ({\it meaning, $J$ is locally complete intersection}) and $k=\codim(X, Y)$. Then $\forall ~{\SF}\in \qcoh(Y)$ there are natural isomorphisms, as follows:
$$
%\forall ~{\SF}\in \qcoh(Y)~\exists~a~natural~map\quad 
\diagram
{\SE}xt^k_Y\LRf{{\CO}_X, {\SF}} \ar[r]^{\sim\quad} & {\CH}om\LRf{\Lambda^k\frac{J}{J^2}, \frac{\SF}{J{\SF}}}\ar[r]^{\sim\quad~}
& {\CH}om\LRf{{\CO}_X, \frac{\SF}{J{\SF}}}\Lambda^k\frac{J}{J^2}^{-1}\\
\enddiagram
$$
}
\eT
%%%%%%%%%%%%%%%%%%%%%%%%%%%%%%%%%%%%%%%%%%%%

\vspace{3mm}
The following is an immediate consequence,
\bC\label{1316Thinht}{\rm 
Let $Y$ be a scheme, and $X\subseteq Y$ be a closed subscheme, defined by an ideal sheaf ${\CJ} \subseteq {\CO}_Y$. Assume that $X$ is regularly embedded in $Y$ ({\it meaning, ${\CJ}$ is locally complete intersection}) and 
$k=\codim(X, Y)$. 
 Let $M\in {\SV}\LRf{X}$.  Then $\forall ~{\SF}\in \qcoh(Y)$ there are natural isomorphisms, as follows:
$$
%\forall ~{\SF}\in \qcoh(Y)~\exists~a~natural~map\quad 
\diagram
{\SE}xt^k\LRf{M, {\SF}} \ar[r]^{\sim\qquad} \ar@/_/[dr]_{\Phi}& {\CH}om\LRf{M\Lambda^k\frac{{\CJ}}{{\CJ}^2}, \frac{\SF}{{\CJ}{\SF}}}\ar[d]^{\wr}\\
& {\CH}om\LRf{M, \frac{\SF}{{\CJ}{\SF}}}\Lambda^k\frac{{\CJ}}{{\CJ}^2}^{-1}\\
\enddiagram
$$
}
\eC
\pf Note that the vertical isomorphism follows from further generalities \cite[pp 235]{H77}. So, we shall establish the horizontal isomorphism only. 
First, assume $Y=\spec{A}$ is affine, and $\overline{A}=\frac{A}{J}$. Assume $M=\overline{A}^n$ is free. Then the assertion follows immediately from (\ref{AK4570}). 

Now suppose $Y$ is any noetherian scheme. Let $U=\spec{A}\subseteq Y$ be an open subset, such that 
$J=\Gamma\LRf{U, {\CJ}}\subseteq A$ is a complete intersection ideal, and $M_U:=\Gamma\LRf{U\cap X, M}$ is free.
Write $F_U=\Gamma\LRf{U, {\SF}}$. From the affine case, there are natural isomorphisms
\begin{equation}\label{proffendia}
\diagram
{E}xt^k\LRf{M_U, {\SF}_{|U}} \ar[r]^{\psi\qquad}_{\sim\qquad} \ar@/_/[dr]_{\Phi_U}& {H}om\LRf{M_U\Lambda^k\frac{{J}}{{J}^2}, \frac{F_U}{{J}{F}_U}}\ar[d]^{\wr}\\
& {H}om\LRf{M, \frac{F_U}{{J}{F_U}}}\Lambda^k\frac{{J}}{{J}^2}^{-1}\\
\enddiagram
%\quad These~are~direct~sum~of~maps. 
\end{equation}
The map $\psi$ is independent of basis. To see this write $\overline{A}= \frac{A}{J}$ and $M_U=\bigoplus_{i=1}^n\overline{A}e_i= \bigoplus_{i=1}^n\overline{A}\varepsilon_i$, where $e_1, \ldots, e_n$ and $\varepsilon_1, \ldots, 
\varepsilon_n$ are   two sets of  bases of $M_U$. 
Let $\hat{e}_1, \ldots, \hat{e}_n$ and $\hat{\varepsilon}_1, \ldots, 
\hat{\varepsilon}_n$ be the dual basis in
 $Hom\LRf{M_U\Lambda^k\frac{J}{J^2}, \frac{F_U}{JF_U}}$.
Then there are direct sum decompositions 
$$
\LBrace{{l}
{E}xt^k\LRf{M_U, {\SF}_{|U}}=\bigoplus_{i=1}^n{E}xt^k\LRf{\overline{A}, {\SF}_{|U}}e_i\\
{E}xt^k\LRf{M_U, {\SF}_{|U}}=\bigoplus_{i=1}^n{E}xt^k\LRf{\overline{A}, {\SF}_{|U}}\varepsilon_i\\
%}
%\LBrace{{l}
{H}om\LRf{M_U\Lambda^k\frac{{J}}{{J}^2}, \frac{F_U}{{J}{F}_U}}=\bigoplus_{i=1}^n {H}om\LRf{\Lambda^k\frac{{J}}{{J}^2}, \frac{F_U}{{J}{F}_U}}\hat{e}_i\\
{H}om\LRf{M_U\Lambda^k\frac{{J}}{{J}^2}, \frac{F_U}{{J}{F}_U}}=\bigoplus_{i=1}^n {H}om\LRf{\Lambda^k\frac{{J}}{{J}^2}, \frac{F_U}{{J}{F}_U}}\hat{\varepsilon}_i\\
}
$$
We can write 
$$
\matrix{{ccc} e_1 & \ldots & e_n}= \matrix{{ccc} \varepsilon_1 & \ldots & \varepsilon_n}T
\quad where\quad T=\LRf{\t_{ij}}\in GL_n\LRf{\overline{A}}
$$ 
We borrow the notation $\varphi$ from \cite[Thm I.4.5]{AK70}. 
Define two maps (isomorphisms)
$$
\diagram
{E}xt^k\LRf{M_U, {\SF}_{|U}} \ar[r]^{\psi_1\quad}_{ \psi_2\quad}& 
{H}om\LRf{M_U\Lambda^k\frac{{J}}{{J}^2}, \frac{F_U}{{J}{F}_U}}\\
%{E}xt^k\LRf{M_U, {\SF}_{|U}} \ar@/_/[ru]_{\psi_2} & \\
%{H}om\LRf{M_U\Lambda^k\frac{{J}}{{J}^2}, \frac{F_U}{{J}{F}_U}}\\
\enddiagram
$$
as follows
$$
\LBrace{{l}
\psi_1\LRf{e_i}=\varphi_{|\overline{A}e_i}\\
\psi_2\LRf{\varepsilon_i}=\varphi_{|\overline{A}e\varepsilon_i}\\
}
\diagram 
\overline{A}e_i \ar[r]^{\sim}\ar@/_/[dr]^{\psi_1} & \overline{A}\ar[d]^{\varphi}\\
& {H}om\LRf{\Lambda^k\frac{{J}}{{J}^2}, \frac{F_U}{{J}{F}_U}}\hat{e}_i\\
\enddiagram
\quad
\diagram 
\overline{A}\varepsilon_i \ar[r]^{\sim}\ar@/_/[dr]^{\psi_2} & \overline{A}\ar[d]^{\varphi}\\
& {H}om\LRf{\Lambda^k\frac{{J}}{{J}^2}, \frac{F_U}{{J}{F}_U}}\hat{\varepsilon}_i\\
\enddiagram
$$
So, ${\bf x}\in {E}xt^k\LRf{M_U, {\SF}_{|U}}$ can be written as
$$
{\bf x}=\matrix{{ccc} e_1 & \ldots & e_n}\matrix{{c}\a_1\\ \cdots \\ \a_n\\}
= \matrix{{ccc} \varepsilon_1 & \ldots & \varepsilon_n}T\matrix{{c}\a_1\\ \cdots \\ \a_n\\}
\quad \a_i\in {E}xt^k\LRf{\overline{A}, {\SF}_{|U}}
$$
So,
$$
\psi_1\LRf{{\bf x}} =
\matrix{{ccc} \hat{e}_1 & \ldots & \hat{e}_n}\matrix{{c}\Phi\LRf{\a_1}\\ \cdots \\ \Phi\LRf{\a_n}\\}
\quad 
\psi_2\LRf{{\bf x}} =
\matrix{{ccc} \hat{\varepsilon}_1 & \ldots & \hat{\varepsilon}_n}T\matrix{{c}\Phi\LRf{\a_1}\\ \cdots \\ \Phi\LRf{\a_n}\\}
$$
However, 
$$
\matrix{{ccc} \hat{e}_1 & \ldots & \hat{e}_n}= \matrix{{ccc} \hat{\varepsilon}_1 & \ldots & \hat{\varepsilon}_n}T
$$
Therefore,
$$
\psi_1\LRf{{\bf x}} =\psi_2\LRf{{\bf x}} 
$$
Therefore, it is established that the map $\psi$ (\ref{proffendia}) is independent of the bases. 
 So, $\psi$ extends to a global isomorphism on $Y$. $\eop$ 
 
%%%%%%%%%%%%%%%%%%%%%%%%%%%%%%%%%%%%%%%%%%%%
\bC\label{subsech} %{ZXYcontexdt}
{\rm 
Let $Y$ be a Cohen Macaulay scheme. Let   $Z\subseteq X\subseteq Y$ be  closed subschemes of $Y$ with 
$\codim(X, Y)=r$ and $\codim(Z, Y)=r+s$.  Assume      $Z\subseteq X$ and $X\subseteq Y$ are regularly embedded. Let ${\SI}_X\subseteq {\SI}_Z\subseteq {\CO}_Y$ be the ideal sheaves of $X$ and $Z$, respectively.
Denote
$$
\LBrace{{ll}
{\CO}_Z= \frac{{\CO}_Y}{{\SI}_Z} &N_{Y/Z}=\Lambda^{r+s}\frac{{\SI}_Z}{{\SI}_Z^2} \\
{\CO}_X= \frac{{\CO}_Y}{{\SI}_X} &N_{Y/X}=\Lambda^r\frac{{\SI}_X}{{\SI}_X^2} \\
%\overline{{\SJ}}
{\SI}=\frac{{\SJ}_Z}{{\SI}_X} &N_{X/Z}=\Lambda^s\frac{{\SI}}{{\SI}^2}=\Lambda^s\frac{{\SI}_Z}{{\SI}_Z^2+{\SI}_X} \\
}
$$
Let ${\CL}$ be an invertible sheaf on $Y$.
Then the following  %We have a commutative diagram of functors:
$$
\diagram 
& {\SV}\LRf{Z, {\CL}N_{Y/Z}^{-1}}\ar@/^/[dr]^{\varphi_Y} \ar@/_/[dl]_{\varphi_X} &C{\BM}^Z\LRf{Z, {\SE}xt^0\LRf{-, {\CL}N_{Y/Z}^{-1}}} \ar@{..>}[l]\\
C{\BM}^Z\LRf{X, {\SE}xt^s_X\LRf{-, N_{Y/X}^{-1}{\CL}}} \ar[rr]_{\varphi}
&& C{\BM}^{Z}\LRf{Y, {\SE}xt^{r+s}_Y\LRf{-, {\CL}}} \\
\enddiagram 
$$
is a commutative diagram of duality preserving functors. 
}
\eC
\pf At the object level, the functor is defined by $\varphi(M)=M$.  We need prove that
\begin{equation}\label{toPrkl}
\varphi\LRf{M^{\vee_X}}={\SE}xt^s_X\LRf{M, N_{Y/X}^{-1}{\CL}}\cong {\SE}xt^{r+s}_Y\LRf{M, {\CL}}\\
\end{equation}
is an isomorphism. 
Note $s=\codim\LRf{Z, X}$. Let $M\in C{\BM}^Z(X)$. Consider the following diagram of resolutions.
$$
\diagram 
& 0 \ar[d] &0 \ar[d]  &  & 0 \ar[d] &0\ar[d] &  &  \\
0 \ar[r] & F_{sr} \ar[r]\ar[d] & F_{s-1,r} \ar[r]\ar[d]  & \cdots \ar[r] & F_{1r} \ar[r] \ar[d] & F_{0r} \ar[r] \ar[d] & 0 &  F_{ir}\in {\SV}(Y)\\
0 \ar[r] & F_{s,r-1} \ar[r]\ar[d] & F_{s-1,r-1} \ar[r]\ar[d]  & \cdots \ar[r] & F_{1r-1} \ar[r] \ar[d] & F_{0r-1} \ar[r] \ar[d] & 0 &  F_{i, r-1}\in {\SV}(Y)\\
0\ar[r]& \cdots \ar[d] \ar[r]&  \cdots \ar[r]\ar[d]  & \cdots \ar[r] &  \cdots \ar[r] \ar[d] &  \cdots\ar[r] \ar[d] & 0  & \\
0 \ar[r] & F_{s1} \ar[r]\ar[d] & F_{s-1,1} \ar[r]\ar[d]  & \cdots \ar[r] & F_{11} \ar[r] \ar[d] & F_{01} \ar[r] \ar[d] & 0  & F_{i1}\in {\SV}(Y)\\
0 \ar[r] & F_{s0} \ar[r]\ar[d] & F_{s-1,0} \ar[r]\ar[d]  & \cdots \ar[r] & F_{10} \ar[r] \ar[d] & F_{00} \ar[r] \ar[d] & 0  & F_{i0}\in {\SV}(Y)\\
0 \ar[r] & P_s \ar[r] & P_{s-1} \ar[r] & \cdots \ar[r] & P_1 \ar[r] & P_0 \ar[r] & M \ar[r] &0  P_i\in {\SV}(X)\\
\enddiagram 
$$ 
The double complex is a resolution of $M$ in ${\SV}(Y)$. 
Now, we have  the following commutative diagram %This leads to 
{\scalefont{.7}
\begin{equation}\label{1394Hola}
\diagram 
0 \ar[r] & {\SE}xt^r\LRf{P_0, {\CL}}\ar[r] \ar[d]_{\wr} &  {\SE}xt^r\LRf{P_1, {\CL}}\ar[r] \ar[d]_{\wr} & \cdots \ar[r] &  {\SE}xt^r\LRf{P_{s-1}, {\CL}}\ar[r] \ar[d]_{\wr} & {\SE}xt^r\LRf{P_s, {\CL}}\ar[r]\ar[d]_{\wr}  & 0\\
0 \ar[r] & {\CH}om\LRf{P_0, N_{Y/X}^{-1}{\CL}}\ar[r] &  {\CH}om\LRf{P_1, N_{Y/X}^{-1}{\CL}}\ar[r] & \cdots \ar[r] &  
 {\CH}om\LRf{P_{s-1}, N_{Y/X}^{-1}{\CL}}\ar[r] & {\CH}om\LRf{P_s, N_{Y/X}^{-1}{\CL}}\ar[r] & 0\\
\enddiagram 
\end{equation}
}
Here the vertical maps are isomorphisms by  (\ref{1316Thinht}). Now dualize the double complex, by ${\CH}om\LRf{-, {\CL}}$. 
Since $P_i$ is a resolution of $M$ in ${\SV}(X)$, from the $2^{nd}$ line of (\ref{1394Hola}), we obtain  the exact sequence 
$$
\diagram 
0 \ar[r] & {\CH}om\LRf{P_0, N_{Y/X}^{-1}{\CL}}\ar[r] %&  {\CH}om\LRf{P_1, N_{Y/X}^{-1}{\CL}}\ar[r] 
& \cdots \ar[r] &  
 %{\CH}om\LRf{P_{s-1}, N_{Y/X}^{-1}{\CL}}\ar[r] &
  {\CH}om\LRf{P_s, N_{Y/X}^{-1}{\CL}}\ar[r] & {\SE}xt^{r}_X\LRf{M, N_{Y/X}^{-1}{\CL}}\ar[r] & 0\\
\enddiagram 
$$
Dualizing the full double complex by ${\CH}om\LRf{-, {\CL}}$, and considering the total complex, we obtain
$$
{\SE}xt^{r+s}_Y\LRf{M, {\CL}}\cong {\SE}xt^{r}_X\LRf{M, N_{Y/X}^{-1}{\CL}}
$$
This establishes (\ref{toPrkl}). 
\bE
\item Definition and duality issues of $\varphi_Y$ is a particular case of that of $\varphi$. 
\item By the same argument there is natural functor 
$$
\diagram 
{\SV}\LRf{Z, N_{Y/X}^{-1}{\CL} N_{X/Z}^{-1}} \ar[r] & C{\BM}^Z\LRf{X, {\SE}xt^s_X\LRf{-, N_{Y/X}^{-1}{\CL}}} \\
\enddiagram 
$$
However, it is easy to check
$$
  N_{X/Z}N_{Y/X}\cong \Lambda^s\frac{{\SI}_Z}{{\SI}_Z^2+{\SI}_X}\otimes \Lambda^r \frac{{\SI}_X}{{\SI}_X^2}
  \cong \Lambda^{r+s} \frac{{\SI}_Z}{{\SI}_Z^2}
  \cong N_{Y/Z}
$$
\eE 
\pic $\eop$

\end{document}